\documentclass[12pt,frenchb,english]{article}
\usepackage[T1]{fontenc}
\usepackage{babel}
\usepackage{graphics,hyperref}

\makeatletter

\providecommand{\LyX}{L\kern-.1667em\lower.25em\hbox{Y}\kern-.125emX\@}

\NoAutoSpaceBeforeFDP

\usepackage{amsthm}
\theoremstyle{plain}    
\newtheorem{thm}{Th{\'e}or{\`e}me}[section]
\newenvironment{lyxlist}[1]
  {\begin{list}{}
    {\settowidth{\labelwidth}{#1}
     \setlength{\leftmargin}{\labelwidth}
     \addtolength{\leftmargin}{\labelsep}
     }}
  {\end{list}}
\theoremstyle{plain}    
\newtheorem{cor}[thm]{Corollaire} 
\theoremstyle{plain}    
\newtheorem{lem}[thm]{Lemme} 
\theoremstyle{definition}
 \newtheorem{problem}[thm]{Probl{\`e}me}
\theoremstyle{remark}
\newtheorem{rem}[thm]{Remarque}
\theoremstyle{remark}    
\newtheorem*{claim*}{Affirmation}

\@addtoreset{equation}{section}

\newcommand{\Diff}{\mathrm{Diff}}
\renewcommand{\Re}{\mathrm{Re}\,}
\renewcommand{\Im}{\mathrm{Im}\,}
\renewcommand{\hom}{\mathrm{Hom}}

\usepackage{amssymb,eucal}

\NoAutoSpaceBeforeFDP
\makeatother

\begin{document}

\selectlanguage{frenchb}

\title{M{\'e}triques autoduales sur la boule}

\author{Olivier Biquard}

\maketitle
\selectlanguage{english}

\begin{abstract}
\renewcommand{\thefootnote}{}\footnote{
2000 Mathematics Subject Classification: 53C25, 32G07, 53C28, 53C26
}A conformal metric on a 4-ball induces on the boundary 3-sphere a conformal
metric and a trace-free second fundamental form. Conversely, such a data on
the 3-sphere is the boundary of a unique selfdual conformal metric, defined
in a neighborhood of the sphere. In this paper we characterize the conformal
metrics and trace-free second fundamental forms on the 3-sphere (close to the
standard round metric) which are boundaries of selfdual conformal metrics on
the whole 4-ball.

When the data on the boundary is reduced to a conformal metric (the trace-free
part of the second fundamental form vanishes), one may hope to find in the conformal
class of the filling metric an Einstein metric, with a pole of order 2 on the
boundary. We determine which conformal metrics on the 3-sphere are boundaries
of such selfdual Einstein metrics on the 4-ball. In particular, this implies
the Positive Frequency Conjecture of LeBrun.

The proof uses twistor theory, which enables to translate the problem in terms
of complex analysis; this leads us to prove a criterion for certain integrable
CR structures of signature (1,1) to be fillable by a complex domain.

Finally, we solve an analogous, higher dimensional problem: selfdual Einstein
metrics are replaced by quaternionic-K{\"a}hler metrics, and conformal structures
on the boundary by quaternionic contact structures (previously introduced by
the author); in contrast with the 4-dimensional case, we prove that any small
deformation of the standard quaternionic contact structure on the \( (4m-1) \)-sphere
is the boundary of a quaternionic-K{\"a}hler metric on the \( (4m) \)-ball.
\end{abstract}
\selectlanguage{frenchb}\renewcommand{\thefootnote}{\fnsymbol{footnote}}\setcounter{footnote}{0}Cet
article est consacr{\'e} {\`a} la construction de m{\'e}triques riemanniennes autoduales.
En dimension 4, le tenseur de Weyl d'une m{\'e}trique riemannienne \( g \) se d{\'e}compose
en composantes autoduale et antiautoduale, \( W=W_{+}+W_{-} \), et la m{\'e}trique
est dite autoduale si \( W_{-}=0 \) ; cette propri{\'e}t{\'e} ne d{\'e}pend d'ailleurs
que de la classe conforme \( [g] \) de la m{\'e}trique.

Une m{\'e}trique conforme sur une vari{\'e}t{\'e} {\`a} bord induit sur le bord la donn{\'e}e conforme
\( [g,\mathcal{Q}] \), o{\`u} \( g \) est la m{\'e}trique et \( \mathcal{Q} \) la
partie {\`a} trace nulle de la seconde forme fondamentale. Dans cet article, on
r{\'e}sout, au voisinage de la m{\'e}trique ronde, le probl{\`e}me suivant :

\begin{problem}
\label{pro-1}Quelles donn{\'e}es \( [g,\mathcal{Q}] \) sur la sph{\`e}re \( S^{3} \)
sont les bords de m{\'e}triques autoduales \( [g] \) sur la boule \( B^{4} \)
?
\end{problem}
LeBrun \cite{LeB:82,LeB:85} a montr{\'e} que toute donn{\'e}e de ce type est \emph{localement}
le bord d'une unique m{\'e}trique autoduale ; de plus, dans le cas o{\`u} \( \mathcal{Q}=0 \),
la m{\'e}trique peut {\^e}tre prise d'Einstein, avec un p{\^o}le d'ordre 2 sur le bord :
on dit alors que la m{\'e}trique conforme sur le bord est l'\emph{infini conforme}
de la m{\'e}trique qui remplit (ainsi, la m{\'e}trique ronde sur \( S^{3} \) est l'infini
conforme de la m{\'e}trique hyperbolique r{\'e}elle sur \( B^{4} \)). Dans le cas d'une
m{\'e}trique invariante {\`a} gauche sur \( S^{3} \), une m{\'e}trique autoduale d'Einstein
explicite, globale, est connue : Pedersen \cite{Ped:86} a trait{\'e} le cas des
m{\'e}triques de Berger, et Hitchin \cite{Hit:95} le cas g{\'e}n{\'e}ral. Enfin, LeBrun
\cite{LeB:91} a montr{\'e} qu'il existe une famille de dimension infinie de m{\'e}triques
sur \( S^{3} \) munies de remplissages par des m{\'e}triques autoduales d'Einstein
sur \( B^{4} \).

Cependant, il existe des obstructions {\`a} un tel remplissage : Hitchin \cite{Hit:97}
a ainsi montr{\'e} que l'invariant \( \eta  \) d'une m{\'e}trique riemannienne sur
\( S^{3} \) doit {\^e}tre n{\'e}gatif pour qu'existe un remplissage autodual d'Einstein.

Gr{\^a}ce {\`a} la construction twistorielle de Penrose, la construction de m{\'e}triques
autoduales se ram{\`e}ne {\`a} la construction de certaines vari{\'e}t{\'e}s complexes de dimension
3. Une version de cette construction twistorielle, due {\`a} LeBrun \cite{LeB:84, LeB:85},
existe en dimension 3 pour les donn{\'e}es \( [g,\mathcal{Q}] \) : l'espace des
twisteurs est alors une vari{\'e}t{\'e} CR (int{\'e}grable) de dimension 5. Par exemple,
l'espace des twisteurs de \( S^{3} \) est l'hypersurface r{\'e}elle \( |z^{1}|^{2}+|z^{2}|^{2}=|z^{3}|^{2}+|z^{4}|^{2} \)
de l'espace projectif complexe \( P^{3} \), qui borde le domaine complexe 
\[
D=\{|z^{1}|^{2}+|z^{2}|^{2}\leq |z^{3}|^{2}+|z^{4}|^{2}\}\subset P^{3},\]
lui-m{\^e}me espace des twisteurs de \( B^{4} \). De cette mani{\`e}re, la question
suivante g{\'e}n{\'e}ralise le probl{\`e}me \ref{pro-1}.

\begin{problem}
\label{pro-2}Quelles structures CR sur \( \partial D \) sont les bords de
d{\'e}formations complexes de \( D \) ?
\end{problem}
Cette question naturelle de g{\'e}om{\'e}trie complexe est difficile, car \( \partial D \)
est de signature (1,1), un cas o{\`u} le th{\'e}or{\`e}me d'extension de Kiremidjian \cite{Kir:79}
ne s'applique pas. Dans le cas pseudoconvexe, en dimension 3, par exemple pour
la sph{\`e}re \( S^{3} \), le probl{\`e}me a {\'e}t{\'e} beaucoup {\'e}tudi{\'e}, et largement r{\'e}solu,
dans une succession d'articles, notamment par Burns-Epstein \cite{Bur-Eps:90},
Lempert \cite{Lem:92}, Epstein \cite{Eps:92}, et Bland \cite{Bla:94}. Nous
donnons, pour les structures CR int{\'e}grables proches de la structure standard
sur \( \partial D \), une condition n{\'e}cessaire et suffisante, analogue {\`a} celle
trouv{\'e}e par ces auteurs ; elle est donn{\'e}e en termes des coefficients de Fourier
par rapport {\`a} une action de \( S^{1} \).

\begin{thm}
\label{th-A}Une structure CR sur \( \partial D \), proche de la structure
standard, est le bord d'une d{\'e}formation complexe de \( D \) si et seulement
s'il existe une action de contact de \( S^{1} \) sur \( \partial D \), par
rapport {\`a} laquelle la structure complexe s'{\'e}crit seulement avec des coefficients
de Fourier positifs.
\end{thm}
Pour plus de pr{\'e}cisions, voir le lemme \ref{lem-fill-J} et le corollaire \ref{cor-ex-disques}.
Ajoutons qu'en r{\'e}alit{\'e}, on d{\'e}montre cet {\'e}nonc{\'e} en toute dimension, pour le bord
d'un domaine qui est l'espace total d'un fibr{\'e} en disques holomorphes, {\`a} forme
de Levi non d{\'e}g{\'e}n{\'e}r{\'e}e. Dans notre situation, cela permet de d{\'e}crire l'espace
tangent aux structures CR remplissables parmi toutes les structures CR int{\'e}grables
(th{\'e}or{\`e}me \ref{th-CR-rem}) ; toutes ne le sont pas ; en particulier, il existe
une famille de dimension infinie de structures CR non remplissables (car la
sous-vari{\'e}t{\'e} des structures CR remplissables est de codimension infinie), voir
le corollaire \ref{cor-CR-non-rem}.

On montre aussi (corollaire \ref{cor-CR-non-def}) qu'existe sur \( \partial D \)
une famille de dimension infinie de structures CR remplissables qui ne sont
pas induites par une d{\'e}formation de \( \partial D \) comme hypersurface de
\( P^{3} \) : il s'agit d'un ph{\'e}nom{\`e}ne nouveau, qui est d{\^u} {\`a} la pr{\'e}sence d'une
valeur propre n{\'e}gative dans la forme de Levi, et tranche avec le th{\'e}or{\`e}me de
Hamilton \cite{Ham:77}, ou le th{\'e}or{\`e}me de stabilit{\'e} de Lempert \cite{Lem:92}
sur \( S^{3}\subset \mathbf{C}^{2} \).

Les espaces de twisteurs des m{\'e}triques autoduales ont des propri{\'e}t{\'e}s suppl{\'e}mentaires
: structure r{\'e}elle, existence d'une famille de courbes rationnelles particuli{\`e}res,
et structure de contact holomorphe si la m{\'e}trique est d'Einstein. Cependant,
on montre que la seule obstruction au remplissage autodual d'une donn{\'e}e \( [g,\mathcal{Q}] \)
est le remplissage de la structure CR.

\begin{thm}
\label{th-B}La donn{\'e}e \( [g,\mathcal{Q}] \) d'une m{\'e}trique et d'une seconde
forme fondamentale conformes sur \( S^{3} \), proche de la structure standard,
est le bord d'une m{\'e}trique autoduale sur \( B^{4} \) si et seulement si l'espace
des twisteurs de \( [g,\mathcal{Q}] \) satisfait les conditions du th{\'e}or{\`e}me
\ref{th-A}.

Une m{\'e}trique conforme \( [g] \) sur \( S^{3} \) est l'infini conforme d'une
m{\'e}trique autoduale d'Einstein sur \( B^{4} \) si et seulement si l'espace des
twisteurs de \( [g] \) satisfait les conditions du th{\'e}or{\`e}me \ref{th-A}.
\end{thm}
Ce th{\'e}or{\`e}me donne une condition n{\'e}cessaire et suffisante de remplissage : cette
condition n'est pas facile {\`a} tester ; en particulier, on souhaiterait savoir
s'il existe, en dehors des exemples connus, des \( [g,\mathcal{Q}] \) qui soient
remplissables par une m{\'e}trique autoduale sur \( B^{4} \). Ce probl{\`e}me est r{\'e}solu
en appliquant le th{\'e}or{\`e}me \ref{th-B} {\`a} une description (th{\'e}or{\`e}me \ref{th-gQ-rem})
de la vari{\'e}t{\'e} \( \mathcal{B}_{+}/\Diff _{0}S^{3} \), o{\`u} \( \mathcal{B}_{+} \)
d{\'e}signe les \( [g,\mathcal{Q}] \) remplissables (on ignore dans cette introduction
les probl{\`e}mes de r{\'e}gularit{\'e}) et \( \Diff _{0}S^{3} \) d{\'e}signe le groupe des
diff{\'e}omorphismes de \( S^{3} \) qui fixent deux points dans le fibr{\'e} tangent
projectivis{\'e}, de sorte que 
\[
\Diff S^{3}/\Diff _{0}S^{3}=SO_{1,4},\]
 le groupe des transformations conformes de \( S^{3} \) ; on obtient notamment
une description explicite de l'espace tangent \( T(\mathcal{B}_{+}/\Diff _{0}S^{3}) \)
en la m{\'e}trique ronde. 

En particulier, notons \( \mathcal{B} \) l'espace de toutes les donn{\'e}es conformes
\( [g,\mathcal{Q}] \), et \( \mathcal{B}_{-} \) l'espace des \( [g,\mathcal{Q}] \)
remplissables par une m{\'e}trique antiautoduale, alors on d{\'e}duit du th{\'e}or{\`e}me \ref{th-gQ-rem}
le  r{\'e}sultat suivant.

\begin{thm}
\label{th-C0}Pr{\`e}s de la m{\'e}trique ronde, \( \mathcal{B}_{+}/\Diff _{0}S^{3} \)
et \( \mathcal{B}_{-}/\Diff _{0}S^{3} \) sont deux sous-vari{\'e}t{\'e}s transverses
de \( \mathcal{B}/\Diff _{0}S^{3} \), dont l'intersection est constitu{\'e}e des
\( [g,\mathcal{Q}] \) induits par une d{\'e}formation de la sph{\`e}re \( S^{3} \)
dans l'espace plat \( \mathbf{R}^{4} \).
\end{thm}
Remarquons que le th{\'e}or{\`e}me affirme en particulier que pour \( [g,\mathcal{Q}] \)
proche de la m{\'e}trique ronde dans l'intersection \( \mathcal{B}_{+}\cap \mathcal{B}_{-} \),
la m{\'e}trique autoduale et la m{\'e}trique antiautoduale qui remplissent co{\"\i}ncident,
ce qui n'{\'e}tait pas {\'e}vident a priori.

Restreignons-nous {\`a} pr{\'e}sent au cas d'une m{\'e}trique conforme sur le bord : notons
\( \mathcal{M} \) l'espace des m{\'e}triques conformes sur \( S^{3} \) (en oubliant
encore les questions de r{\'e}gularit{\'e}). Rappelons que, dans ce cas, le th{\'e}or{\`e}me
de Graham-Lee \cite{Gra-Lee:91} fournit sur \( B^{4} \) un unique remplissage
d'Einstein ; notons \( \mathcal{M}_{+} \) (resp. \( \mathcal{M}_{-} \)) l'espace
des m{\'e}triques conformes sur \( S^{3} \) dont le remplissage d'Einstein est
autodual (resp. antiautodual). Le th{\'e}or{\`e}me pr{\'e}c{\'e}dent s'applique aux m{\'e}triques
de la mani{\`e}re suivante (voir le th{\'e}or{\`e}me \ref{th-g-rem} pour une description
explicite des espaces tangents).

\begin{thm}
\label{th-C}Dans \( \mathcal{M}/\Diff _{0}S^{3} \), les espaces \( \mathcal{M}_{+}/\Diff _{0}S^{3} \)
et \( \mathcal{M}_{-}/\Diff _{0}S^{3} \) de m{\'e}triques remplissables par une
m{\'e}trique d'Einstein autoduale ou antiautoduale sur \( B^{4} \), sont deux sous-vari{\'e}t{\'e}s,
telles qu'en la m{\'e}trique ronde on ait la d{\'e}composition de l'espace tangent :
\[
T(\mathcal{M}/\Diff _{0}S^{3})=T(\mathcal{M}_{+}/\Diff _{0}S^{3})\oplus T(\mathcal{M}_{-}/\Diff _{0}S^{3}).\]

\end{thm}
En particulier, cet {\'e}nonc{\'e} implique imm{\'e}diatement la conjecture de fr{\'e}quence
positive de LeBrun \cite{LeB:91} :

\begin{cor}
\textbf{\label{cor-D}} \textbf{\emph{Conjecture de fr{\'e}quence positive de LeBrun.}}
Une m{\'e}trique \( g \) sur \( S^{3} \), proche de la m{\'e}trique ronde, se d{\'e}compose
en \( g=g_{-}+g_{0}+g_{+} \), o{\`u} \( g_{-}+g_{0} \) (resp. \( g_{+}+g_{0} \))
est l'infini conforme d'une m{\'e}trique d'Einstein antiautoduale (resp. autoduale).
\end{cor}
Enfin, nous {\'e}tudions une version naturelle de ce probl{\`e}me en dimension sup{\'e}rieure.
Il est bien connu qu'existe en dimension \( 4m \) un analogue des m{\'e}triques
antiautoduales d'Einstein, les m{\'e}triques \emph{quaternion-k{\"a}hl{\'e}riennes}, dont
le prototype est la m{\'e}trique hyperbolique quaternionienne sur la boule \( B^{4m} \).
Dans \cite{Biq:00} (voir aussi le survey \cite{Biq:99}), j'ai {\'e}tudi{\'e} les infinis
conformes des m{\'e}triques quaternion-k{\"a}hl{\'e}riennes, et ai appel{\'e} ces objets, qui
vivent sur des vari{\'e}t{\'e}s de dimension \( 4m-1 \), des \emph{structures de contact
quaternioniennes} (voir section \ref{sec-11}) ; {\`a} partir de la dimension 11,
une structure de contact quaternionienne, analytique r{\'e}elle, est l'infini conforme
d'une unique m{\'e}trique quaternion-k{\"a}hl{\'e}rienne, localement d{\'e}finie ; en dimension
7, une condition additionnelle d'int{\'e}grabilit{\'e} (l'existence d'un espace de twisteurs)
est requise.

L'analogie entre la dimension 4 et la dimension sup{\'e}rieure est due {\`a} la persistance,
pour les m{\'e}triques quaternion-k{\"a}hl{\'e}riennes, d'une construction twistorielle
(\cite{Sal:82}, voir aussi \cite{Bes:87} ; on trouve la construction inverse
dans \cite{LeB:89}) ; ainsi, la structure de contact quaternionienne standard
de la sph{\`e}re \( S^{4m-1} \) a pour espace des twisteurs l'hypersurface r{\'e}elle
de \( P^{2m+1} \) d'{\'e}quation 
\[
|z^{1}|^{2}+|z^{2}|^{2}+\cdots +|z^{2m}|^{2}=|z^{2m+1}|^{2}+|z^{2m+2}|^{2},\]
 bord du domaine complexe 
\[
D=\{|z^{1}|^{2}+|z^{2}|^{2}+\cdots +|z^{2m}|^{2}<|z^{2m+1}|^{2}+|z^{2m+2}|^{2}\}\subset P^{2m+1}\]
 qui est l'espace des twisteurs de la boule \( B^{4m} \) munie de la m{\'e}trique
hyperbolique quaternionienne.

La diff{\'e}rence avec le cas de la dimension 4 est que le probl{\`e}me d'extension
a toujours une solution positive.

\begin{thm}
\label{th-E}Pour \( 4m\geq 12, \) toute structure de contact quaternionienne
sur la sph{\`e}re \( S^{4m-1} \), proche de la structure standard, est l'infini
conforme d'une m{\'e}trique quaternion-k{\"a}hl{\'e}rienne sur la boule \( B^{4m} \). Pour
\( 4m=8 \), le m{\^e}me r{\'e}sultat est valable, pourvu que la structure de contact
quaternionienne admette un espace de twisteurs (ce qui est toujours vrai en
dimension sup{\'e}rieure).
\end{thm}
{\'E}voquons maintenant les techniques utilis{\'e}es dans cet article.

La m{\'e}thode de remplissage (section \ref{sec2}) d'une structure CR est une m{\'e}thode
classique d'analyse complexe : le remplissage holomorphe suivant des disques
(d'o{\`u} la condition de positivit{\'e} des coefficients de Fourier) ; cette technique
est utilis{\'e}e notamment dans l'article de Bland-Duchamp \cite{Bla-Duc:91} sur
la description de l'espace des modules des domaines de \( \mathbf{C}^{n} \)
en termes de la repr{\'e}sentation circulaire de Lempert \cite{Lem:81}.

Pour montrer que le remplissage de la structure CR suffit pour obtenir le remplissage
de \( [g,\mathcal{Q}] \), la principale difficult{\'e} est, dans le cas \( \mathcal{Q}=0 \),
le remplissage de la structure de contact holomorphe (qui fournit la m{\'e}trique
d'Einstein) : il s'agit de montrer que la donn{\'e}e sur le bord est n{\'e}cessairement
{\`a} coefficients de Fourier positifs ; ce probl{\`e}me est trait{\'e} dans la section
\ref{sec3} : la forme particuli{\`e}re de l'{\'e}quation (non lin{\'e}aire) permet de ramener
le probl{\`e}me {\`a} l'analyse de l'op{\'e}rateur \( \overline{\partial } \) sur les (1,0)-formes.

Ce remplissage serait de peu d'int{\'e}r{\^e}t si la condition de positivit{\'e} des coefficients
de Fourier n'{\'e}tait pas n{\'e}cessaire : on montre que tel est le cas dans la section
\ref{sec4} ; le probl{\`e}me consiste, partant d'un domaine complexe proche du
mod{\`e}le, {\`a} reconstruire l'action de \( S^{1} \) sur le bord, ce qui se fait
en recherchant les disques extr{\'e}maux qui les remplissent ; on montre ici un
th{\'e}or{\`e}me g{\'e}n{\'e}ral de d{\'e}formation de tels disques dans l'espace total d'un fibr{\'e}
en disques holomorphes, dont le bord a une forme de Levi non d{\'e}g{\'e}n{\'e}r{\'e}e. A priori,
il est difficile de produire de tels disques extr{\'e}maux loin du mod{\`e}le standard
dans notre situation de signature (1,1), car le domaine n'est pas pseudo-convexe
ou pseudo-concave, donc ces disques ne semblent pas maximiser une fonctionnelle
; cependant, dans le cas d'un espace de twisteurs, leur interpr{\'e}tation en termes
de la g{\'e}om{\'e}trie conforme de dimension 4 (lemme \ref{lem-tw/conf}) donne un
espoir de construction plus g{\'e}n{\'e}rale---cette construction est le probl{\`e}me essentiel
{\`a} r{\'e}soudre pour donner une condition de remplissage pour les m{\'e}triques qui ne
sont pas proches de la m{\'e}trique ronde.

{\`A} partir de la section \ref{sec5}, on s'int{\'e}resse {\`a} une question l{\'e}g{\`e}rement
diff{\'e}rente : d{\'e}crire les espaces de modules des objets satisfaisant la condition
de coefficients de Fourier positifs, et en particulier, voir {\`a} quel point cette
condition contraint la structure. Un contrexemple est fournit par le cas des
sph{\`e}res \( S^{2n-1} \) dans \( \mathbf{C}^{n} \) : le th{\'e}or{\`e}me \ref{th-A}
est valable, mais la condition est vide pour \( n>2 \), puisqu'alors toutes
les d{\'e}formations CR se remplissent ; cela signifie que l'on peut toujours trouver
un contactomorphisme ramenant une structure CR donn{\'e}e sur une structure {\`a} coefficients
de Fourier positifs. Un ph{\'e}nom{\`e}ne semblable explique le th{\'e}or{\`e}me \ref{th-E}.

Il faut donc {\'e}tudier l'action des contactomorphismes sur les structures CR ;
la section \ref{sec5} propose une param{\'e}trisation du groupe des contactomorphismes
par des fonctions. Ici, il faut entrer un peu dans la technique : les espaces
fonctionnels naturels en g{\'e}om{\'e}trie de contact sont les espaces de Folland-Stein
de fonctions avec des d{\'e}riv{\'e}es uniquement dans les directions de contact ; le
th{\'e}or{\`e}me \ref{th-par-contact} propose une param{\'e}trisation des contactomorphismes
{\`a} \( k \) d{\'e}riv{\'e}es horizontales par des fonctions {\`a} \( k+1 \) d{\'e}riv{\'e}es horizontales.
Cette param{\'e}trisation permet, dans la section \ref{sec6}, de montrer l'existence
d'une jauge de Coulomb par rapport {\`a} l'action des contactomorphismes, dans laquelle
le maximum de coefficients n{\'e}gatifs sont tu{\'e}s ; de cette mani{\`e}re, la possibilit{\'e}
de remplissage de la structure CR se lit sur la jauge de Coulomb, comme en dimension
3 (\cite{Bla:94}, voir aussi \cite{Bur-Eps:90}).

Pour passer {\`a} la description explicite des espaces tangents des structures remplissables,
il faut calculer l'image de l'action infinit{\'e}simale des contactomorphismes sur
les structures CR ; cette action est donn{\'e}e par un hessien complexe agissant
sur les fonctions r{\'e}elles. L'analyse, dans la section \ref{sec7}, est bas{\'e}e
sur le fait que le groupe d'isom{\'e}tries \( SO_{4} \) de \( S^{3} \) agit de
mani{\`e}re homog{\`e}ne sur son espace des twisteurs \( \{|z^{1}|^{2}+|z^{2}|^{2}=|z^{3}|^{2}+|z^{4}|^{2}\}\subset P^{3} \)
; par cons{\'e}quent, la d{\'e}composition des objets suivant les repr{\'e}sentations irr{\'e}ductibles
de \( SO_{4} \) permet de se ramener {\`a} un probl{\`e}me alg{\'e}brique pour chaque repr{\'e}sentation.
On peut alors faire, dans la section \ref{sec8}, une th{\'e}orie de d{\'e}formation
compl{\`e}te pour les structures CR, dans la section \ref{sec9} pour les donn{\'e}es
conformes \( [g,\mathcal{Q}] \) : cela signifie qu'on dispose bien d'un proc{\'e}d{\'e}
fabriquant {\`a} partir d'un {\'e}l{\'e}ment de l'espace tangent indiqu{\'e} dans le th{\'e}or{\`e}me
\ref{th-gQ-rem} une donn{\'e}e \( [g,\mathcal{Q}] \) remplissable ; enfin, dans
la section \ref{sec10}, on attaque le probl{\`e}me pour les m{\'e}triques conformes
\( [g] \) : une difficult{\'e} surprenante est qu'il est difficile de garantir
que le proc{\'e}d{\'e} utilis{\'e} pour \( [g,\mathcal{Q}] \) fournisse, {\`a} partir d'une
m{\'e}trique infinit{\'e}simale, une m{\'e}trique conforme (avec \( \mathcal{Q}=0 \)) ;
d'un autre c{\^o}t{\'e}, il n'est pas {\'e}vident a priori que l'intersection \( \mathcal{B}_{+}\cap \mathcal{M} \)
soit transverse ; pour contourner la difficult{\'e}, on est oblig{\'e} de reprendre
toute une th{\'e}orie de d{\'e}formation du couple constitu{\'e} d'une structure CR int{\'e}grable
avec une forme de contact holomorphe.

Enfin, dans la section \ref{sec-11}, on traite le cas, moins difficile, de
la dimension sup{\'e}rieure. Le probl{\`e}me essentiel est similaire : le remplissage
d'une petite d{\'e}formation de la structure CR, de signature \( (2m-1,1) \). On
montre que la lin{\'e}arisation de l'action des contactomorphismes sur les structures
CR atteint tous les coefficients strictement n{\'e}gatifs, de sorte qu'il est toujours
possible de pr{\'e}senter une structure CR avec des coefficients positifs, et donc
de la remplir ; une autre d{\'e}monstration possible est probablement l'application
du th{\'e}or{\`e}me de Kiremidjian \cite{Kir:79}.

\section{Construction twistorielle}

Dans cette section, on rappelle la construction de l'espace des twisteurs d'une
vari{\'e}t{\'e} conforme de dimension 3, due {\`a} LeBrun \cite{LeB:84, LeB:85}.

\subsection{\label{sec11}Twisteurs d'une vari{\'e}t{\'e} de dimension 3}

Soit \( X^{3} \) une vari{\'e}t{\'e}. L'espace \( \mathcal{T} \) des twisteurs de
\( X \) est le fibr{\'e} des grassmanniennes en 2-plans orient{\'e}s de \( TX \).
On notera \( p:\mathcal{T}\rightarrow X \) la projection naturelle. La fibre
en \( x\in X \) de la projection \( p \) est donc la 2-sph{\`e}re \( Gr_{2}^{+}(T_{x}X) \),
et on dispose d'une involution \( \tau  \) de \( \mathcal{T} \), qui associe
{\`a} un 2-plan orient{\'e} le m{\^e}me 2-plan muni de l'orientation oppos{\'e}e. Enfin, \( \mathcal{T} \)
est muni d'une distribution de contact canonique \( H\subset T\mathcal{T} \),
d{\'e}finie au point \( P\in Gr_{2}^{+}(T_{x}X) \) par
\[
H_{P}=p_{*}^{-1}(P).\]

Supposons {\`a} pr{\'e}sent \( X \) muni d'une m{\'e}trique riemannienne \( g \). Alors
la distribution de contact \( H \) h{\'e}rite d'une structure presque complexe
\( J \), d{\'e}finie de la mani{\`e}re suivante : la connexion de Levi-Civita de \( g \)
permet de d{\'e}composer l'espace tangent en \( P\in Gr_{2}^{+}(T_{x}X) \) en
\[
T_{P}\mathcal{T}=\textrm{Ver}_{P}\oplus \textrm{Hor}_{P},\]
 o{\`u} \( \textrm{Ver}_{P} \) est l'espace tangent {\`a} la fibre de \( p \) et \( \textrm{Hor}_{P} \)
est l'espace horizontal induit par la connexion ; on peut identifier \( \textrm{Ver}_{P}=\hom (P,T_{x}X/P) \),
et \( \textrm{Hor}_{P}=T_{x}X \) via la projection \( p \), et la distribution
de contact est alors \( H_{P}=\textrm{Ver}_{P}\oplus P \) ; la m{\'e}trique \( g \)
induit une structure complexe \( j \) sur le 2-plan orient{\'e} \( P \), et la
structure complexe \( J \) sur \( H_{P} \) est d{\'e}finie par 
\begin{equation}
\label{def-J}
\begin{array}{ll}
J(u)=-u\circ j,\quad  & u\in \hom (P,T_{x}X/P),\\
J(v)=jv, & v\in P.
\end{array}
\end{equation}
Si on multiplie la m{\'e}trique par une fonction, on peut v{\'e}rifier que \( J \)
n'est pas modifi{\'e}e, de sorte que \( J \) ne d{\'e}pend que de la structure conforme
\( [g] \). La structure complexe \( J \) est en r{\'e}alit{\'e} une structure CR int{\'e}grable
sur \( \mathcal{T} \), et les fibres de la projection twistorielle \( p \)
sont holomorphes.

Il y a une variante de la construction pr{\'e}c{\'e}dente dans le cas o{\`u} est de plus
donn{\'e}e sur \( X \) une forme quadratique sans trace, \( \mathcal{Q} \) ; g{\'e}om{\'e}triquement,
la donn{\'e}e \( (g,\mathcal{Q}) \) appara{\^\i}t naturellement au bord d'une vari{\'e}t{\'e}
riemannienne comme la donn{\'e}e de la premi{\`e}re forme fondamentale et de la partie
sans trace de la seconde forme fondamentale ; le facteur conforme \( f \) agit
par
\begin{equation}
\label{2-inv}
(g,\mathcal{Q})\longrightarrow (f^{2}g,f\mathcal{Q}),
\end{equation}
 et on notera \( [g,\mathcal{Q}] \) la donn{\'e}e d'un tel couple, {\`a} l'action pr{\`e}s
du facteur conforme. La donn{\'e}e de \( \mathcal{Q} \) induit, pour un 2-plan
orient{\'e} \( P\subset T_{x}X \), un morphisme \( \widetilde{\mathcal{Q}}\in \hom (\textrm{Hor}_{P},\textrm{Ver}_{P}) \),
car
\begin{eqnarray*}
\hom (\textrm{Hor}_{P},\textrm{Ver}_{P}) & = & \hom (T_{x}X,\hom (P,T_{x}X/P))\\
 & = & \hom (T_{x}X\otimes P,T_{x}X/P),
\end{eqnarray*}
et on identifie \( T_{x}X/P=\mathbf{R} \) via la m{\'e}trique ; compte tenu de
l'invariance conforme (\ref{2-inv}), cette identification ne d{\'e}pend pas du
facteur conforme, et l'{\'e}l{\'e}ment de \( \hom (\textrm{Hor}_{P},\textrm{Ver}_{P}) \)
est donc bien d{\'e}fini. Finalement, \( \widetilde{\mathcal{Q}} \) d{\'e}finit un
nouvel espace horizontal
\[
\textrm{Hor}_{P}^{\mathcal{Q}}=\{u+\widetilde{\mathcal{Q}}u,\, u\in \textrm{Hor}_{P}\},\]
et la construction de \( J \) par les formules (\ref{def-J}) fournit {\`a} nouveau
une structure CR int{\'e}grable sur \( \mathcal{T} \), ne d{\'e}pendant que de la classe
conforme \( [g,\mathcal{Q}] \).

Dans tous les cas, l'involution \( \tau  \) est compatible avec la structure
CR, elle est anti-holomorphe : 
\begin{equation}
\label{tau-CR}
\tau ^{*}T^{0,1}\mathcal{T}=T^{1,0}\mathcal{T}.
\end{equation}

On peut remarquer que si on a une vari{\'e}t{\'e} riemannienne autoduale \( M^{4} \)
{\`a} bord \( X^{3} \), alors l'espace des twisteurs de \( M \) ---une vari{\'e}t{\'e}
complexe {\`a} bord--- a pr{\'e}cis{\'e}ment pour bord l'espace des twisteurs de \( X \),
muni des premi{\`e}re et seconde formes fondamentales induites.

Un th{\'e}or{\`e}me de LeBrun indique que toutes les vari{\'e}t{\'e}s CR feuillet{\'e}es en fibres
holomorphes proviennent de la construction pr{\'e}c{\'e}dente :

\begin{thm}
\textbf{\emph{\cite[th{\'e}or{\`e}me 7]{LeB:85}}} Soit \( \mathcal{T}^{5} \) une vari{\'e}t{\'e}
CR int{\'e}grable munie d'une projection \( p:\mathcal{T}\rightarrow X \) {\`a} fibres
\( P^{1} \) holomorphes ; alors \( \mathcal{T} \) est l'espace des twisteurs
d'un couple \( [g,\mathcal{Q}] \) sur \( X \). De plus, \( \mathcal{Q}=0 \)
si et seulement si \emph{\( \mathcal{T} \)} admet une structure de contact
holomorphe, {\`a} laquelle soient tangents les \( P^{1} \).
\end{thm}
Rappelons que sur une vari{\'e}t{\'e} CR int{\'e}grable, le fibr{\'e} \( T'=T^{\mathbf{C}}/T^{0,1} \)
est holomorphe, avec un op{\'e}rateur \( \overline{\partial } \) donn{\'e} par la formule
(\( u\in T^{0,1} \), \( v\in T^{\mathbf{C}} \))
\begin{equation}
\label{TC/T01}
\overline{\partial }_{u}v=[u,v],
\end{equation}
 la formule ayant un sens sur \( T^{\mathbf{C}}/T^{0,1} \) gr{\^a}ce {\`a} \( [T^{0,1},T^{0,1}]\subset T^{0,1} \)
; une structure de contact holomorphe est un sous-fibr{\'e} holomorphe, de codimension
1, non int{\'e}grable, de \( T' \).

Dans le cas \( \mathcal{Q}=0 \), la distribution de contact holomorphe est
construite ainsi : en un 2-plan orient{\'e} \( P\subset T_{x}X \), la projection
\( \pi _{P}\circ p_{*} \) (o{\`u} \( \pi _{P} \) est la projection orthogonale
\( T_{x}X\rightarrow P \)) peut se complexifier en
\[
\pi _{P}\circ p_{*}:T^{\mathbf{C}}_{P}\mathcal{T}/T^{0,1}_{P}\mathcal{T}\longrightarrow P^{1,0},\]
et le noyau de cette projection est la structure de contact holomorphe.

\subsection{\label{sec12}Twisteurs de l'espace hyperbolique r{\'e}el}

L'espace des twisteurs de la sph{\`e}re ronde \( S^{4} \) est le projectif complexe
\( P^{3} \) ; pour d{\'e}crire la fibration twistorielle \( p \), il est utile
d'identifier la sph{\`e}re avec le projectif quaternionien : \( S^{4}=\mathbf{H}P^{1} \);
la fibration twistorielle s'{\'e}crit alors en coordonn{\'e}es homog{\`e}nes
\[
p([z^{1}:z^{2}:z^{3}:z^{4}])=[z^{1}+jz^{2}:z^{3}+jz^{4}].\]
 L'espace hyperbolique r{\'e}el se r{\'e}alise comme une demi-sph{\`e}re : 
\[
\mathbf{R}H^{4}=\{[q^{1}:q^{2}]\in \mathbf{H}P^{1},\, |q^{1}|^{2}<|q^{2}|^{2}\}\]
 et sa fibration twistorielle est la restriction de \( p \) au domaine 
\[
\mathcal{N}=\{[z^{1}:z^{2}:z^{3}:z^{4}],\, |z^{1}|^{2}+|z^{2}|^{2}<|z^{3}|^{2}+|z^{4}|^{2}\}.\]
 Comme espace des twisteurs d'une vari{\'e}t{\'e} autoduale, \( \mathcal{N} \) est
muni :

\begin{itemize}
\item d'une structure r{\'e}elle compatible aux autres structures :
\begin{equation}
\label{def-tau-N}
\tau ([z^{1}:z^{2}:z^{3}:z^{4}])=[-\overline{z^{2}}:\overline{z^{1}}:-\overline{z^{4}}:\overline{z^{3}}];
\end{equation}

\item d'une famille de \( P^{1} \) holomorphes, r{\'e}els, {\`a} fibr{\'e} normal \( \mathcal{O}(1)\oplus \mathcal{O}(1) \)
: les fibres de \( p \).
\end{itemize}
Rappelons que, par la construction twistorielle inverse, une vari{\'e}t{\'e} complexe
\( \mathcal{N} \) munie d'une telle structure est l'espace des twisteurs d'une
m{\'e}trique (conforme) autoduale, d{\'e}finie sur l'espace des \( P^{1} \) r{\'e}els.

Jusqu'ici, on a utilis{\'e} seulement la structure conforme de \( \mathbf{R}H^{4} \),
{\'e}gale {\`a} celle de la boule \( B^{4} \) ; mais l'espace hyperbolique r{\'e}el est
muni en outre d'une m{\'e}trique d'Einstein, ce qui est {\'e}quivalent {\`a} l'existence
sur \( \mathcal{N} \) d'une structure de contact holomorphe \( \eta ^{c} \),
transverse aux \( P^{1} \) r{\'e}els ; dans notre cas,
\begin{equation}
\label{contact-hol}
\eta ^{c}=z^{1}dz^{2}-z^{2}dz^{1}-(z^{3}dz^{4}-z^{4}dz^{3}).
\end{equation}

La sph{\`e}re \( S^{3} \) appara{\^\i}t comme bord {\`a} l'infini de l'espace hyperbolique
r{\'e}el, et son espace des twisteurs \( \mathcal{T} \) est le bord de \( \mathcal{N} \),
donc se d{\'e}crit comme l'hypersurface r{\'e}elle de \( P^{3} \) :
\begin{equation}
\label{dec-T}
\mathcal{T}=\{[z^{1}:z^{2}:z^{3}:z^{4}],\, |z^{1}|^{2}+|z^{2}|^{2}=|z^{3}|^{2}+|z^{4}|^{2}\}.
\end{equation}
En particulier, \( \mathcal{T} \) se projette sur \( P^{1}\times P^{1} \)
par la projection 
\[
\pi ([z^{1}:z^{2}:z^{3}:z^{4}])=([z^{1}:z^{2}],[z^{3}:z^{4}])\]
 et appara{\^\i}t ainsi comme le \( S^{1} \)-fibr{\'e} \( \mathcal{O}(-1,1) \) sur
\( P^{1}\times P^{1} \), si on choisit comme action de \( u\in S^{1} \) 
\[
u([z^{1}:z^{2}:z^{3}:z^{4}])=[uz^{1}:uz^{2}:z^{3}:z^{4}];\]
 cette action se prolonge en une action du disque \( \Delta  \) sur \( \mathcal{N} \)
avec points fixes \( \{[0:0:z^{3}:z^{4}]\} \), c'est-{\`a}-dire la fibre twistorielle
au-dessus du point \( [0:1]\in \mathbf{H}P^{1} \); en {\'e}clatant cette fibre,
on obtient un fibr{\'e} \( \pi :\widetilde{\mathcal{N}}\rightarrow P^{1}\times P^{1} \)
qui est exactement le fibr{\'e} holomorphe en disques \( \mathcal{O}(-1,1) \) au-dessus
de \( P^{1}\times P^{1} \).

Dans cette image, les structures dont est muni l'espace des twisteurs \( \mathcal{T} \)
apparaissent naturellement : 

\begin{itemize}
\item la structure de contact \( \eta  \) est fourni par la connexion standard sur
le fibr{\'e} \( \mathcal{O}(-1,1) \), avec courbure
\begin{equation}
\label{def-eta}
-d\eta =-\omega _{1}+\omega _{2},
\end{equation}
o{\`u} \( \omega _{1} \) et \( \omega _{2} \) sont les formes de Fubini-Study
des deux \( P^{1} \) ; 
\item la structure complexe \( J \) sur \( \mathcal{T} \) est le tir{\'e} en arri{\`e}re
de la structure complexe de \( P^{1}\times P^{1} \) dans la distribution de
contact, elle est \( S^{1} \)-invariante ;
\item la structure r{\'e}elle anti-commute {\`a} l'action de \( S^{1} \), 
\begin{equation}
\label{tau-S1}
\tau \circ e^{i\theta }=e^{-i\theta }\circ \tau ;
\end{equation}

\item la forme de contact holomorphe est la restriction de (\ref{contact-hol}) au
bord ; plus concr{\`e}tement, dans des coordonn{\'e}es \( [u:uz^{2}:1:z^{4}] \), o{\`u}
\( z^{2} \) et \( z^{4} \) sont des coordonn{\'e}es sur \( P^{1}\times P^{1} \),
on obtient 
\begin{equation}
\label{def-etac}
\eta ^{c}=u^{2}dz^{2}-dz^{4};
\end{equation}
 on d{\'e}duit que \( \eta ^{c} \) est une forme {\`a} valeurs dans le fibr{\'e} \( \mathcal{O}(0,2) \).
\end{itemize}
On remarquera que, par rapport {\`a} la distribution de contact holomorphe, les
\( P^{1} \) r{\'e}els {\`a} l'int{\'e}rieur de \( \mathcal{N} \) sont transverses, tandis
que les \( P^{1} \) r{\'e}els au bord de \( \mathcal{N} \) sont tangents ; cela
traduit le fait que la m{\'e}trique hyperbolique explose au bord, avec un p{\^o}le d'ordre
2. En revanche, la structure conforme s'{\'e}tend de mani{\`e}re lisse au bord. Toutes
les m{\'e}triques que nous construirons auront ce comportement.

\section{\label{sec2}Extension holomorphe suivant les disques}

Dans cette section, nous mettons en \oe{}uvre une construction pour d{\'e}former
la structure complexe de l'espace des twisteurs \( \mathcal{N} \), en utilisant
l'id{\'e}e classique de prolongement holomorphe le long de disques.

\subsection{D{\'e}formation de la structure complexe}

Commen{\c c}ons par quelques rappels sur la th{\'e}orie de d{\'e}formation. Fixons une vari{\'e}t{\'e}
complexe \( D \), une d{\'e}formation de la structure presque complexe de \( D \)
peut {\^e}tre repr{\'e}sent{\'e}e par une (0,1)-forme \( \phi  \) {\`a} valeurs dans \( T^{1,0} \),
de sorte que le nouveau \( T^{0,1} \) soit {\'e}gal {\`a} 
\[
\{X+\phi _{X},\, X\in T^{0,1}D_{0}\}.\]
 La d{\'e}formation est int{\'e}grable si et seulement si 
\begin{equation}
\label{int-int}
\overline{\partial }\phi +\frac{1}{2}[\phi ,\phi ]=0,
\end{equation}
 o{\`u} le crochet de deux (0,1)-formes {\`a} valeurs dans \( T^{1,0} \) est une (0,2)-forme
{\`a} valeurs dans \( T^{1,0} \) d{\'e}finie par
\[
\frac{1}{2}[\phi ,\phi ]_{X,Y}=[\phi _{X},\phi _{Y}]-\phi _{[X,\phi _{Y}]^{0,1}}-\phi _{[Y,\phi _{X}]^{0,1}}.\]

Supposons {\`a} pr{\'e}sent que \( D \) soit l'espace total d'un fibr{\'e} holomorphe en
disques sur une vari{\'e}t{\'e} complexe \( X \), et que le fibr{\'e} soit muni d'une m{\'e}trique
\( h:D\rightarrow \mathbf{R}_{+} \) telle que \( h(ux)=|u|^{2}h(x) \) si \( u\in \Delta  \)
; alors on dispose d'une 1-forme de connexion, en notant \( d^{C}=i(\overline{\partial }-\partial ) \),
\begin{equation}
\label{eta-dCh}
\eta =\frac{1}{2}d^{C}\ln h,
\end{equation}
 se restreignant sur chaque disque {\`a} la forme angulaire \( d\theta  \), de
sorte que la courbure de \( h \) est la 2-forme sur \( X \) d{\'e}finie par 
\[
F=id\eta ;\]
sur chaque fibr{\'e} en cercles \( h=cst \), le noyau de la forme \( \eta  \)
d{\'e}finit un espace horizontal \( H \), si bien qu'en chaque point \( x\in D \)
on a la d{\'e}composition (\( V_{x} \) d{\'e}signant l'espace tangent {\`a} la fibre de
\( D\rightarrow X \))
\[
T_{x}D=V_{x}\oplus H_{x}.\]

Sur le bord \( M=\partial D \), on a une description similaire des d{\'e}formations
CR : puisque deux structures de contact proches sont diff{\'e}omorphes, on peut
supposer que la structure de contact \( H \) est fixe, et les d{\'e}formations
de la structure CR sont param{\'e}tr{\'e}es par les (0,1)-formes horizontales \( \phi  \)
{\`a} valeurs dans \( T^{1,0} \), et, {\`a} nouveau, la d{\'e}formation est int{\'e}grable
si et seulement si 
\[
\overline{\partial }\phi +\frac{1}{2}[\phi ,\phi ]=0;\]
 ici il faut noter que le fibr{\'e} \( T^{1,0} \) n'est pas holomorphe, voir (\ref{TC/T01}),
si bien que \( \overline{\partial }\phi  \) est en r{\'e}alit{\'e} {\`a} valeurs dans le
fibr{\'e} holomorphe \( T^{\mathbf{C}}M/T^{0,1}M \) ; la pr{\'e}sence de la forme de
contact \( \eta  \) permet de choisir un champ de Reeb \( R \) ---l'action
infinit{\'e}simale de \( S^{1} \)---, et la projection de l'{\'e}quation sur \( T^{1,0} \)
et \( \mathbf{C}R \) donne les deux conditions 
\begin{equation}
\label{int-bord}
\overline{\partial }_{H}\phi +\frac{1}{2}[\phi ,\phi ]=0,\qquad \phi \lrcorner F=0,
\end{equation}
 o{\`u}, par abus de notation, \( \lrcorner  \) est la contraction suivie du produit
ext{\'e}rieur : 
\[
\Omega ^{0,1}\otimes T^{1,0}\otimes \Omega ^{1,1}\rightarrow \Omega ^{0,1}\otimes \Omega ^{0,1}\rightarrow \Omega ^{0,2};\]
 plus loin, on utilisera aussi l'op{\'e}rateur plus g{\'e}n{\'e}ral 
\[
\lrcorner :(\Omega ^{0,1}\otimes T^{1,0})\otimes \Omega ^{i,j}\rightarrow \Omega ^{i-1,j+1};\]
 on peut remarquer que la seconde condition dans (\ref{int-bord}) dit que \( F \)
reste de type (1,1) pour la nouvelle structure complexe.

Bland-Duchamp {\'e}tudient les d{\'e}formations \( \phi  \) de \( D \) qui ne modifient
pas la structure complexe des disques et pr{\'e}servent \( H \) :

\begin{thm}
\label{th-BD}\textbf{\emph{\cite[th{\'e}or{\`e}mes 7.2 et 9.3]{Bla-Duc:91}}} Supposons
que \( \phi  \), d{\'e}fini sur le fibr{\'e} en disques \( D \), s'annule sur \( V \)
et satisfait \( \phi (H^{0,1})\subset H^{1,0} \), alors \( \phi  \) est une
d{\'e}formation int{\'e}grable de \( D \) si et seulement si
\begin{enumerate}
\item \( \phi  \) est holomorphe le long de chaque disque du fibr{\'e} ;
\item \( \phi  \) satisfait les {\'e}quations (\ref{int-bord}) sur le bord \( \partial D \).
\end{enumerate}
\end{thm}
Remarquons que la premi{\`e}re condition a un sens, puisque le fibr{\'e} \( H_{x} \)
sur chaque disque est {\'e}gal {\`a} l'espace fixe \( T_{\pi (x)}X \) : elle signifie
que dans une base constante de \( H_{x} \), les coefficients de \( \phi  \)
sont holomorphes le long du disque.

Le th{\'e}or{\`e}me dit en particulier qu'{\'e}tant donn{\'e} \( \phi  \) int{\'e}grable sur le
bord \( \partial D \), {\`a} coefficients de Fourier positifs ou nuls par rapport
{\`a} l'action de \( S^{1} \), suffisamment petite, l'extension holomorphe de \( \phi  \)
le long de chaque disque fournit une d{\'e}formation int{\'e}grable de \( D \). La
d{\'e}monstration est simple : d'une part, si \( \phi  \) est assez petite au bord,
l'extension holomorphe {\`a} l'int{\'e}rieur est petite et correspond donc bien {\`a} une
d{\'e}formation de la structure complexe ; d'autre part, la condition que \( \phi  \)
soit holomorphe le long des disques et la condition (\ref{int-bord}) impliquent
qu'en r{\'e}alit{\'e}, la condition (\ref{int-bord}) est v{\'e}rifi{\'e}e sur toutes les hypersurfaces
\( h=cst \), ce qui donne l'int{\'e}grabilit{\'e} (\ref{int-int}).

\subsection{Application aux d{\'e}formations autoduales de la boule}

Appliquons le th{\'e}or{\`e}me \ref{th-BD} {\`a} l'espace des twisteurs \( \mathcal{N} \)
de l'espace hyperbolique r{\'e}el d{\'e}fini en \ref{sec12} ; l'{\'e}clatement d'un \( P^{1} \)
nous a donn{\'e} le fibr{\'e} en disques \( \widetilde{\mathcal{N}}=\mathcal{O}(-1,1) \)
au-dessus de \( P^{1}\times P^{1} \). Rappelons que \( \mathcal{N} \) et \( \widetilde{\mathcal{N}} \)
sont munis de la structure r{\'e}elle \( \tau  \) d{\'e}finie en (\ref{def-tau-N}).

\begin{lem}
\label{lem-fill-J}Si on a une petite d{\'e}formation \( \phi \in \hom (H^{0,1},H^{1,0}) \)
de la structure CR int{\'e}grable de \( \mathcal{T}=\partial \mathcal{N} \), {\`a}
coefficients de Fourier positifs ou nuls par rapport {\`a} l'action de \( S^{1} \)
sur \( \mathcal{T} \), alors \( \phi  \) s'{\'e}tend en une d{\'e}formation de la
structure complexe de \( \mathcal{N} \).

De plus, si \( \phi  \) est compatible {\`a} la structure r{\'e}elle \( \tau  \),
alors le prolongement reste compatible {\`a} \( \tau  \).
\end{lem}
\begin{proof}
En appliquant le th{\'e}or{\`e}me \ref{th-BD}, on peut prolonger \( \phi  \) en une
d{\'e}formation int{\'e}grable sur \( \widetilde{\mathcal{N}} \). Supposons dans un
premier temps que les coefficients de Fourier de \( \phi  \) soient strictement
positifs, alors les coefficients de l'extension s'annulent sur la section nulle,
qui demeure donc \( P^{1}\times P^{1} \) avec fibr{\'e} normal \( \mathcal{O}(-1,1) \)
; on peut alors contracter le premier \( P^{1} \) pour r{\'e}cup{\'e}rer une d{\'e}formation
de \( \mathcal{N} \). Dans le cas g{\'e}n{\'e}ral, les coefficients de \( \phi  \)
ne s'annulent pas sur la section nulle, mais comme \( \phi  \) est {\`a} valeurs
dans \( H^{1,0} \), la section nulle reste une sous-vari{\'e}t{\'e} holomorphe pour
la nouvelle structure complexe ; comme les d{\'e}formations holomorphes de \( P^{1}\times P^{1} \)
sont triviales, elle reste biholomorphe {\`a} \( P^{1}\times P^{1} \) ; le fibr{\'e}
normal lui aussi ne peut {\^e}tre d{\'e}form{\'e}, donc reste {\'e}gal {\`a} \( \mathcal{O}(-1,1) \),
et la m{\^e}me conclusion s'applique.

Pour la seconde partie de l'{\'e}nonc{\'e}, il est facile de voir que si au bord \( \phi  \)
est compatible {\`a} la structure r{\'e}elle \( \tau  \), c'est-{\`a}-dire \( \tau ^{*}\overline{\phi }=-\phi  \),
alors cette propri{\'e}t{\'e} est pr{\'e}serv{\'e}e dans le prolongement holomorphe.
\end{proof}
\begin{cor}
\label{cor-fill-g2}Si on a une petite d{\'e}formation \( [g,\mathcal{Q}] \) de
la m{\'e}trique standard de \( S^{3} \), dont l'espace des twisteurs satisfait
les hypoth{\`e}ses du lemme \ref{lem-fill-J}, alors \( [g,\mathcal{Q}] \) est
le bord d'une m{\'e}trique autoduale sur la boule \( B^{4} \).
\end{cor}
\begin{proof}
La construction r{\'e}alis{\'e}e dans la section \ref{sec11} permet d'associer {\`a} une
d{\'e}formation de \( [g,\mathcal{Q}] \) une d{\'e}formation de \( \mathcal{T} \)
en gardant fixe la structure de contact et la structure r{\'e}elle ; le lemme \ref{lem-fill-J}
permet d'{\'e}tendre cette d{\'e}formation en une d{\'e}formation de \( \mathcal{N} \).
Pour avoir une structure d'espace des twisteurs, il reste {\`a} fabriquer la famille
des \( P^{1} \) r{\'e}els : remarquons que cette famille est donn{\'e}e sur \( \mathcal{T}=\partial N \),
et il faut la prolonger {\`a} l'int{\'e}rieur ; comme les \( P^{1} \) initiaux ont
pour fibr{\'e} normal \( \mathcal{O}(1)\oplus \mathcal{O}(1) \), la th{\'e}orie de
d{\'e}formation de Kodaira permet de suivre la d{\'e}formation. On a ainsi construit
l'espace des twisteurs du prolongement souhait{\'e} de \( [g,\mathcal{Q}] \) sur
\( B^{4} \).
\end{proof}
\begin{rem}
La construction de l'espace des twisteurs de \( [g,\mathcal{Q}] \) faite dans
la section \ref{sec11} n'a en r{\'e}alit{\'e} aucune chance de fournir directement
une d{\'e}formation \( \phi  \) de la structure CR de \( \mathcal{T} \) {\`a} coefficients
de Fourier positifs ou nuls ; l'id{\'e}e pour appliquer le corollaire consiste {\`a}
trouver {\'e}ventuellement un contactomorphisme de \( \mathcal{T} \), pr{\'e}servant
\( \tau  \), ramenant \( \phi  \) {\`a} un tenseur {\`a} coefficients de Fourier positifs
ou nuls.
\end{rem}

\section{\label{sec3}Extension de la forme de contact holomorphe}

Si l'on part sur \( S^{3} \) d'une m{\'e}trique dont l'espace des twisteurs \( \mathcal{T} \)
satisfait les hypoth{\`e}ses du lemme \ref{lem-fill-J}, le corollaire \ref{cor-fill-g2}
fournit un remplissage autodual, avec espace des twisteurs \( \mathcal{N} \)
; pour obtenir un remplissage autodual Einstein, il faut de plus prolonger {\`a}
\( \mathcal{N} \) la structure de contact holomorphe d{\'e}finie sur \( \mathcal{T}=\partial \mathcal{N} \).
Ce probl{\`e}me peut {\^e}tre repos{\'e} sous la forme : {\'e}tant donn{\'e}e sur \( \partial \mathcal{N} \)
une section de la grassmannienne des 2-plans complexes de \( \mathcal{N} \),
peut-on {\'e}tendre cette section en une section holomorphe sur \( \mathcal{N} \)
? nous proposons ci-dessous une solution de ce probl{\`e}me de Plateau complexe,
pour les structures complexes proches du \( \mathcal{N} \) standard.

\subsection{Probl{\`e}me infinit{\'e}simal}

Rappelons que l'espace des twisteurs \( \mathcal{T} \) de \( S^{3} \) est
le fibr{\'e} en cercles \( \pi :\mathcal{O}(-1,1)\rightarrow P^{1}\times P^{1} \)
; on peut identifier, en notant \( P^{1}_{1} \) et \( P^{1}_{2} \) les deux
facteurs du produit \( P^{1}\times P^{1} \), 
\begin{equation}
\label{dec-T'}
T'\mathcal{T}:=T^{\mathbf{C}}\mathcal{T}/T^{0,1}\mathcal{T}=\mathbf{C}R\oplus T^{1,0}\mathcal{T}=\mathbf{C}R\oplus T^{1,0}P^{1}_{1}\oplus T^{1,0}P^{1}_{2}.
\end{equation}
 Pr{\'e}cisons la structure holomorphe du fibr{\'e} dual 
\begin{equation}
\label{dec-O'}
\Omega '\mathcal{T}=(T'\mathcal{T})^{*}=\mathbf{C}\eta \oplus \Omega ^{1,0}P^{1}_{1}\oplus \Omega ^{1,0}P^{1}_{2}.
\end{equation}
 Les deux sous-fibr{\'e}s \( \Omega ^{1,0}P^{1}_{i} \), tir{\'e}s en arri{\`e}re depuis
la base \( P^{1}\times P^{1} \), sont holomorphes, mais pas \( \eta  \) :
compte tenu de (\ref{def-eta}), on obtient l'op{\'e}rateur \( \overline{\partial } \)
dans la d{\'e}composition (\ref{dec-O'}) : 
\begin{equation}
\label{db-O'}
\overline{\partial }^{\Omega '\mathcal{T}}=\left( \begin{array}{ccc}
\overline{\partial } &  & \\
\omega _{1} & \overline{\partial } & \\
-\omega _{2} &  & \overline{\partial }
\end{array}\right) .
\end{equation}

La forme de contact holomorphe \( \eta ^{c} \) d{\'e}finie en (\ref{def-etac})
est une section de \( \Omega '\mathcal{T}\otimes L \), o{\`u} \( L=\mathcal{O}(0,2) \).
L'action de \( S^{1} \) sur \( \mathcal{T} \) se rel{\`e}ve naturellement au fibr{\'e}
\( \Omega '\mathcal{T}\otimes L \), car celui-ci provient de la base ; son
op{\'e}rateur \( \overline{\partial } \) est {\'e}galement \( S^{1} \)-invariant.
Cela permet d'{\'e}crire le lemme suivant.

\begin{lem}
\label{lem-31}L'op{\'e}rateur \( \overline{\partial } \) du fibr{\'e} \( \Omega '\mathcal{T}\otimes L \)
n'a pas de noyau sur les sections {\`a} coefficients de Fourier strictement n{\'e}gatifs
; sur les sections \( S^{1} \)-invariantes, son noyau est r{\'e}duit aux sections
constantes de \( \Omega ^{1,0}P^{1}_{2}\otimes \mathcal{O}(0,2)=\mathbf{C} \).
\end{lem}
\begin{proof}
Le fibr{\'e} \( \mathcal{O}(-1,1) \) sur \( \mathcal{T} \) a une section tautologique
\( t \), holomorphe, de poids 1 pour l'action de \( S^{1} \) ; {\'e}tant donn{\'e}
un fibr{\'e} holomorphe \( S^{1} \)-invariant \( E \), la donn{\'e}e d'une section
holomorphe \( s \) de \( E \) avec poids \( k \) par rapport {\`a} l'action de
\( S^{1} \) est {\'e}quivalente {\`a} la donn{\'e}e de la section holomorphe \( S^{1} \)-invariante
\( s\otimes t^{-k} \) du fibr{\'e} \( E\otimes \mathcal{O}(-1,1)^{\otimes (-k)}=E\otimes \mathcal{O}(k,-k) \).
Appliquons cette observation au fibr{\'e} \( E=\Omega '\mathcal{T}\otimes L \)
muni de la structure (\ref{dec-O'}) : une section holomorphe de poids \( k \)
fournit une section holomorphe invariante de \( E\otimes \mathcal{O}(k,-k)=\Omega '\mathcal{T}\otimes \mathcal{O}(k,-k+2) \)
; sous la d{\'e}composition (\ref{dec-O'}), une telle section se d{\'e}compose en \( s=s_{0}+s_{1}+s_{2} \)
avec \( s_{0},s_{1},s_{2} \) sections respectives de \( \mathcal{O}(k,-k+2) \),
\( \mathcal{O}(k-2,-k+2) \) et \( \mathcal{O}(k,-k) \), et la section est
holomorphe si 
\[
\overline{\partial }s_{0}=0,\quad \overline{\partial }s_{1}+\omega _{1}s_{0}=0,\quad \overline{\partial }s_{2}-\omega _{2}s_{0}=0;\]
 comme les sections sont \( S^{1} \)-invariantes, elles descendent sur \( P^{1}\times P^{1} \)
et les {\'e}quations ont lieu sur \( P^{1}\times P^{1} \). 

Si \( k<0 \), alors \( H^{0}(P^{1}\times P^{1},\mathcal{O}(k,-k+2))=0 \) donc
\( s_{0}=0 \), et les {\'e}quations deviennent simplement 
\[
\overline{\partial }s_{1}=0,\quad \overline{\partial }s_{2}=0,\]
 mais, {\`a} nouveau, \( H^{0}(P^{1}\times P^{1},\mathcal{O}(k-2,-k+2))=H^{0}(P^{1}\times P^{1},\mathcal{O}(k,-k))=0 \)
si \( k<0 \) donc \( s=0 \).

Si \( k=0 \), alors \( s_{0}\in H^{0}(P^{1}\times P^{1},\mathcal{O}(0,2))=H^{0}(P^{1}_{2},\mathcal{O}(2)) \),
mais la projection de l'{\'e}quation \( \overline{\partial }s_{1}+\omega _{1}s_{0}=0 \)
sur \( H^{1}(P^{1}\times P^{1},\mathcal{O}(-2,2))=H^{0}(P^{1}_{2},\mathcal{O}(2)) \)
implique en r{\'e}alit{\'e} \( s_{0}=0 \) ; on d{\'e}duit alors de la m{\^e}me mani{\`e}re \( s_{1}=0 \)
et \( s_{2}=cst \).
\end{proof}
Ce lemme implique en particulier que, par prolongement holomorphe le long des
disques, toutes les sections holomorphes du fibr{\'e} \( \Omega '\mathcal{T}\otimes L \)
au bord se prolongent sur le fibr{\'e} en disques ; en outre, sur la section nulle,
seule la composante sur \( \Omega ^{1,0}P^{1}_{2} \) peut {\^e}tre non triviale,
ce qui signifie que la section descend en r{\'e}alit{\'e} sur l'espace des twisteurs
\( \mathcal{N} \) (obtenu en contractant le \( P^{1}_{1} \) de la section
nulle). Cette observation est {\`a} la base de la d{\'e}monstration de l'existence d'un
prolongement {\`a} l'int{\'e}rieur de la forme de contact holomorphe.

\begin{cor}
\label{cor-31}Pour une section \( s \) de \( \Omega '\otimes L \) sur \( \mathcal{T} \),
{\`a} coefficients de Fourier strictement n{\'e}gatifs, on a une estimation 
\[
\left\Vert \overline{\partial }s\right\Vert _{L^{2}}\geq c(\left\Vert \nabla _{H}s\right\Vert _{L^{2}}+\left\Vert s\right\Vert _{L^{2}});\]
 la m{\^e}me estimation reste valable pour une section {\`a} coefficients de Fourier
n{\'e}gatifs ou nuls, s'annulant sur \( T^{1,0}P^{1}_{2} \).
\end{cor}
\begin{proof}
La vari{\'e}t{\'e} CR \( \mathcal{T} \) est de signature (1,1), l'op{\'e}rateur \( \overline{\partial } \)
est donc hypoelliptique sur les 0-formes, et l'absence de noyau, cons{\'e}quence
du lemme (\ref{lem-31}), implique les estimations du corollaire.

Bien qu'il s'agisse d'une th{\'e}orie bien connue, il peut {\^e}tre utile au lecteur
de voir de mani{\`e}re {\'e}l{\'e}mentaire comment ces estimations sont une cons{\'e}quence
de la signature (1,1) du fibr{\'e} \( \mathcal{O}(-1,1) \). Expliquons cela bri{\`e}vement
dans le cas d'une fonction \( s \) d'un fibr{\'e} \( E \) provenant de la base
\( P^{1}\times P^{1} \). On peut d{\'e}composer \( s \) suivant en s{\'e}ries de Fourier
sur chaque cercle : 
\[
s=\sum _{k}s_{k};\]
 comme on l'a vu dans la d{\'e}monstration du lemme \ref{lem-31}, si \( t \) est
la section tautologique du fibr{\'e} \( \mathcal{O}(-1,1) \), alors \( s_{k}t^{-k} \)
est une section \( S^{1} \)-invariante de \( \mathcal{O}(k,-k) \) ; rappelons
la formule de Weitzenb{\"o}ck, pour une section \( \sigma  \) d'un fibr{\'e} \( E \)
sur la vari{\'e}t{\'e} k{\"a}hl{\'e}rienne \( X \) : 
\begin{equation}
\label{for-Wei}
\int _{X}\left| \overline{\partial }\sigma \right| ^{2}=\frac{1}{2}\int _{X}\left| \nabla \sigma \right| ^{2}-\left\langle i\Lambda F^{E}\sigma ,\sigma \right\rangle ,
\end{equation}
 o{\`u} \( \Lambda  \) est la contraction par la forme de K{\"a}hler ; maintenant,
remarquons que 
\[
iF^{\mathcal{O}(k,-k)}=k(-\omega _{1}+\omega _{2});\]
 si \( k\geq 0 \), nous appliquons la formule de Weitzenb{\"o}ck sur chaque \( \{x\}\times P^{1}_{2} \)
: 
\[
\int _{X}\left| \overline{\partial }(s_{k}t^{-k})\right| ^{2}\geq \int _{X}\left| \overline{\partial }_{P^{1}_{2}}(s_{k}t^{-k})\right| ^{2}\geq \int _{X}(k-\sup |F^{E}|)\left| s_{k}\right| ^{2},\]
 tandis que pour \( k\leq 0 \), le m{\^e}me r{\'e}sultat est obtenu en utilisant plut{\^o}t
le signe de la courbure de \( \mathcal{O}(k,-k) \) sur \( P^{1}_{1} \), et
par cons{\'e}quent, 
\[
\left\Vert \overline{\partial }s\right\Vert _{L^{2}}^{2}+\left\Vert s\right\Vert _{L^{2}}^{2}\geq c\sum (1+|k|)\left\Vert s_{k}\right\Vert ^{2};\]
 en r{\'e}utilisant (\ref{for-Wei}) pour obtenir un contr{\^o}le sur \( \nabla (s_{k}t^{-k}) \),
on obtient 
\[
\left\Vert \overline{\partial }s\right\Vert _{L^{2}}^{2}+\left\Vert s\right\Vert _{L^{2}}^{2}\geq \left\Vert \nabla _{H}s\right\Vert _{L^{2}}^{2};\]
 ce contr{\^o}le implique aussi le contr{\^o}le d'une demi-d{\'e}riv{\'e}e dans la direction
transverse {\`a} la distribution de contact, c'est-{\`a}-dire l'estimation, d{\'e}montr{\'e}e
directement, sur \( \sum k\left\Vert s_{k}\right\Vert ^{2} \).

Par un argument standard d'analyse fonctionnelle, cette estimation implique
que le noyau de \( \overline{\partial } \) est de dimension finie, et, dans
le cas o{\`u} il est nul, l'estimation du corollaire.
\end{proof}

\subsection{\label{sec32}Probl{\`e}me non lin{\'e}aire}

La d{\'e}formation de la structure complexe, \( \phi :T^{0,1}\rightarrow T^{1,0} \),
induit au dual \( \phi ^{t}:\Omega ^{1,0}\rightarrow \Omega ^{0,1} \), que
l'on peut aussi {\'e}crire comme la contraction \( \Omega ^{0,1}\otimes T^{1,0}\otimes \Omega ^{1,0}\rightarrow \Omega ^{0,1} \),
 
\[
\phi ^{t}(\xi )=\phi \lrcorner \xi ;\]
les (1,0)-formes pour la nouvelle structure complexe sont param{\'e}tr{\'e}es par les
\[
\xi -\phi \lrcorner \xi ,\]
 pour \( \xi \in \Omega ^{1,0} \) ; notons que cette formule a {\'e}galement un
sens sur \( \Omega '=\mathbf{C}\eta \oplus \Omega ^{1,0} \). Une structure
de contact complexe pour la nouvelle structure complexe peut donc {\^e}tre repr{\'e}sent{\'e}e
par la 1-forme {\`a} valeurs dans \( L \),
\[
\varpi ^{c}=(\eta ^{c}+\xi )-\phi \lrcorner (\eta ^{c}+\xi ),\]
 o{\`u} \( \xi  \) est une section de \( \Omega '\mathcal{T}\otimes L \), que
l'on peut soumettre {\`a} la normalisation 
\[
\xi |_{T^{1,0}P^{1}_{2}}=0;\]
 le noyau de la forme \( \varpi ^{c} \) est holomorphe si et seulement si,
pour tout \( X\in T^{0,1} \), 
\begin{equation}
\label{oc-hol}
\iota _{X+\phi _{X}}d\varpi ^{c}=\alpha _{X}\varpi ^{c},\quad \alpha \in \Omega ^{0,1};
\end{equation}
 on remarquera que l'op{\'e}rateur \( d \) fait a priori intervenir une connexion
sur le fibr{\'e} \( L \), mais un choix diff{\'e}rent de connexion ne modifie le membre
de gauche que par un terme proportionnel {\`a} \( \varpi ^{c} \), donc sans influence
sur l'{\'e}quation.

\subsubsection*{Digression : d{\'e}composition de la diff{\'e}rentielle ext{\'e}rieure}

Avant de poursuivre l'{\'e}tude de l'{\'e}quation (\ref{oc-hol}), nous avons besoin
de pr{\'e}ciser quelques d{\'e}compositions de la diff{\'e}rentielle ext{\'e}rieure \( d \)
sur la vari{\'e}t{\'e} CR \( \mathcal{T} \).

Pour une fonction \( f \), on peut d{\'e}composer la diff{\'e}rentielle en 
\[
df=\overline{\partial }f+\partial 'f,\quad \overline{\partial }f\in \Omega ^{0,1},\quad \partial 'f\in \Omega '.\]

Passons maintenant aux 1-formes : puisque \( [T^{1,0},T^{1,0}]\subset T^{1,0} \),
on a 
\[
d\Omega ^{0,1}\subset \Omega ^{0,2}+\Omega ^{1,1}+\mathbf{C}\eta \otimes \Omega ^{0,1};\]
 pour une (0,1)-forme \( \alpha \in \Omega ^{0,1} \), on peut donc d{\'e}composer
\[
d\alpha =\overline{\partial }\alpha +\partial '\alpha ,\quad \left\{ \begin{array}{l}
\overline{\partial }\alpha \in \Omega ^{0,2},\\
\partial '\alpha \in \Omega ^{1,1}+\mathbf{C}\eta \otimes \Omega ^{0,1}=\Omega '\otimes \Omega ^{0,1}.
\end{array}\right. \]
 De mani{\`e}re similaire, pour une (1,0)-forme \( \alpha \in \Omega ^{1,0} \),
on a la d{\'e}composition 
\[
d\alpha =\overline{\partial }\alpha +\partial '\alpha ,\quad \left\{ \begin{array}{l}
\overline{\partial }\alpha \in \Omega ^{1,1},\\
\partial '\alpha \in \Omega ^{2,0}+\mathbf{C}\eta \otimes \Omega ^{1,0}=\Lambda ^{2}\Omega '.
\end{array}\right. \]
 Enfin, pour une fonction \( f \), on a 
\[
d(f\eta )=(\overline{\partial }f\wedge \eta +fd\eta )+\partial f\wedge \eta \in \Omega ^{0,1}\otimes \Omega '+\Lambda ^{2}\Omega '\]
 et on peut d{\'e}finir \( \overline{\partial }(f\eta ) \) et \( \partial (f\eta ) \)
comme les deux morceaux de cette diff{\'e}rentielle.

En particulier, nous avons donc la diff{\'e}rentielle sur \( \Omega ' \) qui se
d{\'e}compose en \( d=\overline{\partial }+\partial ' \), avec 
\[
\overline{\partial }\Omega '\subset \Omega ^{0,1}\otimes \Omega ',\quad \partial '\Omega '\subset \Lambda ^{2}\Omega '.\]
Ces observations m{\`e}nent facilement au lemme plus g{\'e}n{\'e}ral suivant.

\begin{lem}
\label{dec-d}Notons \( \tilde{\Omega }^{i,j}=\Lambda ^{i}\Omega '\otimes \Omega ^{0,j}\subset \Omega ^{i+j} \),
alors la diff{\'e}rentielle ext{\'e}rieure se d{\'e}compose en \( d=\overline{\partial }+\partial ' \),
avec 
\[
\overline{\partial }\tilde{\Omega }^{i,j}\subset \tilde{\Omega }^{i+1,j},\quad \partial '\tilde{\Omega }^{i,j}\subset \tilde{\Omega }^{i,j+1};\]
 en particulier, on a \( \overline{\partial }\partial '+\partial '\overline{\partial }=0. \)\qed 
\end{lem}

\subsubsection*{Forme de contact holomorphe}

Revenons {\`a} l'{\'e}quation (\ref{oc-hol}) : les deux membres sont dans l'espace
\( \Omega ' \) pour la nouvelle structure complexe induite par \( \phi  \),
l'{\'e}quation est donc {\'e}quivalente {\`a} sa projection sur \( \Omega ' \) pour la
structure initiale, ce qui nous donne 
\begin{equation}
\label{etac+xi-hol0}
\iota _{X}\overline{\partial }\xi -\iota _{X}\partial '\phi \lrcorner (\eta ^{c}+\xi )+\iota _{\phi _{X}}(d\eta ^{c}+\partial '\xi )=\alpha _{X}(\eta ^{c}+\xi ).
\end{equation}
Notons l'identit{\'e}, pour toute forme \( \xi \in \Omega ' \), 
\begin{equation}
\label{id-dphi-phid}
\partial '(\phi \lrcorner \xi )-\phi \lrcorner \partial '\xi =(\partial '\phi )\lrcorner \xi -\phi \lrcorner \nabla ^{1,0}\xi ;
\end{equation}
 cette identit{\'e} ne fait intervenir qu'en apparence une d{\'e}riv{\'e}e covariante \( \nabla ^{1,0} \)
sur \( T^{1,0} \), car, en r{\'e}alit{\'e}, le membre de droite ne d{\'e}pend pas du choix
pr{\'e}cis effectu{\'e}. Nous obtenons finalement l'{\'e}quation 
\begin{equation}
\label{etac+xi-hol}
\overline{\partial }\xi -(\partial '\phi )\lrcorner (\eta ^{c}+\xi )+\phi \lrcorner \nabla ^{1,0}(\eta ^{c}+\xi )=\alpha \otimes (\eta ^{c}+\xi ).
\end{equation}
 En {\'e}valuant l'{\'e}quation sur un vecteur dans \( T^{1,0}P^{1}_{2} \), on obtient
\( \alpha =\alpha (\xi ) \), avec \( \alpha (\xi ) \) ne d{\'e}pendant de \( \xi  \)
qu'{\`a} l'ordre 0 ; on peut alors r{\'e}crire l'{\'e}quation comme 
\[
\overline{\partial }\xi =(\partial '\phi )\lrcorner (\eta ^{c}+\xi )-\phi \lrcorner \nabla ^{1,0}(\eta ^{c}+\xi )+\alpha (\xi )\otimes (\eta ^{c}+\xi ).\]
 Supposons \( \phi  \) {\`a} coefficients de Fourier positifs ou nuls, et projetons
cette {\'e}quation sur les coefficients de Fourier strictement n{\'e}gatifs : comme
\( \eta ^{c} \) a un coefficient de poids 0 sur le \( \Omega ^{1,0}P^{1}_{2} \)
et un coefficient de poids 2 sur le \( \Omega ^{1,0}P^{1}_{1} \), on obtient,
en notant \( \xi _{<0} \) la partie de \( \xi  \) {\`a} coefficients n{\'e}gatifs
ou nuls, 
\[
\overline{\partial }\xi _{<0}=\pi _{<0}\left( (\partial '\phi )\lrcorner \xi -\phi \lrcorner \nabla ^{1,0}\xi +\alpha (\xi )\otimes (\eta ^{c}+\xi )\right) ;\]
 compte tenu de la forme des {\'e}quations, dont les termes non lin{\'e}aires sont seulement
des multiplications, on obtient 
\begin{equation}
\label{est-xi}
\left\Vert \overline{\partial }\xi _{<0}\right\Vert _{L^{2}}\leq c\left( \left\Vert \xi _{<0}\right\Vert _{L^{2}}+\left\Vert \nabla _{H}\xi _{<0}\right\Vert _{L^{2}}\right) \left( \left\Vert \xi \right\Vert _{L^{\infty }}+\left\Vert \phi \right\Vert _{C^{1}}\right) .
\end{equation}
De cette estimation, on d{\'e}duit le r{\'e}sultat suivant.

\begin{thm}
\label{th-fill-etac}Supposons que l'on ait une perturbation \( \mathcal{T}_{\phi } \)
de la structure complexe de \( \mathcal{T} \) par un \( \phi  \) suffisamment
petit, {\`a} coefficients de Fourier positifs ou nuls (de sorte que la perturbation
s'{\'e}tende {\`a} une perturbation \( \mathcal{N}_{\phi } \) de \( \mathcal{N} \)).
Supposons que \( \mathcal{T}_{\phi } \) admette une structure de contact holomorphe,
proche de la structure standard \( \eta ^{c} \). Alors celle-ci s'{\'e}tend en
une structure de contact holomorphe sur \( \mathcal{N}_{_{\phi }} \).
\end{thm}
\begin{proof}
D'apr{\`e}s le corollaire (\ref{cor-31}), on a une estimation 
\[
\left\Vert \overline{\partial }\xi _{<0}\right\Vert _{L^{2}}\geq c(\left\Vert \xi _{<0}\right\Vert _{L^{2}}+\left\Vert \nabla _{H}\xi _{<0}\right\Vert _{L^{2}});\]
 compte tenu de l'estimation (\ref{est-xi}), si \( \phi  \) et \( \xi  \)
sont assez petits, on d{\'e}duit \( \xi _{<0}=0 \), donc \( \xi  \) est {\`a} coefficients
de Fourier positifs ; on peut alors prolonger \( \xi  \) de mani{\`e}re holomorphe
le long de chaque disque : le prolongement obtenu v{\'e}rifie encore l'{\'e}quation
(\ref{etac+xi-hol}), ce qui signifie qu'on a finalement obtenu dans \( \widetilde{\mathcal{N}}_{\phi } \)
une 1-forme \( \varpi  \) {\`a} valeurs dans \( L \), solution de l'{\'e}quation 
\[
\overline{\partial }\varpi =\alpha \otimes \varpi ,\quad \alpha \in \Omega ^{0,1}_{\phi },\]
 et, quitte {\`a} int{\'e}grer \( \alpha  \) dans l'op{\'e}rateur \( \overline{\partial } \)
du fibr{\'e} \( L \), une solution de 
\[
\overline{\partial }\varpi =0.\]
 Rappelons (voir le lemme \ref{lem-fill-J}) que la base du fibr{\'e} en disques
est \( P^{1}\times P^{1} \) avec fibr{\'e} normal \( \mathcal{O}(-1,1) \). La
forme \( \varpi  \) d{\'e}finit au-dessus de \( P^{1}\times P^{1} \) une section
holomorphe de \( \Omega ^{1}_{\widetilde{\mathcal{N}}_{\phi }}\otimes L \),
avec \( L \) ne pouvant {\^e}tre {\'e}gal qu'{\`a} \( \mathcal{O}(0,2) \), mais, comme
dans la d{\'e}monstration du lemme \ref{lem-31}, \( \varpi  \) doit s'annuler
sur \( T^{1,0}P^{1}_{1} \), et descend {\`a} la contraction \( \mathcal{N}_{\phi } \)
de \( P_{1}^{1} \). On obtient ainsi une distribution holomorphe sur \( \mathcal{N}_{\phi } \),
qui est bien de contact car c'est une petite d{\'e}formation de la structure de
contact de \( \mathcal{N} \).
\end{proof}
\begin{cor}
\label{cor-fill-g}Si on a une petite d{\'e}formation \( g \) de la m{\'e}trique standard
de \( S^{3} \), dont l'espace des twisteurs satisfait les hypoth{\`e}ses du lemme
\ref{lem-fill-J}, alors \( g \) est le bord {\`a} l'infini d'une m{\'e}trique autoduale
Einstein (compl{\`e}te) sur \( B^{4} \).
\end{cor}
\begin{proof}
Le corollaire (\ref{cor-fill-g2}) nous dit que \( g \) est le bord d'une m{\'e}trique
autoduale sur \( B^{4} \), et le th{\'e}or{\`e}me pr{\'e}c{\'e}dent indique que son espace
des twisteurs \( \mathcal{N} \) admet une forme de contact holomorphe ; les
\( P^{1} \) r{\'e}els sont tangents {\`a} la distribution de contact holomorphe au
bord, transverses {\`a} l'int{\'e}rieur, ce qui fournit la m{\'e}trique d'Einstein voulue
dans la classe conforme.
\end{proof}

\section{\label{sec4}Construction de disques extr{\'e}maux}

Dans les sections \ref{sec2} et \ref{sec3}, on a vu une proc{\'e}dure pour fabriquer
des m{\'e}triques autoduales, ou autoduales Einstein, {\`a} partir d'une donn{\'e}e conforme
\( [g,\mathcal{Q}] \) sur le bord \( S^{3} \) : la m{\'e}thode consiste {\`a} {\'e}tendre
les structures le long des disques holomorphes, pourvu que la structure complexe
de l'espace des twisteurs au bord puisse se pr{\'e}senter avec des coefficients
de Fourier positifs. Dans cette section, on montre la r{\'e}ciproque, {\`a} savoir que
si l'extension existe, alors de tels disques holomorphes doivent toujours exister.

\subsection{Existence d'une fibration en disques}

Dans le cas d'un domaine pseudo-convexe, on dispose de m{\'e}thodes pour produire
des disques extr{\'e}maux \cite{Lem:81}. Dans le cas g{\'e}n{\'e}ral o{\`u} la forme de Levi
n'a pas de signe fix{\'e}, mais demeure non d{\'e}g{\'e}n{\'e}r{\'e}e, nous construisons, par d{\'e}formation,
des disques extr{\'e}maux au sens de Lempert \cite[lemme 7.1, th{\'e}or{\`e}me 10.1]{Lem:92}.

\begin{thm}
\label{th-disques}Soit \( D \) un fibr{\'e} en disques holomorphes au-dessus de
la vari{\'e}t{\'e} complexe \( X \), tel que la forme de Levi de \( \partial D \)
soit non d{\'e}g{\'e}n{\'e}r{\'e}e. Soit \( D_{\phi } \) une petite perturbation de la structure
complexe de \( D \), telle que \( X\subset D_{\phi } \) demeure une sous-vari{\'e}t{\'e}
complexe. Alors, pour tout point \( p\in \partial D_{\phi } \), il existe un
unique plongement \( f:\overline{\Delta }\rightarrow \overline{D_{\phi }} \),
et un unique sous-fibr{\'e} \( E \) de \( f^{*}T^{1,0}\overline{D_{\phi }} \),
tels que :
\begin{lyxlist}{00.00.0000}
\item [(a)]\( f \) et \( E \) sont holomorphes sur \( \Delta  \) ;
\item [(b)]\( f(1)=p \), \( f(S^{1})\subset \partial D \), \( f(0)\in X \) ;
\item [(c)]\( E|_{S^{1}}=T^{1,0}\partial D_{\phi } \), et \( E_{0} \) est tangent
{\`a} \( X \) ;
\end{lyxlist}
De plus, les \( f=f^{p} \) satisfont
\begin{lyxlist}{00.00.0000}
\item [(d)]\( T^{1,0}\Delta  \) et \( E \) sont transverses.
\item [(e)]\( f^{up}(\zeta )=f^{p}(u\zeta ) \) pour tout \( u\in S^{1} \) ;
\item [(f)]\( F:\partial D\times \overline{\Delta }\rightarrow \overline{D_{\phi }} \)
d{\'e}finie par \( F(p,\zeta )=f^{p}(\zeta ) \) est lisse, et \( F(\cdot ,\zeta ):\partial D\rightarrow \overline{D_{\phi }} \)
est une immersion pour \( \zeta \neq 0 \). 
\end{lyxlist}
\end{thm}
Nous donnerons la d{\'e}monstration du th{\'e}or{\`e}me dans la section \ref{sec42}. Pour
le moment, nous en tirons quelques cons{\'e}quences. Remarquons que la propri{\'e}t{\'e}
(e) dans le th{\'e}or{\`e}me indique qu'on a en r{\'e}alit{\'e} construit une action de \( S^{1} \)
(et m{\^e}me de \( \overline{\Delta } \)) sur \( D_{\phi } \). Le corollaire suivant
indique que la construction faite dans le lemme \ref{lem-fill-J} capture bien
toutes les d{\'e}formations de \( D \).

\begin{cor}
\label{cor-ex-disques}Sous les hypoth{\`e}ses du th{\'e}or{\`e}me \ref{th-disques}, l'action
de \( S^{1} \) d{\'e}finie par les \( f^{p} \) sur \( \partial D_{\phi } \) est
de contact ; plus g{\'e}n{\'e}ralement, la distribution \( \Re (E+\overline{E})\subset TD_{\phi } \)
est bien d{\'e}finie, et invariante sous \( \overline{\Delta } \), ce qui signifie
que \( D_{\phi } \) se pr{\'e}sente par rapport {\`a} \( D \) par un tenseur \( \phi  \)
satisfaisant les hypoth{\`e}ses du lemme \ref{lem-fill-J}.
\end{cor}
\begin{proof}
Ce corollaire est essentiellement \cite[Lemme 7.1]{Lem:92} et nous ne refaisons
pas la d{\'e}monstration : le point crucial est de montrer que si \( Y\in T_{p}\partial D_{_{\phi }} \)est
un vecteur dans la distribution de contact, alors pour tout \( \zeta \in \overline{\Delta } \),
on a 
\[
F(\cdot ,\zeta )_{*}Y\in E+\overline{E},\]
 ce qui est une cons{\'e}quence sans difficult{\'e} des conditions au bord et du fait
que \( E \) est holomorphe le long des disques.
\end{proof}
Nous d{\'e}duisons maintenant de ce r{\'e}sultat une caract{\'e}risation des \( [g,\mathcal{Q}] \)
qui admettent un remplissage autodual.

\begin{thm}
\label{th-fill-g2}Soit \( [g,\mathcal{Q}] \) la donn{\'e}e d'une premi{\`e}re et d'une
seconde formes fondamentales sur \( S^{3} \), {\`a} changement conforme pr{\`e}s. Supposons
que \( [g,\mathcal{Q}] \) soit suffisamment proche de la m{\'e}trique standard,
alors \( [g,\mathcal{Q}] \) est le bord d'une m{\'e}trique autoduale sur la boule
\( B^{4} \) si et seulement si l'espace des twisteurs \( \mathcal{T} \) de
\( [g,\mathcal{Q}] \) admet une action de contact de \( S^{1} \), anti-commutant
{\`a} la structure r{\'e}elle de \( \mathcal{T} \), par rapport {\`a} laquelle la structure
CR de \( \mathcal{T} \) soit {\`a} coefficients de Fourier positifs ou nuls.
\end{thm}
Nous pouvons aussi donner un r{\'e}sultat sur le remplissage d'une classe conforme
\( [g] \) par une m{\'e}trique d'Einstein autoduale.

\begin{thm}
\label{th-fill-g}Soit \( [g] \) une m{\'e}trique conforme sur \( S^{3} \), suffisamment
proche de la m{\'e}trique standard, alors \( [g] \) est l'infini conforme d'une
m{\'e}trique d'Einstein autoduale sur la boule \( B^{4} \) si et seulement si l'espace
des twisteurs \( \mathcal{T} \) de \( [g] \) admet une action de contact de
\( S^{1} \), anti-commutant {\`a} la structure r{\'e}elle de \( \mathcal{T} \), par
rapport {\`a} laquelle la structure CR de \( \mathcal{T} \) soit {\`a} coefficients
de Fourier positifs ou nuls.
\end{thm}
\noindent \emph{D{\'e}monstration des th{\'e}or{\`e}mes \ref{th-fill-g2} et \ref{th-fill-g}.}
Si l'on suppose que \( [g,\mathcal{Q}] \) admet un remplissage autodual, alors
l'espace des twisteurs \( \mathcal{T} \) admet un remplissage par l'espace
des twisteurs \( \mathcal{N} \) de la m{\'e}trique autoduale. Au-dessus du point
\( 0\in B^{4} \), nous avons une fibre twistorielle \( \mathcal{N}_{0} \)
qui est un \( P^{1} \) ; en {\'e}clatant cette fibre, on obtient une vari{\'e}t{\'e} \( \widetilde{\mathcal{N}} \)
avec au-dessus de \( \mathcal{N}_{0} \) un \( P^{1}\times P^{1}\subset \widetilde{\mathcal{N}} \),
et on garde \( \mathcal{T}=\partial \widetilde{\mathcal{N}} \) ; on peut appliquer
le th{\'e}or{\`e}me \ref{th-disques} et le corollaire \ref{cor-ex-disques} avec \( X=P^{1}\times P^{1} \)
pour d{\'e}duire que \( \mathcal{T} \) admet effectivement l'action de \( S^{1} \)
demand{\'e}e. D'autre part, l'assertion d'unicit{\'e} dans le th{\'e}or{\`e}me \ref{th-disques}
contraint l'action de \( S^{1} \) qui s'en d{\'e}duit {\`a} satisfaire la condition
de compatibilit{\'e} (\ref{tau-S1}) avec la structure r{\'e}elle.

R{\'e}ciproquement, {\'e}tant donn{\'e} \( [g,\mathcal{Q}] \) satisfaisant les hypoth{\`e}ses
du th{\'e}or{\`e}me, l'existence du remplissage autoduale r{\'e}sulte du corollaire \ref{cor-fill-g2}
; si de plus \( \mathcal{Q}=0 \), alors le prolongement autodual admet bien
dans sa classe conforme un repr{\'e}sentant d'Einstein par le corollaire \ref{cor-fill-g}.\qed

\subsection{\label{sec42}D{\'e}monstration du th{\'e}or{\`e}me \ref{th-disques}}

Il est clair que \( D \), avec sa structure de fibr{\'e} en disques \( p:D\rightarrow X \),
satisfait les conclusions du th{\'e}or{\`e}me. On va construire les structures pour
\( D_{\phi } \) par un argument de d{\'e}formation.

Puisque la forme de Levi du bord est non d{\'e}g{\'e}n{\'e}r{\'e}e, on a une 1-forme de connexion
\( \eta  \) sur \( \partial D, \) se restreignant {\`a} la forme angulaire \( d\theta  \)
sur chaque disque, telle que 
\[
(\ker \eta )\otimes \mathbf{C}=T^{0,1}\partial D+T^{1,0}\partial D\]
 soit une distribution de contact sur \( \partial D \) ; comme deux structures
de contact proches sont diff{\'e}omorphes, on peut supposer que la m{\^e}me chose demeure
valable pour \( D_{\phi } \), ce qui se traduit sur \( \phi  \) par la condition
\[
\phi (T^{0,1}\partial D)\subset T^{1,0}\partial D.\]

Le r{\'e}sultat {\'e}tant local au voisinage d'un disque donn{\'e} au-dessus d'un point
\( x\in X \), on va choisir des coordonn{\'e}es locales \( (z^{i})_{i\geq 1} \)
sur \( X \) pr{\`e}s de \( x \). Nous trivialisons le fibr{\'e} \( \Omega ^{1,0}D \)
pr{\`e}s du disque \( p^{-1}(x) \) par la base 
\[
e^{0}=2i\eta ^{1,0},\quad e^{i}=p^{*}dz^{i};\]
 il faut noter que sur chaque disque, \( e^{0} \) s'identifie simplement {\`a}
\( d\zeta /\zeta  \) ; en particulier, nous avons choisi \( e^{0} \) singulier
sur la section nulle ; nous avons la base duale constitu{\'e}e du vecteur d'homoth{\'e}tie
dans les fibres, \( \partial _{0}=\zeta \partial /\partial \zeta  \), et des
relev{\'e}s horizontaux \( \partial _{i} \) des vecteurs \( \partial /\partial z^{i} \)
de la base.

{\'E}tant donn{\'e}e une perturbation \( \phi  \) de la structure complexe de \( D \),
nous trivialisons le fibr{\'e} \( \Omega ^{1,0}D_{\phi } \) par les \( e^{i}-\phi ^{t}e^{i} \),
et nous cherchons \( (f,\varpi ) \) avec \( f:\overline{\Delta }\rightarrow \overline{D} \)
et 
\[
\varpi =\sum _{i\geq 0}\varpi _{i}(e^{i}-\phi ^{t}e^{i}),\]
 o{\`u} les \( \varpi _{i} \) sont des fonctions sur \( \overline{\Delta } \)
et \( E=\ker \varpi  \) ; bien entendu, \( \varpi  \) n'est d{\'e}fini qu'{\`a} multiplication
pr{\`e}s par une fonction, et nous choisissons la normalisation 
\[
\varpi _{0}=\zeta ;\]
 compte tenu de la singularit{\'e} de \( e^{0} \) sur la section nulle, ce choix
d{\'e}finit bien une forme lisse y compris sur la section nulle ; nous notons {\'e}galement
\[
\pi ^{\varpi }(\sum _{i\geq 0}\vartheta _{i}e^{i})=\sum _{i>0}(\vartheta _{i}-\vartheta _{0}\varpi _{i})e^{i}\]
 la projection sur \( \Omega ^{1,0}D/\mathbf{C}\varpi \approx p^{*}\Omega ^{1,0}X \).

Le couple \( (f,\varpi ) \) doit satisfaire aux {\'e}quations suivantes : d'une
part, la diff{\'e}rentielle \( df:T\overline{\Delta }\rightarrow T\overline{D}_{\phi } \)
se d{\'e}compose en \( \overline{\partial }_{\phi }f+\partial _{\phi }f:T^{0,1}\overline{\Delta }\rightarrow T^{1,0}\overline{D}_{\phi }+T^{0,1}\overline{D}_{\phi } \),
et \( f \) est holomorphe si 
\[
\overline{\partial }_{\phi }f=0;\]
 comme on peut trivialiser \( T^{0,1}\overline{\Delta } \) par le vecteur \( \partial /\partial \overline{\zeta } \),
on peut consid{\'e}rer \( \overline{\partial }_{\phi }f \) comme une section de
\( f^{*}T^{1,0}\overline{D}_{\phi }\approx f^{*}T^{1,0}\overline{D} \) via
\( \phi  \) ; d'autre part, \( E \) doit {\^e}tre holomorphe le long de \( \overline{\Delta } \),
ce qui se traduit par \( \overline{\partial }_{\phi }\varpi =\alpha \otimes \varpi  \)
; en identifiant {\`a} nouveau \( \Omega ^{1,0}\overline{D}_{\phi } \) et \( \Omega ^{1,0}\overline{D} \)
via \( \phi  \), on obtient finalement l'{\'e}quation 
\[
\pi ^{\varpi }\overline{\partial }_{\phi }\varpi =0,\]
 o{\`u}, {\`a} nouveau en trivialisant \( T^{0,1}\overline{\Delta } \) par \( \partial /\partial \overline{\zeta } \),
on peut consid{\'e}rer \( \pi ^{\varpi }\overline{\partial }_{\phi }\varpi  \)
comme une section de \( f^{*}\Omega ^{1,0}X \). 

Les conditions au bord se traduisent par 
\begin{eqnarray}
f(1) & = & p,\nonumber \\
f(S^{1}) & \subset  & \partial D,\label{bord-f} \\
f(0) & \in  & X,\nonumber 
\end{eqnarray}
 et, pour \( i>0 \), 
\begin{equation}
\label{bord-varpi}
\begin{array}{rl}
\varpi _{i}|_{S^{1}} & =0,\\
\varpi _{i}(0) & =0.
\end{array}
\end{equation}

Pour bien poser le probl{\`e}me d'inversion locale, il nous reste {\`a} poser 
\begin{eqnarray*}
Map(\overline{\Delta },\overline{D}) & = & \{f:\overline{\Delta }\rightarrow \overline{D},\textrm{ de classe }C^{1+\alpha },\textrm{ satisfaisant }(\ref {bord-f})\},\\
\mathcal{E}(\overline{\Delta },\overline{D}) & = & \{\varpi =\zeta e^{0}+\sum _{i>0}\varpi _{i}e^{i},\textrm{ }\varpi _{i}\textrm{ de classe C}^{1+\alpha },\textrm{ satisfaisant }(\ref {bord-varpi})\}.
\end{eqnarray*}
 Nous regardons l'op{\'e}rateur 
\[
P(f,\varpi ,\phi )=(\overline{\partial }_{\phi }f,\pi ^{\varpi }\overline{\partial }_{\phi }\varpi ),\]
 d{\'e}fini sur \( Map(\overline{\Delta },\overline{D})\times \mathcal{E}(\overline{\Delta },\overline{D}) \),
{\`a} valeurs dans les espaces de sections (rappelons que nous avons trivialis{\'e}
\( T^{1,0}\overline{D} \) et \( \Omega ^{1,0}X \) dans un voisinage du disque
consid{\'e}r{\'e}) 
\[
C^{\alpha }(\overline{\Delta },T^{1,0}\overline{D})\times C^{\alpha }(\overline{\Delta },\Omega ^{1,0}X).\]
 Il est manifeste que ces espaces sont des vari{\'e}t{\'e}s de Banach, et l'op{\'e}rateur
\( P \) est lisse entre ces espaces (on ne regarde que des \( \phi  \) lisses,
mais on peut, pour avoir un espace de Banach, se contenter de \( \phi  \) de
classe \( C^{k} \) pour \( k \) assez grand).

\begin{lem}
La diff{\'e}rentielle en \( (f(\zeta )=\zeta p,\varpi =\zeta e^{0},\phi =0) \)
de \( P \) par rapport {\`a} \( (f,\varpi ) \) est un isomorphisme. 
\end{lem}
Admettons le lemme quelques instants. On en d{\'e}duit que si la perturbation \( \phi  \)
est assez petite, alors le probl{\`e}me \( P(f,\varpi ,\phi )=0 \) a une solution
unique dans \( Map(\overline{\Delta },\overline{D})\times \mathcal{E}(\overline{\Delta },\overline{D}) \)
; cela fournit le couple \( (f,E) \) du th{\'e}or{\`e}me \ref{th-disques} satisfaisant
les propri{\'e}t{\'e}s (a) {\`a} (c). La propri{\'e}t{\'e} (d) reste vraie apr{\`e}s petite perturbation.
Analysons {\`a} pr{\'e}sent la d{\'e}pendance par rapport au point \( p \) : gr{\^a}ce {\`a} l'unicit{\'e}
de la solution, la propri{\'e}t{\'e} (e) est imm{\'e}diate, tandis que (f) r{\'e}sulte du fait
que la solution d{\'e}pend de mani{\`e}re lisse des param{\`e}tres. Le th{\'e}or{\`e}me \ref{th-disques}
est d{\'e}montr{\'e}.\qed

\begin{rem}
La condition que \( X \) demeure une sous-vari{\'e}t{\'e} complexe a pour seule utilit{\'e}
de donner un sens {\`a} la condition (c). Mais la r{\'e}solution ci-dessus reste valable
sans cette condition, la diff{\'e}rence {\'e}tant que la condition \( \varpi _{i}(0)=0 \)
n'a plus d'interpr{\'e}tation g{\'e}om{\'e}trique naturelle.
\end{rem}
\begin{proof}
Il nous reste {\`a} calculer la diff{\'e}rentielle de l'op{\'e}rateur \( P \). Un vecteur
tangent {\`a} \( Map(\overline{\Delta },\overline{D}) \) est un champ de vecteur
\( X \) le long de \( \overline{\Delta } \), mais il sera plus commode d'utiliser
\[
X=f+\overline{f},\]
 avec \( f \) vecteur\footnote{
il est commode ici d'utiliser le m{\^e}me symbole pour l'application \( f \) et
pour sa version infinit{\'e}simale qui est un champ de vecteurs
} de type (1,0), repr{\'e}sent{\'e} par 
\[
f=f^{0}\frac{e_{0}}{\zeta }+\sum _{i>0}f^{i}e_{i},\quad f^{0},f^{1},\ldots \in C^{1+\alpha };\]
 le facteur \( \zeta  \) s'explique en se rappelant que \( e_{0} \) s'annule
sur \( X \), mais \( e_{0}/\zeta  \) est lisse. Les conditions au bord (\ref{bord-f})
se traduisent par 
\begin{eqnarray}
\Re \frac{f^{0}}{\zeta } & = & 0\textrm{ sur le bord }S^{1},\nonumber \\
f^{i}(1) & = & 0,\label{bord-fi} \\
f^{0}(0) & = & 0.\nonumber 
\end{eqnarray}
 D'autre part, rappelons que \( \varpi =\zeta e^{0}+\sum _{i>0}\varpi _{i}e^{i} \),
avec \( \varpi _{i} \) de classe \( C^{1+\alpha } \), et conditions au bord
(\ref{bord-varpi}).

En \( \phi =0 \), l'op{\'e}rateur \( P \) est simplement \( P(f,\varpi )=(\overline{\partial }f,\pi ^{\varpi }\overline{\partial }\varpi ) \)
; plus pr{\'e}cis{\'e}ment, le second op{\'e}rateur est 
\[
\pi ^{\varpi }\iota _{f_{*}\frac{\partial }{\partial \overline{\zeta }}}d\varpi ;\]
 compte tenu que la courbure du fibr{\'e} ---une 2-forme non d{\'e}g{\'e}n{\'e}r{\'e}e sur \( X \)---,
est 
\[
de^{0}=2i\overline{\partial }\eta ^{1,0}=F=F_{\overline{\imath }j}e^{\overline{\imath }}\wedge e^{j},\]
 la diff{\'e}rentielle de \( P \) suivant \( \varpi  \) s'{\'e}crit 
\[
\pi ^{\zeta e^{0}}(\iota _{[X,\frac{\partial }{\partial \overline{\zeta }}]}d(\zeta e^{0})+\iota _{\frac{\partial }{\partial \overline{\zeta }}}d\varpi )=\pi ^{\zeta e^{0}}(\iota _{[f+\overline{f},\frac{\partial }{\partial \overline{\zeta }}]}\zeta F+\iota _{\frac{\partial }{\partial \overline{\zeta }}}d\varpi _{i}\wedge e^{i}).\]
 Nous sommes maintenant en mesure d'{\'e}crire explicitement le probl{\`e}me infinit{\'e}simal
: {\'e}tant donn{\'e}s le long de \( \overline{\Delta } \) un (1,0)-vecteur \( g=g^{0}\frac{e_{0}}{\zeta }+\sum _{i>0}g^{i}e_{i} \),
et une (1,0)-forme \( \rho =\sum _{i>0}\rho _{i}e^{i}, \) les {\'e}quations {\`a} r{\'e}soudre
deviennent 
\begin{eqnarray}
\frac{\partial f^{i}}{\partial \overline{\zeta }} & = & g^{i},\quad i\geq 0,\label{eq-fi} \\
\frac{\partial \varpi _{i}}{\partial \overline{\zeta }}-\zeta \frac{\partial \overline{f^{j}}}{\partial \overline{\zeta }}F_{\overline{\jmath }i} & = & \rho _{i},\quad i>0,\label{eq-varpii} 
\end{eqnarray}
 avec les conditions au bord (\ref{bord-fi}) et (\ref{bord-varpi}).

L'{\'e}quation (\ref{eq-fi}) d{\'e}finit les \( f^{i} \) {\`a} fonction holomorphe pr{\`e}s
sur \( \overline{\Delta } \). 

Pour ce qui est de \( f^{0} \), les conditions au bord (\ref{bord-fi}) la
d{\'e}finissent uniquement : en effet, commen{\c c}ons par choisir la solution \( \varphi _{1}^{0} \)
de 
\[
\frac{\partial \varphi ^{0}_{1}}{\partial \overline{\zeta }}=g^{0}\]
 telle qu'au bord seuls les coefficients de Fourier strictement n{\'e}gatifs soient
non nuls : 
\[
\varphi ^{0}_{1}|_{S^{1}}=\sum _{k<0}a_{k}e^{ik\theta };\]
 pour obtenir la condition \( \Re (f^{0}/\zeta )=0 \) au bord, nous rajoutons
{\`a} \( \varphi ^{0}_{1} \) la fonction holomorphe 
\[
\varphi ^{0}_{2}=-\sum _{k>2}\overline{a_{-k+2}}\zeta ^{k};\]
 nous avons encore la libert{\'e} de rajouter {\`a} \( f^{0}=\varphi ^{0}_{1}+\varphi ^{0}_{2} \)
un terme 
\[
\varphi _{3}=a_{0}+a_{1}\zeta -\overline{a_{0}}\zeta ^{2},\quad \Re a_{1}=0,\]
 mais les coefficients \( a_{0} \) et \( a_{1} \) sont maintenant impos{\'e}s
par les conditions \( f^{0}(0)=f^{0}(1)=0 \). Bien entendu \( f^{0}=\varphi ^{0}_{1}+\varphi ^{0}_{2}+\varphi ^{0}_{3} \)
est de classe \( C^{1+\alpha } \) si \( g^{0} \) est de classe \( C^{\alpha } \).

Passons maintenant aux \( \varpi _{i} \) : un choix des \( f_{i} \) ayant
{\'e}t{\'e} fait, on peut r{\'e}soudre (\ref{eq-varpii}) de sorte que les coefficients
de Fourier positifs ou nuls des \( \varpi _{i} \) au bord soient nuls, mais
cette condition est encore loin de la condition au bord (\ref{bord-varpi})
voulue, {\`a} savoir l'annulation compl{\`e}te des \( \varpi _{i} \) au bord ; examinons
l'influence d'une modification des \( f^{i} \) (\( i>0 \)) par des fonctions
holomorphes \( \varphi ^{i} \) : 
\[
f^{i}\rightarrow f^{i}+\varphi ^{i};\]
 en s{\'e}parant
\[
\varphi ^{i}=a_{0}^{i}+a_{1}^{i}\zeta +\varphi _{1}^{i},\quad \varphi _{1}^{i}=\sum _{k>1}a_{k}^{i}\zeta ^{k},\]
 on observe que les \( \varpi _{i} \) se transforment par la loi 
\[
\varpi _{i}\rightarrow \varpi _{i}+\zeta \overline{\varphi _{1}^{j}}F_{\overline{\jmath }i}+(|\zeta |^{2}-1)\overline{a_{1}^{j}}F_{\overline{\jmath }i};\]
 par cons{\'e}quent, en choisissant \( \varphi ^{j}_{_{1}} \) holomorphe de sorte
que (ici nous utilisons le fait que la forme de Levi est non d{\'e}g{\'e}n{\'e}r{\'e}e, donc
\( (F_{\overline{\imath }j}) \) admet un inverse \( (F^{\overline{\imath }j}) \)
) 
\[
\varphi _{1}^{j}|_{S^{1}}=-\zeta F^{\overline{\imath }j}\overline{\varpi _{i}}|_{S^{1}},\]
 ce qui est possible puisque les \( \varpi _{i}|_{S^{1}} \) sont {\`a} coefficients
de Fourier strictement n{\'e}gatifs, on peut obtenir la condition au bord \( \varpi |_{S^{1}}=0 \)
; la condition \( \varpi _{i}(0)=0 \) est obtenue en fixant les coefficients
\( a_{1}^{i} \), et la seule libert{\'e} demeurant sur les \( f^{i} \) est d'ajouter
une constante \( a_{0}^{i} \), ce qui permet d'assurer la condition au bord
\( f^{i}(1)=0 \). De nouveau, la solution \( (f^{i},\varpi _{i}) \) est de
classe \( C^{1+\alpha } \) si les \( g^{i} \) et les \( \rho _{i} \) sont
de classe \( C^{\alpha } \), ce qui ach{\`e}ve la d{\'e}monstration du lemme.
\end{proof}

\subsection{\label{sec43}Les disques holomorphes et la g{\'e}om{\'e}trie conforme}

Il est int{\'e}ressant de noter que les disques holomorphes dans l'espace des twisteurs
d'une vari{\'e}t{\'e} antiautoduale ont une interpr{\'e}tation en termes de disques dans
la vari{\'e}t{\'e} conforme. Si l'on dispose d'une m{\'e}trique conforme \( [g] \) sur
la vari{\'e}t{\'e} \( X^{4} \), alors la seconde forme fondamentale \( B \) d'une
surface se d{\'e}compose en trois composantes : 
\begin{eqnarray*}
B & = & B^{+}+B^{-}\\
 & = & B^{+}+B^{2,0}+B^{0,2};
\end{eqnarray*}
 en effet, le plan tangent \( P\subset TX \) {\`a} la surface d{\'e}termine une structure
presque-complexe hermitienne \( J \) sur \( TX=P\oplus P^{\perp } \); alors
\( B^{+}=B^{1,1} \) est la partie \( J \)-invariante de \( B \), satisfaisant
\[
B^{+}_{JU,JV}=B_{U,V}^{+},\quad U,V\in P,\]
et \( B^{-} \) est la partie \( J \)-anti-invariante, se d{\'e}composant elle-m{\^e}me
en deux composantes \( B^{2,0} \) et \( B^{0,2} \) satisfaisant respectivement
\[
B^{2,0}_{JU,V}=J^{P^{\perp }}B_{U,V},\quad B^{0,2}_{JU,V}=-J^{P^{\perp }}B_{U,V}.\]

Une surface \( \Sigma  \) dans \( X \) d{\'e}termine, comme on l'a vu, une structure
presque-complexe \( J \) sur \( TX \) le long de \( \Sigma  \), c'est-{\`a}-dire
une section de la fibration twistorielle que l'on peut appeler le \emph{rel{\`e}vement
de \( \Sigma  \) {\`a} l'espace des twisteurs.} L'{\'e}nonc{\'e} suivant est essentiellement
contenu dans l'article \cite{Gau:86} de Gauduchon.

\begin{lem}
\label{lem-tw/conf}Supposons la m{\'e}trique \( [g] \) sur \( X \) antiautoduale,
soit \( p:\mathcal{N}\rightarrow X \) son espace de twisteurs, et \( \Sigma  \)
une surface de Riemann ; le rel{\`e}vement {\`a} l'espace des twisteurs donne une correspondance
entre
\begin{enumerate}
\item les immersions conformes \( \Sigma \rightarrow X \), dont la seconde forme
fondamentale satisfait \( B^{0,2}=0 \) ;
\item les immersions holomorphes \( \Sigma \rightarrow \mathcal{N} \), transverses
{\`a} la fibration en \( P^{1} \).
\end{enumerate}
\end{lem}
Notons que l'{\'e}nonc{\'e} est donn{\'e} dans le cas antiautodual, qui permet de repr{\'e}senter
plus agr{\'e}ablement l'espace des twisteurs comme le fibr{\'e} des structures presque-complexes
hermitiennes, compatibles {\`a} l'orientation ; le cas autodual est obtenu par un
simple changement d'orientation.

\begin{proof}
La correspondance entre la condition \( B^{0,2}=0 \) et l'holomorphie du rel{\`e}vement
{\`a} l'espace des twisteurs est montr{\'e}e dans \cite{Gau:86}. Il reste {\`a} montrer
qu'une immersion holomorphe \( \iota :\Sigma \rightarrow \mathcal{N} \) provient
bien du rel{\`e}vement de \( p\circ \iota  \). Notons \( J \) la structure presque-complexe
sur \( TX \) le long de \( p\iota (\Sigma ) \), induite par la section \( \iota  \)
de l'espace des twisteurs. Soit un point \( \sigma \in \Sigma  \) : par la
d{\'e}finition m{\^e}me de la structure complexe de l'espace des twisteurs, la projection
\( p:T_{\iota (\sigma )}\mathcal{N}\rightarrow (T_{p\iota (\sigma )}X,J) \)
est \( \mathbf{C} \)-lin{\'e}aire, donc, puisque \( \iota  \) est holomorphe,
\( p\circ \iota :T_{\sigma }\Sigma \rightarrow (T_{p\iota (\sigma )}X,J) \)
est aussi \( \mathbf{C} \)-lin{\'e}aire, ce qui impose que la structure complexe
le long de \( p\iota (\Sigma ) \) induite par l'immersion \( p\iota  \) est
aussi \( J \) ; cela montre que \( \iota  \) est bien le rel{\`e}vement {\`a} l'espace
des twisteurs de \( p\circ \iota  \).
\end{proof}

\section{\label{sec5}Le groupe des contactomorphismes}

Il est bien connu qu'un contactomorphisme est essentiellement donn{\'e} par une
fonction r{\'e}elle ; si l'on dispose d'un champ de Reeb \( R \) pour la forme
de contact \( \eta  \), le contactomorphisme infinit{\'e}simal correspondant {\`a}
la fonction r{\'e}elle \( f \) n'est autre que 
\[
X_{f}=fR-\#d_{H}f,\]
 o{\`u} \( \# \) est d{\'e}fini sur les 1-formes horizontales par \( d\eta (\#\alpha ,X)=\langle \alpha ,X\rangle  \)
; en particulier, \( X_{f} \) ne d{\'e}pend que de \( f \) et de \( d_{H}f \).

Il est moins facile de produire une param{\'e}trisation \( f\rightarrow \varphi _{f} \)
des contactomorphismes proches de l'identit{\'e}, de telle mani{\`e}re que \( \varphi _{f} \)
ne d{\'e}pende que de \( f \) et des d{\'e}riv{\'e}es \emph{horizontales} de \( f \).
Par exemple, la param{\'e}trisation naturelle provenant de la symplectification
semble utiliser toutes les d{\'e}riv{\'e}es de \( f \), et la perte de d{\'e}riv{\'e}e qui
en d{\'e}coule contraint ensuite {\`a} utiliser des arguments de type Nash-Moser, voir
\cite{Che-Lee:95}.

En dimension 3, Bland \cite{Bla:94} a propos{\'e} dans le cas de \( S^{3} \) une
param{\'e}trisation du groupe des contactomorphismes qui ne consomme que les d{\'e}riv{\'e}es
horizontales.

Dans cette section, nous proposons une autre param{\'e}trisation, valable en dimension
sup{\'e}rieure ou {\'e}gale {\`a} 5, qui conviendra pour nos applications.

Nous supposons que la vari{\'e}t{\'e} de contact est un \( S^{1} \)-fibr{\'e} au-dessus
d'une vari{\'e}t{\'e} symplectique \( (X,\omega ) \) : 
\[
\begin{array}{ccc}
S^{1} & \rightarrow  & M\\
 &  & \; \downarrow \pi \\
 &  & X
\end{array}\]
 et que la forme de contact \( \eta  \) est une 1-forme de connexion, avec
\( d\eta =\pi ^{*}\omega  \). On notera \( H=\ker \eta  \) la distribution
de contact. L'id{\'e}e de base de la construction qui suit est la simple constatation
que, pour regarder les d{\'e}riv{\'e}es horizontales d'un contactomorphisme \( \varphi  \)
de \( M \), il suffit de regarder les d{\'e}riv{\'e}es horizontales de la composition
\( \pi \circ \varphi :M\rightarrow X \).

Dans notre situation, il est souhaitable d'utiliser plut{\^o}t les espaces fonctionnels
de Folland-Stein \cite{Fol-Ste:74}, c'est-{\`a}-dire des espaces de Sobolev anisotropes,
mieux adapt{\'e}s {\`a} la g{\'e}om{\'e}trie de contact : l'espace \( \mathcal{H}^{k} \) consiste
des fonctions avec \( k \) d{\'e}riv{\{\'e}es suivant les directions de contact (et,
par cons{\'e}quent, \( \frac{k}{2} \) d{\'e}riv{\'e}es dans la direction du champ de Reeb).
Nous choisirons toujours \( k \) suffisamment grand, de sorte qu'on ait l'injection
de Sobolev pour les espaces de Folland-Stein, \( \mathcal{H}^{k}\subset C^{0} \);
cela implique que les espaces de fonctions de classe \( \mathcal{H}^{k} \)
sont des alg{\`e}bres, ce qui permet de faire l'analyse non lin{\'e}aire. Plus pr{\'e}cis{\'e}ment,
on a le lemme suivant :

\begin{lem}
\label{lem-FS-alg}Pour \( k>n \), l'espace des fonctions de r{\'e}gularit{\'e} \( \mathcal{H}^{k} \)
sur la vari{\'e}t{\'e} de contact \( X^{n} \) est une alg{\`e}bre, et agit par multiplication
sur les espaces \( \mathcal{H}^{i} \) pour \( i<k \).
\end{lem}
\begin{proof}
On se ram{\`e}ne facilement aux consid{\'e}rations similaires sur les espaces de Sobolev
usuels \( L^{i,p} \) de fonctions ayant \( i \) d{\'e}riv{\'e}es dans \( L^{p} \).
En effet, rappelons que \( \mathcal{H}^{i}\subset L^{i/2,2} \) et, classiquement,
\( L^{i/2,2}\subset L^{p} \) pour \( 1/p=1/2-i/2n>0 \). Par cons{\'e}quent, prenons
deux fonctions \( f\in \mathcal{H}^{k} \) et \( g\in \mathcal{H}^{i} \), avec
\( i\leq k \), et montrons que le produit \( fg \) a \( i \) d{\'e}riv{\'e}es horizontales
dans \( L^{2} \) : nous devons donc montrer que \( \nabla ^{j}f\nabla ^{i-j}g\in L^{2} \)
; si \( k-j>n \) ou si \( j>n \), alors l'un des deux facteurs est dans \( C^{0} \)
et le produit est clairement dans \( L^{2} \) ; dans les cas restants, c'est-{\`a}-dire
si \( k-n\leq j\leq n \), alors \( \nabla ^{j}f\in L^{p} \) avec \( 1/p=1/2-(k-j)/2n \)
(si \( k-j\neq n) \), et \( \nabla ^{i-j}g\in L^{q} \) avec \( 1/q=1/2-j/2n \)
(si \( j\neq n \)), donc \( \nabla ^{j}f\nabla ^{i-j}g\in L^{r} \) avec 
\[
\frac{1}{r}=\frac{1}{p}+\frac{1}{q}=1-\frac{k}{2n}<\frac{1}{2};\]
 dans le cas o{\`u} \( j=k-n \) ou \( n \), par exemple si \( j=k-n \), on prend
\( p<\infty  \) suffisamment grand et on arrive quand m{\^e}me {\`a} \( r>2 \).
\end{proof}

\subsection{Espace d'applications de \protect\( M\protect \) dans \protect\( X\protect \)}

Remarquons que si \( \varphi :M\rightarrow M \) est de contact, alors \( \psi =\pi \circ \varphi :M\rightarrow X \)
v{\'e}rifie la condition 
\begin{equation}
\label{con-phi*omega}
\psi ^{*}\omega \in \mathcal{I}(\eta ),
\end{equation}
 o{\`u} \( \mathcal{I}(\eta ) \) est l'id{\'e}al diff{\'e}rentiel engendr{\'e} par \( \eta  \).
Nous pouvons aussi r{\'e}crire la condition d'une mani{\`e}re diff{\'e}rente : notons \( \Lambda ^{2}_{0}H \)
l'espace des 2-vecteurs primitifs (dans le noyau de \( d\eta  \)), et 
\[
\Omega ^{2}_{0}H:=\Lambda ^{2}H^{*}/\mathbf{R}d\eta \]
 l'espace des 2-formes horizontales primitives (on peut identifier \( \Omega ^{2}_{0}H \)
au noyau dans \( \Omega ^{2}H \) de l'op{\'e}rateur de contraction par \( d\eta  \)),
alors (\ref{con-phi*omega}) devient {\'e}quivalente {\`a} 
\[
\psi ^{*}\omega =0\in \Omega ^{2}_{0}H.\]
Cela nous am{\`e}ne {\`a} {\'e}tudier les applications \( M\rightarrow X \) satisfaisant
cette condition.

L'espace des applications \( \psi :M\rightarrow X \) ayant \( k \) d{\'e}riv{\'e}es
horizontales dans \( L^{2} \) est une vari{\'e}t{\'e} hilbertienne, \( \mathcal{H}^{k}(M,X) \),
param{\'e}tr{\'e}e dans un voisinage de \( \pi  \) par un champ de vecteurs horizontal
\( \xi  \) sur \( M \), via l'application 
\[
\xi \longrightarrow (x\rightarrow \exp _{\pi (x)}\pi _{*}\xi ),\]
 o{\`u} l'exponentielle sur \( X \) est prise par rapport {\`a} une m{\'e}trique riemannienne
fix{\'e}e. 

On notera \( \mathcal{P}^{k}(M,X)\subset \mathcal{H}^{k}(M,X) \) le sous-espace
des \( \psi  \) satisfaisant (\ref{con-phi*omega}) ; fixons \( \psi _{0}\in \mathcal{P}^{k}(M,X) \),
et remarquons que 
\[
\psi _{0}^{*}\omega =fd\eta +\alpha \wedge \eta ;\]
 comme \( \omega  \) est ferm{\'e}e, il faut que 
\[
(df-\alpha )\wedge d\eta +d\alpha \wedge \eta =0,\]
 ce qui implique \( \alpha =df \) et finalement 
\[
\psi _{0}^{*}\omega =d(f\eta ).\]

Examinons, au voisinage de \( \psi _{0} \), l'application 
\[
P:\mathcal{H}^{k}(M,X)\rightarrow \mathcal{H}^{k-1}(M,\Omega ^{2}_{0}H)\]
 d{\'e}finie par 
\[
P(\psi )=\psi ^{*}\omega .\]
 Puisque \( d\omega =0 \) et \( \psi _{0}^{*}\omega =d(f\eta ), \) la forme
\( \psi ^{*}\omega  \) reste exacte, donc sa projection \( P\psi  \) sur \( \Omega ^{2}_{0}H \)
est {\`a} valeurs dans \( \Im d_{H} \). Or l'op{\'e}rateur 
\[
d_{H}:\Gamma (M,H^{*})\rightarrow \Gamma (M,\Omega _{0}^{2}H)\]
 s'ins{\`e}re dans un complexe hypoelliptique, le complexe de Rumin \cite{Rum:90,Rum:94}
; il en r{\'e}sulte en particulier que \( \Im d_{H} \) est ferm{\'e} dans \( \mathcal{H}^{k-1}(M,\Omega _{0}^{2}H) \).
De plus, la diff{\'e}rentielle de \( P \) en \( \psi _{0} \) n'est autre que 
\[
\xi \longrightarrow \mathcal{L}_{\xi }d(f\eta )=d\iota _{\xi }d(f\eta );\]
en particulier, la diff{\'e}rentielle en \( \pi  \) est (en notant \( \flat =\#^{-1} \))
\[
\xi \longrightarrow d(\flat \xi )\]
qui est une application submersive 
\[
\mathcal{H}^{k}(M,H)\longrightarrow \Im d_{H}\cap \mathcal{H}^{k-1}(M,\Omega _{0}^{2}H).\]
 Plus g{\'e}n{\'e}ralement, la diff{\'e}rentielle reste submersive en \( \psi _{0} \) tant
que la fonction \( f \) plus haut ne s'annule pas ; on en d{\'e}duit finalement
que \( P^{-1}(0) \) est une sous-vari{\'e}t{\'e} hilbertienne de \( H^{k}(X,Y) \)
: 

\begin{lem}
\label{lem-Pssvar}L'espace \( \mathcal{P}_{0}^{k} \) des applications \( \psi :M\rightarrow X \)
satisfaisant (\ref{con-phi*omega}), telles que \( \psi ^{*}\omega  \) demeure
non d{\'e}g{\'e}n{\'e}r{\'e}e sur \( H \), est une sous-vari{\'e}t{\'e} hilbertienne de \( \mathcal{H}^{k}(M,X) \),
avec espace tangent en \( \pi  \) repr{\'e}sent{\'e} par les champs de vecteurs horizontaux
\( \xi  \) de r{\'e}gularit{\'e} \( \mathcal{H}^{k} \) tels que \( d_{H}(\flat \xi )=0 \).
\end{lem}

\subsection{Rel{\`e}vement}

On a vu que si \( \varphi :M\rightarrow M \) est de contact, alors \( \psi =\pi \circ \varphi  \)
satisfait la condition (\ref{con-phi*omega}). Le lemme suivant montre que l'on
peut presque r{\'e}cup{\'e}rer le contactomorphisme \( \varphi  \) {\`a} partir de \( \psi  \).

\begin{lem}
\label{lem-rel}Supposons \( M \) simplement connexe. Si \( \psi :M\rightarrow X \)
satisfait (\ref{con-phi*omega}), alors on peut relever \( \psi  \) en un contactomorphisme
\( \varphi :M\rightarrow M \) ; ce rel{\`e}vement est unique, {\`a} la composition
pr{\`e}s par l'action d'un {\'e}l{\'e}ment de \( S^{1} \).
\end{lem}
\begin{proof}
On fixe un point de base \( m\in M \), et on choisit un rel{\`e}vement \( \mu \in M \)
du point \( \psi (m) \) ; {\'e}tant donn{\'e} un chemin horizontal \( c(t) \) dans
\( M \) allant de \( m \) {\`a} \( p\in M \), il y a un unique rel{\`e}vement horizontal
\( \mu (t) \) {\`a} \( M \) du chemin \( \psi (c(t)) \), tel que \( \mu (0)=\mu  \).
\begin{claim*}
Le rel{\`e}vement \( \mu (1) \) de \( \psi (p) \) ne d{\'e}pend pas du chemin choisi.
\end{claim*}
Supposons cette affirmation acquise : alors on a d{\'e}fini un rel{\`e}vement \( \varphi :M\rightarrow M \)
de \( \psi  \) qui par construction envoie les directions horizontales sur
les directions horizontales. Le seul choix est celui du relev{\'e} de \( \psi (m) \),
et une modification de ce relev{\'e} entra{\^\i}ne simplement la composition de \( \varphi  \)
avec un {\'e}l{\'e}ment de \( S^{1} \).

Passons maintenant {\`a} la d{\'e}monstration de l'affirmation ci-dessus. Supposons
que l'on ait un lacet \( c(t) \) dans \( M \), donc \( c(0)=c(1)=m \). Il
s'agit de montrer que pour le rel{\`e}vement \( \mu  \), on a \( \mu (1)=\mu (0)=p \).
Puisque \( M \) est simplement connexe, il existe (voir par exemple \cite{All:98})
un disque legendrien dont le bord est ce lacet, c'est-{\`a}-dire que l'on peut trouver
une homotopie legendrienne \( c(s,t)\in M \), telle que 
\begin{eqnarray*}
c(1,t) & = & c(t),\\
c(0,t) & = & m,\\
c(s,0)=c(s,1) & = & m,\\
c^{*}\eta  & = & 0.
\end{eqnarray*}
 On peut relever \( t\rightarrow \varphi (c(s,t)) \) en un chemin horizontal
\( t\rightarrow \mu (s,t) \), avec \( \mu (s,0)=\mu  \). Maintenant, on a
\begin{equation}
\label{int-mu*deta}
\int _{[0,1]\times [0,1]}\mu ^{*}d\eta =\int _{\partial ([0,1]\times [0,1])}\mu ^{*}\eta =\int _{[0,1]}\mu (\cdot ,1)^{*}\eta ,
\end{equation}
 le chemin \( \mu (s,1) \) vit dans une fibre \( S^{1} \) donn{\'e}e, donc l'int{\'e}grale
(\ref{int-mu*deta}) s'annule exactement quand \( \mu (0,1)=\mu (1,1) \), c'est-{\`a}-dire
quand \( \mu (1)=p \) ; d'un autre c{\^o}t{\'e}, par hypoth{\`e}se, 
\[
\mu ^{*}d\eta =\mu ^{*}\pi ^{*}\omega =c^{*}\Phi ^{*}\omega =c^{*}(\eta \wedge \alpha +fd\eta )=0,\]
 {\`a} cause de la condition \( c^{*}\eta =0 \), donc l'int{\'e}grale (\ref{int-mu*deta})
est bien nulle.
\end{proof}

\subsection{Param{\'e}trisation des contactomorphismes}

Nous pouvons {\`a} pr{\'e}sent proposer la param{\'e}trisation suivante des contactomorphismes
de \( M \) qui sont proches de l'identit{\'e}. Fixons un point \( m\in M \) et
une section locale \( \sigma :X\rightarrow M \) de \( \pi  \) telle que \( \sigma (\pi (m))=m \).
Pour \( \psi :M\rightarrow X \) satisfaisant la condition (\ref{con-phi*omega}),
notons \( \tilde{\psi } \) le relev{\'e} de \( \psi  \) construit par le lemme
\ref{lem-rel}, tel que \( \tilde{\psi }(m)=\sigma (\psi (m)) \). Nous pouvons
r{\'e}crire le lemme \ref{lem-rel} de la mani{\`e}re suivante :

\begin{lem}
Si \( M \) est simplement connexe, alors tout contactomorphisme \( \varphi  \)
de \( M \), proche de l'identit{\'e}, ayant \( k \) d{\'e}riv{\'e}es horizontales, s'{\'e}crit
de mani{\`e}re unique \( \varphi =u\circ \tilde{\psi } \), avec \( u\in S^{1} \)
et \( \psi \in \mathcal{P}^{k} \).
\end{lem}
Maintenant, par le lemme \ref{lem-Pssvar}, nous avons une bonne param{\'e}trisation
de \( \mathcal{P}^{k} \) : en effet, puisque \( \mathcal{P}^{k} \) est une
vari{\'e}t{\'e} hilbertienne, nous avons une carte locale 
\[
\Psi :T_{\pi }\mathcal{P}^{k}\longrightarrow \mathcal{P}^{k},\]
 d{\'e}finie sur une petite boule de \( T_{\pi }\mathcal{P}^{k}=\{\alpha \in \mathcal{H}^{k}(M,H^{*}),d_{H}\alpha =0\} \)
; or la cohomologie du complexe de Rumin calcule la cohomologie DeRham, donc,
si \( H^{1}(M,\mathbf{R})=0 \), on peut param{\'e}trer les \( \alpha  \) par 
\[
\alpha =d_{H}f,\quad f\in \mathcal{H}^{k+1},\quad \int _{M}f=0.\]
 Finalement, {\'e}tant donn{\'e}e \( f\in \mathcal{H}^{k+1}, \) suffisamment petite,
on peut poser 
\begin{equation}
\label{for-carte}
\begin{array}{rl}
u & =\exp \int _{M}f,\\
\psi  & =\Psi (d_{H}f),
\end{array}\qquad \varphi =u\circ \tilde{\psi },
\end{equation}
 et on a d{\'e}montr{\'e} le r{\'e}sultat suivant :

\begin{thm}
\label{th-par-contact}Si \( M \) est simplement connexe, l'application (\ref{for-carte})
constitue une param{\'e}trisation locale des contactomorphismes de \( M \), proches
de l'identit{\'e}, et ayant \( k \) d{\'e}riv{\'e}es horizontales, par des fonctions de
r{\'e}gularit{\'e} \( \mathcal{H}^{k+1} \). \qed 
\end{thm}
Cette param{\'e}trisation nous fournit sur les contactomorphismes proches de l'identit{\'e}
la structure de vari{\'e}t{\'e} hilbertienne que nous utiliserons dans la suite. On
notera \( \mathcal{G}^{k} \) le groupe des contactomorphismes de \( M \) avec
\( k \) d{\'e}riv{\'e}es horizontales dans \( L^{2} \).

\section{\label{sec6}Jauge de Coulomb}

Revenons {\`a} pr{\'e}sent {\`a} l'espace des twisteurs \( \mathcal{T} \), qui n'est autre
que le \( S^{1} \)-fibr{\'e} \( \mathcal{O}(-1,1) \) au-dessus de \( P^{1}\times P^{1} \).
{\'E}tant donn{\'e}e une perturbation \( \phi \in \Omega ^{0,1}\otimes T^{1,0} \) de
la structure CR, on cherche {\`a} quelle condition existe une action de \( S^{1} \)
du type produit par le corollaire \ref{cor-ex-disques}. Tout d'abord on a le
r{\'e}sultat suivant.

\begin{lem}
Une action de \( S^{1} \) par contactomorphismes sur \( \mathcal{T} \), proche
de l'action donn{\'e}e, lui est conjugu{\'e}e par un contactomorphisme.
\end{lem}
\begin{proof}
L'action de \( S^{1} \) nous fournit une nouvelle forme de contact, \( \eta ' \),
dont la diff{\'e}rentielle provient de la base \( P^{1}\times P^{1} \) : 
\[
d\eta '=\pi ^{*}\omega ',\]
 o{\`u} \( \omega ' \) est une 2-forme ferm{\'e}e, repr{\'e}sentant ({\`a} une constante pr{\`e}s)
le \( c_{1} \) du fibr{\'e}, et donc dans la m{\^e}me classe de cohomologie que \( \omega  \)
; si \( \omega ' \) est proche de \( \omega  \), alors, par le lemme de Moser,
\( \omega ' \) s'en d{\'e}duit par un diff{\'e}omorphisme, qui se rel{\`e}ve aux \( S^{1} \)-fibr{\'e}s
en un diff{\'e}omorphisme envoyant \( \eta  \) sur \( \eta ' \).
\end{proof}
Ce lemme signifie qu'au lieu de chercher une action de \( S^{1} \) pour laquelle
les coefficients de Fourier de \( \phi  \) soient positifs ou nuls, il suffit
de chercher un contactomorphisme \( \varphi  \) de \( \mathcal{T} \) tel que
\( \varphi ^{*}\phi  \) soit {\`a} coefficients de Fourier positifs ou nuls par
rapport {\`a} l'action standard de \( S^{1} \). 

L'action des contactomorphismes sur les structures CR s'apparente {\`a} l'action
d'un groupe de jauge. L'id{\'e}e est alors, {\`a} l'instar de Bland \cite{Bla:94} sur
\( S^{3} \), de trouver une jauge de Coulomb permettant de fixer le contactomorphisme
; cette jauge sera choisie de sorte de tuer autant de coefficients n{\'e}gatifs
que possible dans \( \phi  \).

\begin{lem}
\label{lem-act-lisse}Le groupe \( \mathcal{G}^{k+1} \) agit de mani{\`e}re lisse
sur \( \mathcal{H}^{k}(\Omega ^{0,1}\otimes T^{1,0}) \) par \( \phi \rightarrow \varphi ^{*}\phi  \)
; l'application tangente {\`a} l'action en l'origine est 
\[
f\in \mathcal{H}^{k+2}\longrightarrow \overline{\partial }_{H}\#\overline{\partial }f,\]
 et son image est ferm{\'e}e. 
\end{lem}
\begin{proof}
Notons \( \psi =\pi \circ \varphi  \), o{\`u} \( \pi  \) est la projection \( \mathcal{T}\rightarrow P^{1}\times P^{1} \)
; on sait que \( \varphi \in \mathcal{G}^{\mathcal{k}+1} \) signifie essentiellement
que \( \psi \in \mathcal{H}^{\mathcal{k}+1}(\mathcal{T},P^{1}\times P^{1}) \)
; l'action de \( \varphi  \) sur \( \phi  \) consiste {\`a} d{\'e}placer le graphe
de \( 1+\phi  \), ce qui se traduit par 
\begin{eqnarray*}
(\varphi ^{*}\phi )_{x} & = & \pi _{1,0}\circ (T_{x}\varphi )^{-1}\circ (1+\phi _{\varphi (x)})\circ \big (\pi _{0,1}\circ (T_{x}\varphi )^{-1}\circ (1+\phi _{\varphi (x)})\big )^{-1}\\
 & = & \pi _{1,0}\circ (T_{x}\psi )^{-1}\circ (1+\phi _{\varphi (x)})\circ \big (\pi _{0,1}\circ (T_{x}\psi )^{-1}\circ (1+\phi _{\varphi (x)})\big )^{-1},
\end{eqnarray*}
 o{\`u} \( \phi  \) est interpr{\'e}t{\'e} comme un morphisme \( \pi ^{*}T^{0,1}(P^{1}\times P^{1})\rightarrow \pi ^{*}T^{1,0}(P^{1}\times P^{1}) \)
; on d{\'e}duit alors facilement la r{\'e}gularit{\'e} annonc{\'e}e, {\`a} partir du lemme \ref{lem-FS-alg}.

La diff{\'e}rentielle en \( f \) de l'action de \( \mathcal{G}^{k+1} \) modifie
un vecteur de type (0,1), \( X, \) par 
\[
[X_{f},X]=[fR-\#df,X];\]
 la structure complexe est modifi{\'e}e par 
\[
\phi _{X}=[X_{f},X]^{1,0}=[-\#\overline{\partial }f,X]^{1,0}=(\overline{\partial }_{H}\#\overline{\partial }f)_{X}.\]

L'image de \( \overline{\partial }_{H}\#\overline{\partial } \) est ferm{\'e}e
car l'op{\'e}rateur \( \overline{\partial } \) est hypoelliptique sur les 0-formes
en signature (1,1), voir la d{\'e}monstration du corollaire \ref{cor-31}.
\end{proof}
\begin{cor}
\label{cor-for-W}Supposons donn{\'e} dans \( \mathcal{H}^{k}(\Omega ^{0,1}\otimes T^{1,0}) \)
un suppl{\'e}mentaire \( W \) de \( \Im \overline{\partial }_{H}\#\overline{\partial } \),
invariant sous l'action des CR-automorphismes de \( \mathcal{T} \). Alors pour
tout \( \phi  \) proche de \( 0 \), il existe un contactomorphisme \( \varphi  \),
proche de l'identit{\'e}, tel que \( \varphi ^{*}\phi \in W \) ; de plus, \( \varphi  \)
est unique modulo les CR-automorphismes de \( \mathcal{T} \).
\end{cor}
\begin{proof}
L'existence d'un suppl{\'e}mentaire \( W \) vient de ce que l'image de \( \overline{\partial }_{H}\#\overline{\partial } \)
est ferm{\'e}e par le lemme \ref{lem-act-lisse}. Le noyau de \( \overline{\partial }_{H}\#\overline{\partial } \)
est form{\'e} des CR-automorphismes infinit{\'e}simaux ; dans le cas de \( \mathcal{T}=\partial \mathbf{H}H^{1}=\partial \mathbf{R}H^{4} \),
le groupe des automorphismes CR est le groupe \( Sp_{1,1}=SO_{1,4} \), c'est-{\`a}-dire
le groupe des transformations conformes de \( S^{3} \) (dont \( \mathcal{T} \)
est l'espace des twisteurs), de dimension 10 ; consid{\'e}rons le sous-groupe 
\[
\mathcal{G}_{0}\subset \mathcal{G}\]
 des contactomorphismes de \( \mathcal{T} \) fixant deux points, de sorte que
\( \mathcal{G}/\mathcal{G}_{0}=SO_{1,4} \) ; consid{\'e}rons l'action de \( \mathcal{G}_{0} \)
suivie de la projection sur \( \Im \overline{\partial }_{H}\#\overline{\partial } \)
parall{\`e}lement {\`a} \( W \) : 
\[
F:\mathcal{G}_{0}^{k+1}\times \mathcal{H}^{k}(\Omega ^{0,1}\otimes T^{1,0})\rightarrow \mathcal{H}^{k}(\Omega ^{0,1}\otimes T^{1,0})\rightarrow \Im \overline{\partial }_{H}\#\overline{\partial },\]
 alors la diff{\'e}rentielle de \( F \) par rapport {\`a} la premi{\`e}re variable est
un isomorphisme, donc, par le th{\'e}or{\`e}me des fonctions implicites, l'{\'e}quation
\( F(\varphi ,\phi )=0, \) que l'on peut r{\'e}crire 
\[
\varphi ^{*}\phi \in W,\]
 a une unique solution proche de l'identit{\'e}, \( \varphi \in \mathcal{G}_{0}^{k+1} \),
si \( \phi  \) est assez petit.
\end{proof}
\begin{lem}
\label{lem-64}Comme suppl{\'e}mentaire \( W \) de \( \Im \overline{\partial }_{H}\#\overline{\partial } \),
on peut choisir 
\[
W=(\Im \overline{\partial }_{H}\#\overline{\partial }\pi _{\leq 0})^{\perp },\]
 o{\`u} \( \pi _{\leq 0} \) est la projection sur les coefficients de Fourier n{\'e}gatifs
ou nuls.

Une structure CR int{\'e}grable \( \phi  \), proche de la structure standard, admet
un remplissage holomorphe si et seulement si, dans la jauge de Coulomb \( \varphi ^{*}\phi \in W \)
construite par le corollaire \ref{cor-for-W}, on a \( \varphi ^{*}\phi  \)
{\`a} coefficients de Fourier positifs ou nuls.
\end{lem}
\begin{proof}
Rappelons que nous regardons l'image des fonctions \emph{r{\'e}elles} ; si \( f \)
est r{\'e}elle, alors les coefficients de Fourier positifs de \( f \) sont d{\'e}termin{\'e}s
par ses coefficients de Fourier n{\'e}gatifs, et il en est de m{\^e}me pour l'image
\( \overline{\partial }_{H}\#\overline{\partial }f \) , puisque \( \overline{\partial }_{H}\#\overline{\partial } \)
commute {\`a} l'action de \( S^{1} \) ; on en d{\'e}duit que l'espace \( W \) propos{\'e}
est bien un suppl{\'e}mentaire de \( \Im \overline{\partial }_{H}\#\overline{\partial } \). 

Pour montrer la seconde assertion, notons tout d'abord que si en jauge de Coulomb
on a des coefficients de Fourier positifs ou nuls, alors le remplissage holomorphe
existe par le lemme \ref{lem-fill-J} ; r{\'e}ciproquement, par le corollaire \ref{cor-ex-disques},
si un remplissage existe, alors on a un contactomorphisme \( \varphi  \) tel
que \( \varphi ^{*}\psi  \) soit {\`a} coefficients de Fourier positifs ou nuls
; il est clair que \( \varphi ^{*}\psi  \) est dans l'orthogonal de \( \Im \overline{\partial }_{H}\#\overline{\partial }\pi _{<0} \),
mais a priori il n'y a pas de raison pour que \( \varphi ^{*}\psi  \) soit
orthogonal aux composantes \( S^{1} \)-invariantes de \( \Im \overline{\partial }_{H}\#\overline{\partial }\pi  \)
; cependant, on peut obtenir aussi cette condition, comme dans le corollaire
\ref{cor-for-W}, en composant \( \varphi  \) par un contactomorphisme qui
sera cette fois \( S^{1} \)-invariant (et donc pr{\'e}servera l'annulation des
coefficients strictement n{\'e}gatifs) ; par cons{\'e}quent, apr{\`e}s cette modification,
\( \varphi ^{*}\phi \in W \), ce qui signifie, par unicit{\'e} de la jauge de Coulomb
(aux CR-automorphismes pr{\`e}s), que \( \varphi  \) est le contactomorphisme mettant
\( \phi  \) en jauge de Coulomb ; en particulier, dans cette jauge, seuls subsistent
les coefficients de Fourier positifs ou nuls.
\end{proof}

\section{\label{sec7}L'image du hessien complexe}

Pour utiliser plus concr{\`e}tement le lemme \ref{lem-64}, il faut d{\'e}terminer l'image
du hessien complexe \( \overline{\partial }_{H}\#\overline{\partial } \). Pour
cela, nous utilisons le fait important que l'espace des twisteurs \( \mathcal{T} \)
est un espace homog{\`e}ne sous le groupe \( Sp_{1}Sp_{1} \).

\subsection{L'espace des twisteurs comme espace homog{\`e}ne}

Rappelons que par (\ref{dec-T}), l'espace des twisteurs \( \mathcal{T} \)
est d{\'e}crit comme une hypersurface r{\'e}elle de \( P^{3} \) ; l'action de \( Sp_{1}Sp_{1} \)
sur \( [z_{1}:z_{2}:z_{3}:z_{4}] \) est simplement donn{\'e}e par l'action sur
chaque couple de coordonn{\'e}es \( (z_{1},z_{2}) \) et \( (z_{3},z_{4}) \) ;
ainsi les actions de \( Sp_{1}Sp_{1} \) sur \( \mathcal{T} \) et sur \( P^{1}\times P^{1} \)
sont-elles compatibles avec la projection \( \mathcal{T}\rightarrow P^{1}\times P^{1} \).
En outre, l'action de \( Sp_{1}Sp_{1} \) pr{\'e}serve toutes les structures d'espace
des twisteurs de \( \mathcal{T} \), c'est-{\`a}-dire la structure CR, la structure
r{\'e}elle, la structure de contact holomorphe, et les \( P^{1} \) r{\'e}els ; on remarquera
que, dans ce formalisme, la base de la projection twistorielle s'{\'e}crit comme
le quotient \( S^{3}=Sp_{1}Sp_{1}/SO_{3} \).

D{\'e}composons l'alg{\`e}bre de Lie de \( Sp_{1} \) comme 
\[
\mathfrak {sp}_{1}=\mathfrak {u}_{1}\oplus \mathbf{R}^{2},\quad \mathfrak {u}_{1}=\left( \begin{array}{cc}
* & 0\\
0 & *
\end{array}\right) ,\quad \mathbf{R}^{2}=\left( \begin{array}{cc}
0 & *\\
* & 0
\end{array}\right) .\]
 Le stabilisateur du point base, 
\[
*=[1:0:1:0],\]
 est un \( U_{1} \) plong{\'e} diagonalement dans \( Sp_{1}Sp_{1} \), de sorte
que \( \mathcal{T}=Sp_{1}Sp_{1}/U_{1} \) ; on a donc la d{\'e}composition 
\begin{eqnarray*}
\mathfrak {sp}_{1}\oplus \mathfrak {sp}_{1} & = & \mathfrak {u}_{1}\oplus \mathfrak {u}_{1}\oplus \mathbf{R}^{2}\oplus \mathbf{R}^{2}\\
 & = & \mathfrak {u}^{+}_{1}\oplus \mathfrak {m},\\
\mathfrak {m} & = & \mathfrak {u}^{-}_{1}\oplus \mathbf{R}^{2}\oplus \mathbf{R}^{2},
\end{eqnarray*}
o{\`u} \( \mathfrak {u}^{+}_{1}=\{(x,x)\in \mathfrak {u}_{1}\oplus \mathfrak {u}_{1}\} \)
est le stabilisateur de \( * \) et \( \mathfrak {m} \) son orthogonal ; le
\( \mathfrak {u}^{-}_{1}=\{(x,-x)\in \mathfrak {u}_{1}\oplus \mathfrak {u}_{1}\} \)
repr{\'e}sente la direction du champ de Reeb, et \( \mathbf{R}^{2}\oplus \mathbf{R}^{2} \)
les directions de contact, munies de la structure CR naturelle.

Cette description fournit une identification naturelle du fibr{\'e} tangent de \( \mathcal{T} \)
et de la distribution de contact \( H \) comme fibr{\'e}s homog{\`e}nes : 
\[
T\mathcal{T}=Sp_{1}Sp_{1}\times _{U_{1}}\mathfrak {m},\quad H=Sp_{1}Sp_{1}\times _{U_{1}}(\mathbf{R}^{2}\oplus \mathbf{R}^{2}).\]
 La modification \( \phi  \) de la structure CR s'{\'e}crit comme une section ({\'e}quivariante)
de \( \mathfrak {m}^{0,1}\otimes \mathfrak {m}_{1,0} \), soit, en choisissant
une base \( (e_{1},e_{2}) \) de \( \mathfrak {m}_{1,0} \),
\[
\phi =\phi _{\overline{\imath }}^{j}e^{\overline{\imath }}\otimes e_{j}.\]

La courbure du \( S^{1} \)-fibr{\'e} s'{\'e}crit \( F=-(e^{1}\wedge e^{\overline{1}}-e^{2}\wedge e^{\overline{2}}) \),
de sorte que la condition \( \phi \lrcorner F=0 \) de (\ref{int-bord}) devient
\begin{equation}
\label{Rphi}
\phi _{\overline{2}}^{1}+\phi ^{2}_{\overline{1}}=0.
\end{equation}

Si \( \phi  \) provient de la construction twistorielle, telle que d{\'e}crite
dans la section \ref{sec11}, alors \( \phi  \) doit pr{\'e}server les \( P^{1} \)
r{\'e}els ; or le \( T^{1,0} \) des \( P^{1} \) r{\'e}els est engendr{\'e} par \( e_{1}+e_{2}, \)
ce qui conduit {\`a} la condition 
\begin{equation}
\label{phiP1}
\phi _{\overline{1}}^{1}+\phi _{\overline{2}}^{1}-\phi _{\overline{1}}^{2}-\phi _{\overline{2}}^{2}=0.
\end{equation}

Enfin, l'op{\'e}rateur \( \#:\mathfrak {m}^{0,1}\rightarrow \mathfrak {m}_{1,0} \)
est donn{\'e} par 
\begin{equation}
\label{diese}
\#e^{\overline{1}}=-ie_{1},\quad \#e^{\overline{2}}=ie_{2}.
\end{equation}

\subsection{D{\'e}composition harmonique}

L'id{\'e}e est bien s{\^u}r d'utiliser la d{\'e}composition harmonique sous les repr{\'e}sentations
irr{\'e}ductibles de \( Sp_{1}Sp_{1} \) pour comprendre l'image de l'op{\'e}rateur
homog{\`e}ne \( \overline{\partial }_{H}\#\overline{\partial } \). Notons \( \sigma ^{k}=Sym^{k}\mathbf{C}^{2} \)
les repr{\'e}sentations irr{\'e}ductibles de \( Sp_{1} \) ; les repr{\'e}sentations irr{\'e}ductibles
\( \rho  \) de \( Sp_{1}\times Sp_{1} \) sont simplement les produits des
repr{\'e}sentations des deux facteurs, 
\[
V_{\rho }=\sigma ^{K}\otimes \sigma ^{L},\]
 elles descendent {\`a} \( Sp_{1}Sp_{1}=Sp_{1}\times Sp_{1}/\mathbf{Z}_{2} \) si
\( K+L \) est pair.

Si l'on dispose d'un fibr{\'e} homog{\`e}ne \( E=Sp_{1}Sp_{1}\times _{\rho _{0}}W \),
o{\`u} \( \rho _{0} \) est une repr{\'e}sentation de \( U_{1} \) dans \( W \), alors
l'espace de ses sections \( L^{2} \) se d{\'e}compose comme une somme directe hilbertienne,
\[
L^{2}(E)=\sum _{\rho }V_{\rho }\otimes \hom _{U_{1}}(V_{\rho },W),\]
 o{\`u} un {\'e}l{\'e}ment \( v\otimes w\in V_{\rho }\otimes \hom _{U_{1}}(V_{\rho },W) \)
correspond {\`a} la section de \( E \) sur \( \mathcal{T} \) donn{\'e}e par l'application
{\'e}quivariante \( s:Sp_{1}Sp_{1}\rightarrow W \), d{\'e}finie par 
\[
s(g)=\left\langle w,\rho (g^{-1})v\right\rangle .\]

Regardons comment d{\'e}crire l'action de la structure r{\'e}elle : rappelons que, par
(\ref{def-tau-N}), la structure r{\'e}elle s'{\'e}crit \( \tau ([z^{1}:z^{2}:z^{3}:z^{4}])=[-\overline{z^{2}}:\overline{z^{1}}:-\overline{z^{4}}:\overline{z^{3}}] \)
; en particulier, pour le point base \( * \), on a 
\[
\tau ([1:0:1:0])=[0:1:0:1]=g_{0}([1:0:1:0]),\]
avec 
\[
g_{0}=\left( \begin{array}{cc}
0 & i\\
i & 0
\end{array}\right) \in Sp_{1}^{diag}\subset Sp_{1}Sp_{1};\]
 on notera que \( g_{0}^{2}=-1 \) dans \( Sp_{1} \), donc \( g_{0}^{2}=1 \)
dans \( Sp_{1}Sp_{1} \). Puisque la structure r{\'e}elle commute {\`a} l'action de
\( Sp_{1}Sp_{1} \), on d{\'e}duit que, pour une fonction \( f(g)=\left\langle w,\rho (g^{-1})v\right\rangle  \),
on a 
\begin{eqnarray}
\tau ^{*}f(g*) & = & f(\tau g*)\nonumber \\
 & = & f(g\tau *)\nonumber \\
 & = & f(gg_{0}*)\nonumber \\
 & = & \left\langle w\circ \rho (g_{0}^{-1}),\rho (g^{-1})v\right\rangle .\label{tau*f} 
\end{eqnarray}
 De mani{\`e}re similaire, on obtient, si \( \phi (g*)=g_{*}\left\langle w,\rho (g^{-1})v\right\rangle  \),
\begin{eqnarray}
\tau ^{*}\overline{\phi }(g*) & = & \tau _{*}\overline{\phi }(\tau g*)\nonumber \\
 & = & \tau _{*}(gg_{0})_{*}\left\langle \overline{w}\circ \rho (g_{0})^{-1},\rho (g^{-1})v\right\rangle \nonumber \\
 & = & g_{*}(\tau g_{0})_{*}\left\langle \overline{w}\circ \rho (g_{0})^{-1},\rho (g^{-1})v\right\rangle ;\label{cal-tau-phi} 
\end{eqnarray}
 en prenant des coordonn{\'e}es complexes \( (z^{2},z^{4}) \) sur \( \mathbf{R}^{2}\oplus \mathbf{R}^{2} \),
comme dans le calcul menant {\`a} la formule (\ref{def-etac}), on voit que 
\[
(\tau g_{0})_{*}(z^{2},z^{4})=-(\overline{z^{2}},\overline{z^{4}});\]
 en particulier, si \( w=w_{\overline{\imath }}^{j}e^{\overline{\imath }}\otimes e_{j} \),
alors 
\[
\tau 'w:=(\tau g_{0})_{*}\overline{w}=\overline{w_{\overline{\imath }}^{j}}e^{\overline{\imath }}\otimes e_{j}\]
 est anti-\( \mathbf{C} \)-lin{\'e}aire et anti-\( U_{1} \)-invariant (car \( g_{0} \)
anticommute {\`a} l'action de \( U_{1} \)) ; comme nous regardons \( Sp_{1}Sp_{1} \)
plut{\^o}t que le produit \( Sp_{1}\times Sp_{1} \), la repr{\'e}sentation \( V_{\rho }=\sigma ^{K}\otimes \sigma ^{L} \)
(avec \( K+L \) pair) admet toujours une structure r{\'e}elle \( \tau _{\rho } \),
ce qui nous permet d'{\'e}crire finalement 
\begin{equation}
\label{tau*phi}
\tau ^{*}\overline{\phi }(g*)=g_{*}\left\langle \tau 'w\circ \rho (g_{0}^{-1})\circ \tau _{\rho },\rho (g^{-1})\tau _{\rho }v\right\rangle ;
\end{equation}
 dans cette formule, on notera que \( \tau 'w\circ \rho (g_{0}^{-1})\circ \tau _{\rho } \)
est bien maintenant un {\'e}l{\'e}ment de \( \hom _{U_{1}}(V_{\rho },\mathfrak {m}^{0,1}\otimes \mathfrak {m}_{1,0}) \),
donc on a bien repr{\'e}sent{\'e} \( \tau ^{*}\overline{\phi } \) dans la d{\'e}composition
harmonique. 

De la m{\^e}me mani{\`e}re, si on a une fonction \( f=\left\langle w,\rho (g^{-1})v\right\rangle  \),
alors on peut repr{\'e}senter sa conjugu{\'e}e comme 
\begin{equation}
\label{f-bar}
\overline{f}=\left\langle \overline{w}\circ \tau _{\rho },\rho (g^{-1})\tau _{\rho }v\right\rangle .
\end{equation}

Passons maintenant {\`a} l'op{\'e}rateur \( \overline{\partial } \) : fixons dans chaque
facteur de \( (\mathfrak {sp}_{1}\oplus \mathfrak {sp}_{1})\otimes \mathbf{C} \)
le \( \mathfrak {sl}_{2} \)-triplet \( (H_{i},X_{i},Y_{i})_{i=1,2} \) 
\[
H_{i}=\left( \begin{array}{cc}
1 & 0\\
0 & -1
\end{array}\right) ,\quad X_{i}=\left( \begin{array}{cc}
0 & 1\\
0 & 0
\end{array}\right) ,\quad Y_{i}=\left( \begin{array}{cc}
0 & 0\\
1 & 0
\end{array}\right) ;\]
 on notera que les \( Y_{i} \) engendrent \( \mathfrak {m}_{0,1} \), \( \mathfrak {m}_{1,0} \)
l'{\'e}tant par les \( X_{i} \) ; la notation \( (X_{i},Y_{i}) \) est donc redondante
avec \( (e_{i},e_{\overline{\imath }}) \), mais sera utilis{\'e}e quand on voudra
souligner l'aspect repr{\'e}sentation de \( \mathfrak {sl}_{2} \). 

De ce point de vue, une repr{\'e}sentation \( V_{\rho }=\sigma ^{K}\otimes \sigma ^{L} \)
se d{\'e}compose en sous-espaces propres de dimension 1 pour les valeurs propres
\( ((k=K,K-2,...,-K),(l=L,L-2,...,-L)) \) du couple \( (H_{1},H_{2}) \) :
\[
V_{\rho }=\sum _{k,l}V_{\rho }(k,l);\]
l'action de \( Y_{1} \) (resp. \( Y_{2} \)) diminue le poids \( k \) (resp.
\( l \)) de 2 : 
\[
\begin{array}{ccc}
V(k-2,l) & \stackrel{{Y_{1}}}{\longleftarrow } & V(k,l)\\
 &  & \downarrow Y_{2}\\
 &  & V(k,l-2)
\end{array}.\]

Nous pouvons maintenant traduire la condition infinit{\'e}simale \( \overline{\partial }_{H}\phi =0 \)
issue de (\ref{int-bord}) en termes de cette d{\'e}composition harmonique : {\'e}crivons
\( \phi =\sum \left\langle w,\rho (g^{-1})v\right\rangle  \), o{\`u} \( v\otimes w\in V_{\rho }\otimes \hom _{U_{1}}(V_{\rho },\mathfrak {m}^{0,1}\otimes \mathfrak {m}_{1,0}) \),
alors, pour \( \dot{g}\in \mathfrak {m}_{0,1} \), on a 
\[
\nabla ^{0,1}\phi (g\dot{g})=\sum \left\langle w\circ \rho _{*}(\dot{g}),\rho (g^{-1})v\right\rangle ,\]
 si bien que la condition \( \overline{\partial }_{H}\phi =0 \) devient 
\begin{equation}
\label{dbar-phi}
w_{\overline{2}}\rho _{*}(Y_{1})-w_{\overline{1}}\rho _{*}(Y_{2})=0.
\end{equation}

Passons au hessien complexe \( \overline{\partial }_{H}\#\overline{\partial } \)
agissant sur les fonctions : pour une fonction \( f=\left\langle w,\rho (g^{-1})v\right\rangle  \),
on obtient 
\[
\overline{\partial }_{H}\#\overline{\partial }f=\left\langle \#(w\circ \rho _{*})\circ \rho _{*},\rho (g^{-1})v\right\rangle ,\]
 o{\`u} \( \#(w\circ \rho _{*})\circ \rho _{*}\in \hom _{U_{1}}(V_{\rho },\mathfrak {m}^{0,1}\otimes \mathfrak {m}_{1,0}) \)
; plus concr{\`e}tement, en {\'e}crivant \( \#(w\circ \rho _{*})\circ \rho _{*}=\phi _{\overline{\imath }}^{j}e^{\overline{\imath }}\otimes e_{j} \),
on obtient 
\begin{equation}
\label{cal-phi}
\begin{array}{rcl}
\phi _{\overline{1}}^{1} & = & -iw\circ (\rho _{*}Y_{1})^{2},\\
\phi _{\overline{1}}^{2} & = & iw\circ \rho _{*}Y_{2}\circ \rho _{*}Y_{1},\\
\phi _{\overline{2}}^{1} & = & -iw\circ \rho _{*}Y_{1}\circ \rho _{*}Y_{2},\\
\phi _{\overline{2}}^{2} & = & iw\circ (\rho _{*}Y_{2})^{2}.
\end{array}
\end{equation}
 Du point de vue des poids, \( \phi  \) d{\'e}pend de \( w \) de la mani{\`e}re suivante
: 

\begin{itemize}
\item \( \phi _{\overline{1}}^{1}(k,l) \) est d{\'e}termin{\'e} par \( w(k-4,l) \) ;
\item \( \phi _{\overline{1}}^{2}(k,l)=-\phi _{\overline{2}}^{1}(k,l) \) est d{\'e}termin{\'e}
par \( w(k-2,l-2) \) ;
\item \( \phi _{\overline{2}}^{2}(k,l) \) est d{\'e}termin{\'e} par \( w(k,l-4) \).
\end{itemize}
Notons que la compatibilit{\'e} {\`a} la structure r{\'e}elle s'{\'e}crit \( \tau ^{*}\overline{\phi }=\phi  \)
et est impliqu{\'e}e, si \( \phi =\overline{\partial }_{H}\#\overline{\partial }f \),
par la condition \( \tau ^{*}\overline{f}=-f \) (soit \( \tau ^{*}f=-f \)
puisque \( f \) est r{\'e}elle).

Enfin, nous avons une action de \( S^{1} \) sur \( \mathcal{T} \), donn{\'e}e
par 
\[
u[z^{1}:z^{2}:z^{3}:z^{4}]=[uz^{1}:uz^{2}:z^{3}:z^{4}];\]
 cette action n'est pas fournie par un {\'e}l{\'e}ment du groupe \( Sp_{1}Sp_{1} \)
(contrairement {\`a} \( Sp_{1}Sp_{1} \), elle ne pr{\'e}serve pas la structure r{\'e}elle),
mais elle commute {\`a} l'action de \( Sp_{1}Sp_{1} \), et va donc pr{\'e}server les
d{\'e}compositions harmoniques ; sur le point base \( * \), on a, pour \( u\in S^{1} \),
\[
u*=[u:0:1:0]\]
 qui est {\'e}gal {\`a} l'action de \( u\times 1\in Sp_{1}Sp_{1} \) sur \( * \) ;
par cons{\'e}quent, si on a une section \( s(g*)=g_{*}\left\langle w,\rho (g^{-1})v\right\rangle  \)
d'un fibr{\'e} tensoriel \( E \) sur \( \mathcal{T} \), alors, comme dans le calcul
de (\ref{cal-tau-phi}), on obtient 
\[
u_{*}s(g*)=g_{*}\left\langle u(u\times 1)^{-1}\circ w\circ \rho (u\times 1),\rho (g^{-1})v\right\rangle ;\]
 notons que 
\begin{equation}
\label{act-U1-E}
u(u\times 1)^{-1}[1:z^{2}:1:z^{4}]=[1:u^{2}z^{2}:1:z^{4}];
\end{equation}
 dans la d{\'e}composition 
\[
\hom _{U_{1}}(V_{\rho },E)=\sum \hom _{U_{1}}(V_{\rho }(k,l),E),\]
 les poids de l'action de \( S^{1} \) sont donc les sommes des \( k \) et
des poids de l'action de \( S^{1} \) sur \( E \) par (\ref{act-U1-E}) ; en
particulier, pour \( E=\mathfrak {m}^{0,1}\otimes \mathfrak {m}_{1,0} \), les
poids de l'action (\ref{act-U1-E}) sont 
\begin{eqnarray}
-4 & \textrm{sur} & e^{\overline{1}}\otimes e_{1},\nonumber \\
-2 & \textrm{sur} & e^{\overline{1}}\otimes e_{2},\, e^{\overline{2}}\otimes e_{1},\label{poids-S1} \\
0 & \textrm{sur} & e^{\overline{2}}\otimes e_{2};\nonumber 
\end{eqnarray}
 on v{\'e}rifie ais{\'e}ment dans (\ref{cal-phi}) que, comme il se doit, l'op{\'e}rateur
\( \overline{\partial }_{H}\#\overline{\partial } \), commutant {\`a} l'action
de \( S^{1} \), pr{\'e}serve ses poids.

\subsection{Calcul de l'image}

Nous sommes maintenant pr{\^e}ts {\`a} calculer l'image du hessien complexe, en tenant
compte de l'action de \( S^{1} \) ; nous fixons une repr{\'e}sentation \( V_{\rho }=V^{K,L} \)
de \( Sp_{1}Sp_{1} \) et voulons comprendre l'image des fonctions r{\'e}elles par
\( \overline{\partial }_{H}\#\overline{\partial } \), en nous int{\'e}ressant plus
particuli{\`e}rement aux composantes {\`a} poids n{\'e}gatif par rapport {\`a} l'action de \( S^{1} \). 

Nous partons d'un 
\[
\phi =\left\langle w,\rho (g^{-1})v\right\rangle ,\quad w\in \hom _{U_{1}}(V_{\rho },\mathfrak {m}^{0,1}\otimes \mathfrak {m}_{1,0});\]
 remarquons que le poids de l'action de \( U_{1} \) sur \( \mathfrak {m}^{0,1}\otimes \mathfrak {m}_{1,0} \)
est {\'e}gal {\`a} \( 4 \), donc les seules composantes de \( w \) non triviales sur
la d{\'e}composition \( V_{\rho }=\sum V(k,l) \) sont {\`a} \( k+l=4 \). De plus les
coefficients \( w_{\overline{\imath }}^{j} \) satisfont les conditions (\ref{Rphi})
et (\ref{dbar-phi}) : la seconde condition signifie que \( w_{\overline{2}}(k,l) \)
est d{\'e}termin{\'e} par \( w_{\overline{1}}(k-2,l+2) \), et la premi{\`e}re que \( w_{\overline{2}}^{1}=-w_{\overline{1}}^{2} \)
; autrement dit, sur les quatre coefficients de \( w(k,l) \), les coefficients
\( w_{\overline{2}}^{1}(k,l) \), \( w_{\overline{2}}^{2}(k,l) \) et \( w_{\overline{1}}^{2}(k,l) \)
sont d{\'e}termin{\'e}s par \( w_{\overline{1}}^{1}(k-2,l+2) \) et \( w_{\overline{1}}^{2}(k-2,l+2) \),
pourvu que \( k-2\geq -K \). 

Nous cherchons une fonction r{\'e}elle 
\[
f=\left\langle w',\rho (g^{-1})v\right\rangle ,\quad w'\in \hom _{U_{1}}(V_{\rho },\mathbf{C}),\]
 telle que \( \overline{\partial }_{H}\#\overline{\partial }f=\phi  \), du
moins sur les composantes {\`a} poids n{\'e}gatifs pour \( S^{1} \) (la condition que
\( f \) soit r{\'e}elle implique que les coefficients sur les poids positifs soient
d{\'e}termin{\'e}s par ceux sur les poids n{\'e}gatifs) ; les seules composantes non nulles
de \( w' \) doivent bien entendu se trouver sur la droite \( k+l=0 \), et
l'op{\'e}rateur est fourni alg{\'e}briquement par les formules (\ref{cal-phi}).

\begin{thm}
\label{th-im-dd}{\'E}tant donn{\'e} \( \phi  \) sur \( \mathcal{T} \), satisfaisant
\( \phi \lrcorner F=0 \) et \( \overline{\partial }_{H}\phi =0 \), il existe
une fonction \( f \) telle que
\begin{enumerate}
\item \( \overline{\partial }_{H}\#\overline{\partial }f_{-}=\phi _{-} \), o{\`u} \( f_{-} \)
et \( \phi _{-} \) repr{\'e}sentent les composantes {\`a} coefficients de Fourier strictement
n{\'e}gatifs par rapport {\`a} l'action de \( U_{1} \),
\item la partie de \( f \) invariante sous \( U_{1} \) est nulle,
\item \( f \) est r{\'e}elle, 
\end{enumerate}
si et seulement si la d{\'e}composition \( \phi =\sum \left\langle w,\rho (g^{-1})v\right\rangle  \)
satisfait les conditions :

\begin{enumerate}
\item si \( V_{\rho }=V^{K,K+2} \), alors \( w_{\overline{1}}^{1}(-K+2,K+2)=0 \)
;
\item si \( V_{\rho }=V^{K,L} \) avec \( K\leq L-4 \), alors 
\[
w_{\overline{1}}^{1}(-K,K+4)=w_{\overline{1}}^{2}(-K,K+4)=w_{\overline{2}}^{1}(-K,K+4)=0.\]

\end{enumerate}
De plus, si \( \phi  \) est compatible {\`a} la structure r{\'e}elle (\( \tau ^{*}\overline{\phi }=\phi  \)),
alors \( f \) l'est aussi (\( \tau ^{*}f=-f \)).

\end{thm}
\begin{proof}
Il y a trois cas suivant les valeurs de \( K \) et \( L \).

\noindent \textbf{Premier cas :} \( K\geq L \). Dans ce cas (voir la figure
\ref{fig-1}), pour \( k<0 \), on peut trouver \( w'(k,-k) \) unique tel que
les relations (\ref{cal-phi}) soient satisfaites pour \( k>0 \). En effet,
le terme \( w(k,l) \) pour le plus petit \( k \) possible est \( w(-L+4,L) \),
et on peut trouver \( w'(-L,L) \), \( w'(-L+2,L-2) \) et \( w'(-L+4,L-4) \)
dont l'image tue exactement \( w(-L+4,L) \) ; le terme suivant est \( w(-L+6,L-2) \)
: dans celui-ci, les coefficients \( w_{\overline{1}}^{1} \), \( w_{\overline{1}}^{2} \)
et \( w_{\overline{2}}^{1} \) sont d{\'e}termin{\'e}s, comme on l'a vu, par \( w(-L+4,L) \)
et en r{\'e}alit{\'e} d{\'e}j{\`a} tu{\'e}s par \( w'(-L+2,L-2) \) et \( w'(-L+4,L-4) \), donc
il ne reste que \( w_{\overline{2}}^{2}(-L+6,L-2) \) qui est tu{\'e} par le choix
judicieux de \( w'(-L+6,L-6) \) ; le processus se poursuit pour \( k<0 \).
Le fait que les relations de compatibilit{\'e} (\ref{dbar-phi}) soient exactement
les obstructions {\`a} r{\'e}soudre le probl{\`e}me de proche en proche est une cons{\'e}quence
du fait que l'image de \( \overline{\partial }_{H}\#\overline{\partial } \)
est dans le noyau de \( \overline{\partial } \) ; la v{\'e}rification directe est
un exercice facile sur les repr{\'e}sentations de \( \mathfrak {sl}_{2} \). 

\begin{figure}[hbt]
{\par\centering \includegraphics{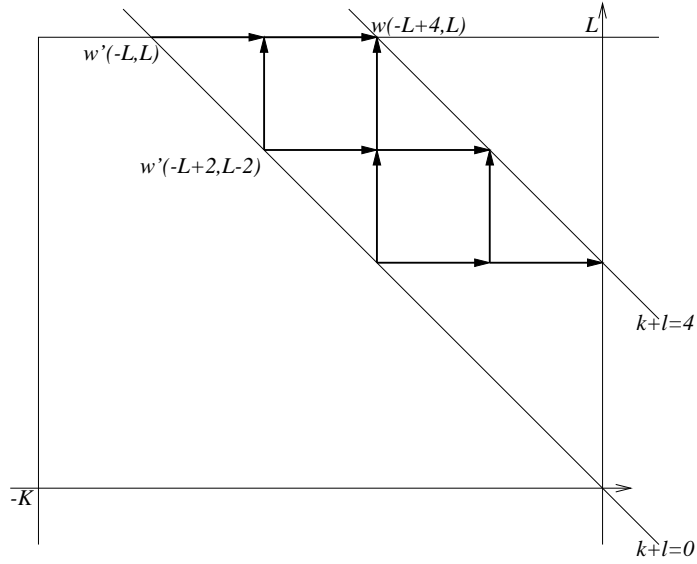} \par}

\caption{\label{fig-1}cas \protect\( K\geq 0\protect \)}
\end{figure}

La condition que \( f \) soit r{\'e}elle s'{\'e}crit, d'apr{\`e}s (\ref{f-bar}), comme
l'invariance par \( w'\otimes v\rightarrow w'\circ \tau _{\rho }\otimes \tau _{\rho }v \)
; comme \( \tau _{\rho } \) envoie \( V(k,l) \) sur \( V(-k,-l) \), on voit
que les \( w'(k,-k) \) pour \( k>0 \) sont d{\'e}termin{\'e}s par les \( w'(k,-k) \)
pour \( k<0 \) ; la seule libert{\'e} restante est donc celle de \( w'(0,0) \),
sur laquelle nous reviendrons plus loin.

\noindent \textbf{Deuxi{\`e}me cas :} \( K=L-2 \). Le processus est exactement
le m{\^e}me (figure \ref{fig-2}), {\`a} la diff{\'e}rence que le terme \( w(k,l) \) avec
\( k \) le plus petit possible est maintenant \( w(-K+2,K+2) \), et on ne
peut pas tuer le terme \( w_{\overline{1}}^{1}(-K+2,K+2) \) par un terme \( w'(-K-2,K+2) \)
qui n'existe pas ; la r{\'e}solution est donc possible si et seulement si \( w_{\overline{1}}^{1}(-K+2,K+2)=0 \).

\begin{figure}[hbt]
{\par\centering \includegraphics{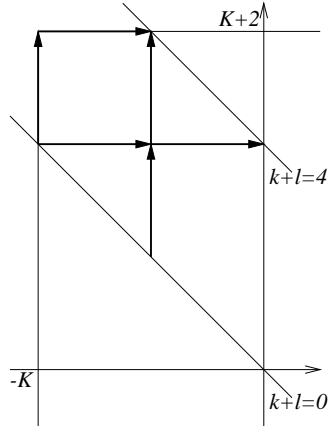} \par}

\caption{\label{fig-2}cas \protect\( K=L-2\protect \)}
\end{figure}

\noindent \textbf{Troisi{\`e}me cas :} \( K\leq L-4 \). Le processus est {\`a} nouveau
identique : pour pouvoir tuer le terme \( w(-K,K+4), \) if faut que 
\[
w_{\overline{1}}^{1}(-K,K+4)=w_{\overline{1}}^{2}(-K,K+4)=w_{\overline{2}}^{1}(-K,K+4)=0;\]
 pour pouvoir tuer le terme \( w(-K+2,K+2), \) il faut que 
\[
w^{1}_{\overline{1}}(-K+2,K+2)=0,\]
 mais cette condition est une cons{\'e}quence de l'annulation de \( w_{\overline{2}}^{1}(-K,K+4) \)
et de la contrainte (\ref{dbar-phi}) ; ensuite on peut r{\'e}soudre de proche en
proche.

Finalement, examinons la compatibilit{\'e} {\`a} la structure r{\'e}elle : rappelons que
par (\ref{tau*phi}), l'op{\'e}ration \( \phi \rightarrow \tau ^{*}\overline{\phi } \)
se traduit par 
\[
w\otimes v\longrightarrow \tau 'w\rho (g^{-1})\tau _{\rho }\otimes \tau _{\rho }v,\]
 et \( f\rightarrow -\tau ^{*}f=-\tau ^{*}\overline{f} \) par 
\[
w'\otimes v\longrightarrow -\overline{w'}\tau _{\rho }\rho (g_{0}^{-1})\otimes \tau _{\rho }v\]
 d'apr{\`e}s (\ref{tau*f}) et (\ref{f-bar}) ; par cons{\'e}quent, il suffit de v{\'e}rifier
que la modification \( w\rightarrow \tau 'w \) impose \( w'\rightarrow -\overline{w'} \)
sur la solution ; comme \( \tau ' \) est simplement la conjugaison sur les
coefficients \( w^{j}_{\overline{\imath }} \), c'est une cons{\'e}quence claire
des formules (\ref{cal-phi}).

\end{proof}

\section{\label{sec8}L'espace des structures CR remplissables}

Notons l'espace des structures CR int{\'e}grables sur \( \mathcal{T} \), de r{\'e}gularit{\'e}
\( \mathcal{H}^{k} \), 
\[
\mathcal{C}^{k}=\{\phi \in \mathcal{H}^{k}(\Omega ^{0,1}\otimes T^{1,0}),\overline{\partial }_{H}\phi +\frac{1}{2}[\phi ,\phi ]=0,\phi \lrcorner F=0\};\]
 bien entendu, pour que \( \phi  \) d{\'e}finisse une structure CR, il faut que
\( T^{0,1}_{\phi }\cap \overline{T^{0,1}_{\phi }}=0 \), mais cette condition
est automatiquement satisfaite pour les petites perturbations auxquelles nous
nous int{\'e}ressons. 

\begin{lem}
L'op{\'e}rateur 
\[
\overline{\partial }_{H}:\{\phi \in \mathcal{H}^{k}(\Omega ^{0,1}\otimes T^{1,0}),\phi \lrcorner F=0\}\longrightarrow \mathcal{H}^{k-1}(\Omega ^{0,2}\otimes T^{1,0})\]
 est surjectif.
\end{lem}
Le lemme est facile {\`a} d{\'e}montrer en l'absence de la condition \( \phi \lrcorner F=0 \),
car il s'agit alors simplement d'un calcul de cohomologie ; la version {\'e}crite
ici est la cons{\'e}quence d'un r{\'e}sultat plus fort dont nous aurons plus loin, le
lemme \ref{lem-92}.

Ce lemme implique imm{\'e}diatement que \( \mathcal{C}^{k} \) est une vari{\'e}t{\'e} lisse
pr{\`e}s de l'origine. Plus pr{\'e}cis{\'e}ment, notons \( P \) un inverse {\`a} droite de
l'op{\'e}rateur \( \overline{\partial }_{H} \) du lemme pr{\'e}c{\'e}dent. Comme l'op{\'e}rateur
commute {\`a} l'action de \( S^{1} \), on peut choisir \( P \) v{\'e}rifiant la m{\^e}me
propri{\'e}t{\'e}. Alors une carte locale \( \Phi :\phi \in \mathcal{C}^{k}\rightarrow \phi _{1}\in T_{0}\mathcal{C}^{k} \)
est fournie par le choix 
\begin{equation}
\label{phi-1}
\phi _{1}=\phi +P\frac{1}{2}[\phi ,\phi ],
\end{equation}
 de sorte que \( \overline{\partial }_{H}\phi _{1}=0 \). Le point important
est que si \( \phi  \) est {\`a} coefficients de Fourier positifs ou nuls, il en
est de m{\^e}me pour \( [\phi ,\phi ] \), et par cons{\'e}quent pour \( P[\phi ,\phi ] \),
donc finalement pour \( \phi _{1} \). La r{\'e}ciproque est vraie, puisqu'on retrouve
\( \phi  \) {\`a} partir de \( \phi _{1} \) par 
\[
\phi =\phi _{1}+\phi _{2}+\cdots ,\quad \phi _{i}=-P\frac{1}{2}\sum _{1\leq j<i}[\phi _{j},\phi _{i-j}].\]
 On en d{\'e}duit le corollaire suivant.

\begin{cor}
\label{cor-82}Pr{\`e}s de l'origine, l'espace \( \mathcal{C}^{k} \) est une vari{\'e}t{\'e}
(hilbertienne) lisse, d'espace tangent {\`a} l'origine 
\[
T_{0}\mathcal{C}^{k}=\{\phi \in \mathcal{H}^{k}(\Omega ^{0,1}\otimes T^{1,0}),\overline{\partial }_{H}\phi =0,\phi \lrcorner F=0\},\]
 et la carte locale \( \phi \in \mathcal{C}^{k}\rightarrow \phi _{1}\in T_{0}\mathcal{C}^{k} \)
d{\'e}finie par (\ref{phi-1}) satisfait la condition : \( \phi  \) est {\`a} coefficients
de Fourier positifs ou nuls si et seulement si \( \phi _{1} \) est {\`a} coefficients
de Fourier positifs ou nuls.\qed 
\end{cor}
Rappelons que \( \mathcal{T} \) est le bord d'un domaine de \( P^{3} \). Plus
g{\'e}n{\'e}ralement, nous regardons maintenant l'espace des structures CR \emph{remplissables},
\[
\mathcal{C}_{+}^{k}\subset \mathcal{C}^{k},\]
 c'est-{\`a}-dire des structures CR qui bordent une vari{\'e}t{\'e} complexe, d{\'e}formation
du domaine de \( P^{3} \) initial. 

\begin{thm}
\label{th-CR-rem}\( \mathcal{C}_{+}^{k} \) est une sous-vari{\'e}t{\'e} de \( \mathcal{C}^{k} \),
dont l'espace tangent {\`a} l'origine est constitu{\'e} des \( \phi _{1}\in T_{0}\mathcal{C}^{k} \),
tels que dans la d{\'e}composition 
\[
\phi _{1}=\sum _{\rho }\left\langle w,\rho (g^{-1})v\right\rangle ,\]
 on ait
\begin{itemize}
\item si \( V_{\rho }=V^{K,K+2} \), alors \( w_{\overline{1}}^{1}(-K+2,K+2)=0 \)
;
\item si \( V_{\rho }=V^{K,L} \) avec \( K\leq L-4 \), alors 
\[
w_{\overline{1}}^{1}(-K,K+4)=w_{\overline{1}}^{2}(-K,K+4)=w_{\overline{2}}^{1}(-K,K+4)=0.\]

\end{itemize}
\end{thm}
\begin{proof}
Par le lemme \ref{lem-64}, on sait que \( \phi  \) est remplissable si et
seulement si la jauge de Coulomb \( \varphi ^{*}\phi  \) est {\`a} coefficients
de Fourier positifs ou nuls ; par cons{\'e}quent, \( \mathcal{C}_{+}^{k} \) est
le lieu des z{\'e}ros de l'application \( \phi \rightarrow \pi _{<0}\varphi ^{*}\phi  \)
d{\'e}finie sur \( \mathcal{C}^{\mathcal{k}} \). Il est malais{\'e} d'analyser directement
l'image de cet op{\'e}rateur, parce que la condition non lin{\'e}aire \( \overline{\partial }_{H}\phi +\frac{1}{2}[\phi ,\phi ]=0 \)
ne donne rien apr{\`e}s projection par \( \pi _{<0} \). Comme la carte \( \Phi  \)
d{\'e}finie par (\ref{phi-1}) pr{\'e}serve la positivit{\'e} des coefficients de Fourier
par le corollaire \ref{cor-82}, on va substituer {\`a} la variable \( \phi  \)
la variable \( \phi _{1}=\Phi (\phi ) \) : pour l'action des contactomorphismes
sur \( \phi _{1} \), nous avons encore une jauge de Coulomb \( J(\phi _{1})=\varphi ^{*}\phi _{1} \)
(attention : \( \Phi  \) n'a pas de raison de pr{\'e}server la jauge de Coulomb,
donc cela donne sur \( \phi  \) une condition l{\'e}g{\`e}rement diff{\'e}rente {\`a} la jauge
de Coulomb pr{\'e}c{\'e}demment utilis{\'e}e), et \( \mathcal{C}_{+}^{k} \) appara{\^\i}t comme
le lieu des z{\'e}ros de l'application 
\[
\pi _{<0}\circ J:T_{0}\mathcal{C}^{k}\longrightarrow (T_{0}\mathcal{C}^{k})_{<0}\cap W;\]
 cet op{\'e}rateur est clairement submersif, donc \( \mathcal{C}^{k}_{+} \) est
bien une vari{\'e}t{\'e}, dont l'espace tangent {\`a} l'origine est le noyau de la lin{\'e}arisation,
d{\'e}termin{\'e} par le th{\'e}or{\`e}me \ref{th-im-dd}.
\end{proof}
En particulier, on d{\'e}duit la cons{\'e}quence suivante.

\begin{cor}
\label{cor-CR-non-rem}Il existe une famille de dimension infinie de structures
CR int{\'e}grables, non remplissables, proches de la structure standard.
\end{cor}
D'autre part, si on note \( \mathcal{C}^{k}_{-} \) l'espace des structures
CR qui sont remplissables de l'autre c{\^o}t{\'e} dans \( P^{3} \), on a une description
de leur espace tangent similaire {\`a} celle du th{\'e}or{\`e}me \ref{th-CR-rem}, avec
une condition sur les \( V^{K,L} \) tel que \( K>L \). On en d{\'e}duit que \( \mathcal{C}^{k}_{+} \)
et \( \mathcal{C}^{k}_{-} \) sont transverses, et que l'intersection \( \mathcal{C}^{k}_{+}\cap \mathcal{C}^{k}_{-} \)
est une sous-vari{\'e}t{\'e} de dimension infinie dans \( \mathcal{C}^{k}_{+} \) et
dans \( \mathcal{C}^{k}_{-} \) : on a donc le corollaire suivant.

\begin{cor}
\label{cor-CR-non-def}Il existe une famille de dimension infinie de structures
CR remplissables, qui ne sont pas obtenues en d{\'e}formant \( \mathcal{T} \) comme
une hypersurface r{\'e}elle dans \( P^{3} \).
\end{cor}

\section{\label{sec9}Espace des bords de m{\'e}triques autoduales}

Nous appliquons {\`a} pr{\'e}sent les m{\^e}mes id{\'e}es {\`a} la description des couples \( [g,\mathcal{Q}] \)
sur la sph{\`e}re \( S^{3} \) qui sont les bords de m{\'e}triques autoduales sur la
boule \( B^{4} \). Nous notons \( \mathcal{M} \) l'espace des m{\'e}triques conformes
sur \( S^{3} \), \( B \) l'espace des couples \( [g,\mathcal{Q}] \) sur \( S^{3} \)
; pour la r{\'e}gularit{\'e}, on notera \( \mathcal{M}^{k} \) et \( \mathcal{B}^{k} \)
les objets qui deviennent de r{\'e}gularit{\'e} \( \mathcal{H}^{k} \) une fois remont{\'e}s
{\`a} l'espace des twisteurs par la construction de la section \ref{sec11}. On
consid{\'e}rera le groupe des diff{\'e}omorphismes de \( S^{3} \) comme un sous-groupe
\( \Diff ^{k+1}S^{3}\subset \mathcal{G}^{k+1} \) du groupe des contactomorphismes
de \( \mathcal{T} \), et on consid{\'e}rera le groupe r{\'e}duit \( \Diff ^{k+1}_{0}S^{3}\subset \mathcal{G}^{k+1}_{0} \)
constitu{\'e} des diff{\'e}omorphismes qui fixent deux points sur l'espace des twisteurs,
de sorte que 
\[
\Diff ^{k+1}S^{3}/\Diff ^{k+1}_{0}S^{3}=SO_{1,4},\]
 le groupe conforme de \( S^{3} \) (qui est aussi le groupe des automorphismes
CR de \( \mathcal{T} \)).

\subsection{\label{sec91}Objets sur \protect\( S^{3}\protect \) et sur \protect\( \mathcal{T}\protect \) }

On peut voir la sph{\`e}re \( S^{3} \) comme un espace homog{\`e}ne sous son groupe
d'isom{\'e}tries \( SO_{4}=Sp_{1}Sp_{1} \) : 
\[
S^{3}=Sp_{1}Sp_{1}/SO_{3},\]
 o{\`u} le \( SO_{3}\hookrightarrow Sp_{1}Sp_{1} \) provient de l'inclusion diagonale
\( Sp_{1}\hookrightarrow Sp_{1}\times Sp_{1} \) ; de ce point de vue, on peut
faire des d{\'e}compositions des objets sur les repr{\'e}sentations de \( Sp_{1}Sp_{1} \). 

En particulier, les couples \( [g,\mathcal{Q}] \) donnent lieu {\`a} la donn{\'e}e
infinit{\'e}simale de \( (\dot{g},\dot{\mathcal{Q}}) \), o{\`u} \( \dot{g} \) et \( \dot{\mathcal{Q}} \)
sont des formes quadratiques {\`a} trace nulle, donc repr{\'e}sent{\'e}es par des 
\[
\sum _{\rho }\left\langle w,\rho (g^{-1})v\right\rangle ,\quad v\otimes w\in [V_{\rho }]\otimes \hom _{SO_{3}}([V_{\rho }],Sym_{0}^{2}\mathbf{R}^{3});\]
 ici, le crochet \( [V_{\rho }] \) repr{\'e}sente une forme r{\'e}elle de \( V_{\rho } \).
Toutes les repr{\'e}sentations de \( Sp_{1}Sp_{1} \) n'interviennent pas : en effet,
si \( V_{\rho }=\sigma ^{K}\otimes \sigma ^{L}, \) la repr{\'e}sentation diagonale
de \( Sp_{1}Sp_{1} \) sur \( V_{\rho } \) est 
\begin{equation}
\label{sKsL}
\sigma ^{K+L}\oplus \sigma ^{K+L-2}\oplus \cdots \sigma ^{|K-L|};
\end{equation}
 comme \( Sym_{0}^{2}\mathbf{C}^{3}=\sigma ^{4} \), seules interviennent les
repr{\'e}sentations avec 
\[
|K-L|=0,2\textrm{ ou }4;\]
 puisque dans la d{\'e}composition (\ref{sKsL}) la repr{\'e}sentation \( \sigma ^{4} \)
n'intervient qu'une fois, l'espace \( \hom _{SO_{3}}([V_{\rho }],Sym_{0}^{2}\mathbf{R}^{3}) \)
est de dimension 1.

Les champs de vecteurs se d{\'e}composent de la m{\^e}me mani{\`e}re, en utilisant la repr{\'e}sentation
\( \mathbf{R}^{3}=[\sigma ^{2}] \) de \( SO_{3} \) : compte tenu {\`a} nouveau
de la formule (\ref{sKsL}), les seules repr{\'e}sentations qui interviennent sont
{\`a} \( |K-L|=0 \) ou \( 2 \) ; il est facile de voir que l'action infinit{\'e}simale
des champs de vecteurs sur les classes conformes de m{\'e}triques, 
\[
[V_{\rho }]\otimes \hom _{SO_{3}}([V_{\rho }],\mathbf{R}^{3})\longrightarrow [V_{\rho }]\otimes \hom _{SO_{3}}([V_{\rho }],Sym_{0}^{2}\mathbf{R}^{3}),\]
 est surjective pour \( |K-L|=0 \) et \( 2 \) ; de l{\`a} on d{\'e}duit le lemme suivant.

\begin{lem}
\label{lem-91}Le quotient de l'espace des classes conformes de m{\'e}triques par
le groupe r{\'e}duit des diff{\'e}omorphismes de la sph{\`e}re, \( \mathcal{M}^{k}/\Diff _{0}^{k+1}S^{3} \),
est (pr{\`e}s de la m{\'e}trique standard) une vari{\'e}t{\'e}, d'espace tangent {\`a} l'origine
constitu{\'e} des \( \dot{g}\in \mathcal{H}^{k} \) tels que 
\[
\dot{g}\in \sum _{K-L=\pm 4}[V^{K,L}]\otimes \hom _{SO_{3}}([V^{K,L}],Sym_{0}^{2}\mathbf{R}^{3}).\]
 
\end{lem}
Par la construction twistorielle de la section \ref{sec11}, la donn{\'e}e de \( [g,\mathcal{Q}] \)
est {\'e}quivalente {\`a} la donn{\'e}e d'une structure CR int{\'e}grable \( \phi  \) sur \( \mathcal{T} \),
compatible {\`a} la structure r{\'e}elle (\( \tau ^{*}\overline{\phi }=\phi  \)), et
pr{\'e}servant les \( P^{1} \) r{\'e}els ; notons que \( \phi  \) pr{\'e}serve les \( P^{1} \)
si sa projection 
\[
\pi _{P^{1}}\phi \in \Omega ^{0,1}P^{1}\otimes (T^{1,0}/T^{1,0}P^{1})\]
 s'annule. Infinit{\'e}simalement, cela se traduit par le fait qu'est associ{\'e} {\`a}
\( (\dot{g},\dot{\mathcal{Q}}) \) un tenseur \( \phi _{1} \) satisfaisant
les {\'e}quations 
\[
\overline{\partial }_{H}\phi _{1}=0,\quad \phi _{1}\lrcorner F=0,\quad \tau ^{*}\overline{\phi _{1}}=\phi _{1},\quad \pi _{P^{1}}\phi _{1}=0;\]
 ces conditions ont {\'e}t{\'e} analys{\'e}es dans les formules (\ref{dbar-phi}),(\ref{Rphi}),(\ref{tau*phi})
et (\ref{phiP1}) ; il n'est pas difficile de voir que ces conditions ne peuvent
{\^e}tre remplies que par des objets concentr{\'e}s sur les repr{\'e}sentations {\`a} \( |K-L|\leq 4 \)
; r{\'e}ciproquement, ces conditions imposent que \( \phi _{1} \) provienne d'un
couple \( (\dot{g},\dot{\mathcal{Q}}) \), mais nous n'aurons pas besoin d'{\'e}crire
explicitement l'application \( (\dot{g},\dot{\mathcal{Q}})\rightarrow \phi _{1} \).

\subsection{Remplissage}

Les probl{\`e}mes de remplissage de \( [g,\mathcal{Q}] \) et de la structure CR
associ{\'e}e \( \phi  \) sont {\'e}quivalents. Dans cette description, nous figeons
la fibration legendrienne en \( P^{1} \) ; cependant, il n'est pas possible
de se restreindre aux contactomorphismes fixant cette fibration (c'est-{\`a}-dire
provenant de diff{\'e}omorphismes de \( S^{3} \)), car le contactomorphisme mettant
\( \phi  \) en jauge de Coulomb n'a pas de raison de la pr{\'e}server.

Pour surmonter cette difficult{\'e}, il est utile de repr{\'e}senter \( [g,\mathcal{Q}] \)
comme la donn{\'e}e de \( (\phi ,\alpha ) \), o{\`u} \( \phi  \) est une structure
CR int{\'e}grable, compatible {\`a} \( \tau  \), et \( \alpha  \) un contactomorphisme,
\( \tau  \)-invariant, tel que \( \alpha ^{*}\phi  \) pr{\'e}serve la fibration
en \( P^{1} \) ; bien s{\^u}r, la donn{\'e}e n'est d{\'e}finie qu'{\`a} contactomorphisme pr{\`e}s.
Maintenant, \( [g,\mathcal{Q}] \) est remplissable s'il existe un repr{\'e}sentant
\( (\phi ,\alpha ) \) avec \( \phi  \) {\`a} coefficients de Fourier positifs
ou nuls ; on peut toujours choisir pour \( \phi  \) la jauge de Coulomb.

\begin{lem}
\label{lem-92}L'op{\'e}rateur \( (\phi ,f)\rightarrow (\overline{\partial }_{H}\phi ,\pi _{P^{1}}(\phi -\overline{\partial }_{H}\#\overline{\partial }f)) \)
entre les espaces 
\[
\{\phi \in \mathcal{H}^{k},\phi \lrcorner F=0\}\times \mathcal{H}^{k+2}(\mathbf{R})\longrightarrow \mathcal{H}^{k-1}\times \mathcal{H}^{k}\]
 est surjectif. De plus, si \( \overline{\partial }_{H}\phi  \) est {\`a} coefficients
de Fourier positifs ou nuls, alors \( \phi  \) peut {\^e}tre choisi {\`a} coefficients
de Fourier positifs ou nuls.
\end{lem}
Nous renvoyons la d{\'e}monstration de ce lemme {\`a} la section \ref{sec93}. 

Comme dans la section \ref{sec8}, on d{\'e}duit du lemme que l'espace des \( (\phi ,\alpha )\in \mathcal{C}^{k}\times \mathcal{G}^{k+1} \),
tels que \( \alpha ^{*}\phi  \) pr{\'e}serve les \( P^{1} \), est une vari{\'e}t{\'e}
model{\'e}e sur les \( (\phi _{1},f_{1}) \) satisfaisant les conditions 
\begin{eqnarray}
\overline{\partial }_{H}\phi _{1} & = & 0,\nonumber \\
\phi _{1}\lrcorner F & = & 0,\label{phi1-f1} \\
\pi _{P^{1}}(\phi _{1}-\overline{\partial }_{H}\#\overline{\partial }f_{1}) & = & 0;\nonumber 
\end{eqnarray}
 surtout, on dispose d'une carte \( (\phi ,\alpha )\rightarrow (\phi _{1},f_{1}) \)
telle que la positivit{\'e} des coefficients de Fourier de \( \phi  \) et \( \phi _{1} \)
soit {\'e}quivalente. Le quotient par l'action de \( \mathcal{G}^{k+1} \) est aussi
une vari{\'e}t{\'e}, model{\'e}e sur une tranche par rapport {\`a} l'action du groupe, c'est-{\`a}-dire
les \( (\phi _{1},f_{1}) \) satisfaisant de plus la condition 
\begin{equation}
\label{phi1-f1-2}
\phi _{1}\in W,\quad f_{1}\perp \ker (\overline{\partial }_{H}\#\overline{\partial }).
\end{equation}
Jusqu'ici, on n'a rien montr{\'e}, puisqu'on a simplement d{\'e}crit localement l'espace
\( \mathcal{B}^{k} \) des couples \( [g,\mathcal{Q}] \) en termes des \( (\phi ,\alpha ) \),
mais avec l'avantage de pouvoir d{\'e}terminer, sur la donn{\'e}e infinit{\'e}simale \( (\phi _{1},f_{1}) \)
si le \( (\phi ,\alpha ) \) correspondant est remplissable. Notons 
\[
\mathcal{B}^{k}_{+}\subset \mathcal{B}^{k}\]
 l'espace des couples \( [g,\mathcal{Q}] \) remplissables par une m{\'e}trique
autoduale.

Passons au quotient par le groupe des diff{\'e}omorphismes de \( S^{3} \) ; {\`a} nouveau,
on voit que \( \mathcal{B}^{k}/{\Diff }_{0}^{k+1}S^{3} \) est une vari{\'e}t{\'e} model{\'e}e
sur les \( (\phi _{1},f_{1}) \) satisfaisant (\ref{phi1-f1}) et la condition
\begin{equation}
\label{phi1-f1-3}
\phi _{1}\in W,\quad f_{1}\perp \mathfrak {Diff}^{k+1}S^{3}.
\end{equation}
 Maintenant, puisque \( \phi _{1} \) est en jauge de Coulomb, \( (\phi _{1},f_{1}) \)
correspond {\`a} un couple \( [g,\mathcal{Q}] \) remplissable si et seulement si
\( \phi _{1} \) est {\`a} coefficients de Fourier positifs ou nuls : on en d{\'e}duit
que \( \mathcal{B}^{k}_{+}/{\Diff }_{0}^{k+1}S^{3} \) est une vari{\'e}t{\'e} lisse
model{\'e}e sur les \( (\phi _{1},f_{1}) \) satisfaisant (\ref{phi1-f1}) et (\ref{phi1-f1-3}),
et tels que 
\[
\pi _{<0}\phi _{1}=0.\]

Il ne reste plus qu'{\`a} traduire ces conditions en termes de \( [g,\mathcal{Q}] \).
On d{\'e}duit le r{\'e}sultat suivant.

\begin{thm}
\label{th-gQ-rem}L'espace \( \mathcal{B}^{k}_{+}/\Diff ^{k+1}_{0}S^{3} \)
des \( [g,\mathcal{Q}] \) remplissables par une m{\'e}trique autoduale sur \( B^{4} \),
modulo le groupe r{\'e}duit des diff{\'e}omorphismes, est pr{\`e}s de la m{\'e}trique standard
une vari{\'e}t{\'e}, dont l'espace tangent en la m{\'e}trique standard est constitu{\'e} des
\( (\dot{g},\dot{\mathcal{Q}})\in \mathcal{H}^{k} \) tels que 
\begin{eqnarray*}
\dot{g} & \in  & \sum _{L\geq 0}[V^{L+4,L}]\otimes \hom _{SO_{3}}([V^{L+4,L}],Sym_{0}^{2}\mathbf{R}^{3}),\\
\dot{\mathcal{Q}} & \in  & \sum _{(K=L,L+2,L+4,L\geq 0)}[V^{K,L}]\otimes \hom _{SO_{3}}([V^{K,L}],Sym_{0}^{2}\mathbf{R}^{3}).
\end{eqnarray*}

\end{thm}
\begin{rem}
L'espace \( \mathcal{B}^{k}_{-}/\Diff _{0}^{k+1}S^{3} \) des \( [g,\mathcal{Q}] \)
remplissables par une m{\'e}trique antiautoduale est, de mani{\`e}re similaire, une
vari{\'e}t{\'e} d'espace tangent les \( (\dot{g},\dot{\mathcal{Q}}) \) tel que \( \dot{g} \)
se d{\'e}compose seulement sur les \( V^{K,K+4} \) et \( \dot{\mathcal{Q}} \)
sur les \( V^{K,K} \), \( V^{K,K+2} \) et \( V^{K,K+4} \) ; on observe donc
que \( \mathcal{B}^{k}_{+}/\Diff _{0}^{k+1}S^{3} \) et \( \mathcal{B}^{k}_{-}/\Diff _{0}^{k+1}S^{3} \)
sont transverses dans \( \mathcal{B}^{k}/\Diff _{0}^{k+1}S^{3} \), et l'intersection
est une vari{\'e}t{\'e} dont l'espace tangent est constitu{\'e} des \( \dot{\mathcal{Q}} \)
qui se d{\'e}composent uniquement sur les \( V^{K,K} \). Cela nous indique que
cette intersection est en r{\'e}alit{\'e} constitu{\'e}e des d{\'e}formations de \( S^{3} \)
comme sous-vari{\'e}t{\'e} de \( \mathbf{R}^{4} \), et les m{\'e}triques qui remplissent
sont donc plates ; en effet, ces d{\'e}formations sont param{\'e}tr{\'e}es par \( fn \),
o{\`u} \( f \) est une fonction et \( n \) le vecteur normal {\`a} \( S^{3} \) ;
infinit{\'e}simalement, elles laissent la classe conforme invariante et ne modifient
donc que la seconde forme fondamentale, ce qui correspond bien {\`a} notre description
infinit{\'e}simale. Cela donne le th{\'e}or{\`e}me \ref{th-C0} comme cons{\'e}quence du th{\'e}or{\`e}me
\ref{th-gQ-rem}.
\end{rem}
\begin{proof}
Le tenseur \( \psi _{1}=\phi _{1}-\overline{\partial }_{H}\#\overline{\partial }f_{1} \)
repr{\'e}sente une donn{\'e}e infinit{\'e}simale \( (\dot{g},\dot{\mathcal{Q}}) \), et
il s'agit de voir dans quelles repr{\'e}sentations il peut vivre, sachant que \( \phi _{1} \)
a seulement des coefficients de Fourier positifs ou nuls. La r{\'e}ponse est fournie
par le th{\'e}or{\`e}me \ref{th-im-dd} : pour \( K\geq L \), on peut toujours, partant
de \( \psi _{1} \), r{\'e}soudre le probl{\`e}me 
\begin{eqnarray*}
(\overline{\partial }_{H}\#\overline{\partial }f_{1})_{<0} & = & -(\psi _{1})_{<0},\\
(\overline{\partial }_{H}\#\overline{\partial }f_{1})_{0} & = & -\pi _{\Im \overline{\partial }_{H}\#\overline{\partial }}(\psi _{1})_{0},
\end{eqnarray*}
 et on aura bien la condition (\ref{phi1-f1-3}) si \( f_{1} \) est orthogonal
aux diff{\'e}omorphismes infinit{\'e}simaux de \( S^{3} \), ce qu'on obtient en prenant
\( \psi _{1} \) lui-m{\^e}me orthogonal {\`a} l'action infinit{\'e}simale des diff{\'e}omorphismes,
c'est-{\`a}-dire \( \dot{g} \) satisfaisant la condition du lemme \ref{lem-91}.
Pour \( K=L-2 \) ou \( L-4 \), il y a des obstructions {\`a} analyser : si \( K=L-2 \),
la condition \( w_{\overline{1}}^{1}(-K+2,K+2)=0 \) du th{\'e}or{\`e}me \ref{th-im-dd}
revient {\`a} s{\'e}lectionner un sous-espace de dimension 1 dans \( \hom _{SO_{3}}([V^{K,L}],Sym^{2}_{0}\mathbf{R}^{3}\oplus Sym^{2}_{0}\mathbf{R}^{3}) \)
qui est de dimension 2 : nous connaissons d{\'e}j{\`a} ce sous-espace, il s'agit des
m{\'e}triques dans l'image infinit{\'e}simale des diff{\'e}omorphismes, et, comme ci-dessus,
nous les {\'e}liminons au quotient. Si \( K=L-4 \), les conditions \( w_{\overline{1}}^{1}(-K,K+4)=w_{\overline{1}}^{2}(-K,K+4) \)
du th{\'e}or{\`e}me \ref{th-im-dd} tuent compl{\`e}tement les deux degr{\'e}s de libert{\'e} pour
\( \psi _{1} \), donc cette composante doit {\^e}tre nulle. Finalement, on d{\'e}duit
qu'il ne reste que les composantes {\`a} \( K\geq L \), d'o{\`u} le th{\'e}or{\`e}me.
\end{proof}

\subsection{\label{sec93}D{\'e}monstration du lemme \ref{lem-92}}

Une des difficult{\'e}s de la d{\'e}monstration du lemme est le fait que l'op{\'e}rateur
\( \pi _{P^{1}} \) ne commute pas {\`a} l'action de \( S^{1} \) : l'assertion
sur les coefficients de Fourier ne sera donc pas cons{\'e}quence de la surjectivit{\'e}
de l'op{\'e}rateur.

On va raisonner dans chaque repr{\'e}sentation \( V^{K,L} \), pour rendre le probl{\`e}me
alg{\'e}brique ; par abus de notation, on confondra \( \phi  \) et \( w \) tel
que \( \phi =\left\langle w,\rho (g^{-1})v\right\rangle  \), donc \( \phi  \)
sera consid{\'e}r{\'e} comme {\'e}l{\'e}ment de \( \hom _{U_{1}}(V^{K,L},\mathfrak {m}^{0,1}\otimes \mathfrak {m}_{1,0}) \)
; on fera la m{\^e}me chose pour les autres tenseurs en jeu.

\begin{figure}[hbt]
{\par\centering \includegraphics{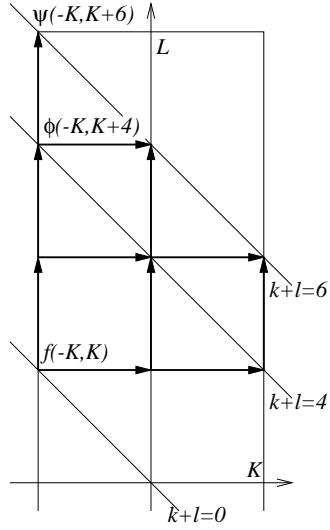} \par}

\caption{\label{fig-3}cas \protect\( K\leq L-6\protect \) }
\end{figure}

\noindent \textbf{Premier cas :} \( K\leq L-6 \). Alors l'op{\'e}rateur 
\begin{equation}
\label{phi-seul}
\phi \textrm{ tel que }\phi \lrcorner F=0\, \longrightarrow (\overline{\partial }_{H}\phi ,\pi _{P^{1}}\phi )
\end{equation}
 est un isomorphisme ; l'inverse est obtenu en regardant la figure \ref{fig-3}
: une solution de
\[
\phi \lrcorner F=0,\quad \overline{\partial }_{H}\phi =\psi ,\quad \pi _{P^{1}}\phi =g,\]
est obtenue de proche en proche ; d'apr{\`e}s (\ref{dbar-phi}), \( \phi _{\overline{1}}(-K,K+4) \)
est d{\'e}termin{\'e} par \( \psi (-K,K+6) \) ; les {\'e}quations \( \phi \lrcorner F=0 \)
et \( \pi _{P^{1}}\phi =g \) se traduisent, d'apr{\`e}s (\ref{phiP1}) et (\ref{Rphi}),
par 
\begin{eqnarray*}
\phi _{\overline{1}}^{1}+\phi _{\overline{2}}^{1}-\phi _{\overline{1}}^{2}-\phi _{\overline{2}}^{2} & = & g,\\
\phi _{\overline{2}}^{1}+\phi ^{2}_{\overline{1}} & = & 0;
\end{eqnarray*}
 donc \( \phi _{\overline{2}}(-K,K+4) \) est d{\'e}termin{\'e} par \( \phi _{\overline{1}}(-K,K+4) \)
et \( g(-K,K+4) \) ; on poursuit alors le raisonnement de proche en proche
pour \( \phi (-K+2,K+2) \),...

On d{\'e}duit que l'op{\'e}rateur du lemme \ref{lem-92} est bien surjectif ; nous devons
maintenant prouver l'{\'e}nonc{\'e} concernant les coefficients de Fourier. Supposons
que \( \psi  \) n'ait que des coefficients de Fourier positifs ou nuls, et
tentons de r{\'e}soudre le probl{\`e}me de sorte que \( \phi  \) v{\'e}rifie la m{\^e}me propri{\'e}t{\'e}.
Commen{\c c}ons par r{\'e}soudre (de proche en proche, comme ci-dessus pour \( \phi  \)),
\[
\pi _{k<0}\pi _{P^{1}}(-\overline{\partial }_{H}\#\overline{\partial }f)=\pi _{k<0}g;\]
 comme \( f \) est r{\'e}elle, cela d{\'e}termine compl{\`e}tement \( f \), {\`a} l'exception
du terme \( f(0,0) \) que nous prenons nul ; notons 
\[
g'=g+\pi _{P^{1}}(-\overline{\partial }_{H}\#\overline{\partial }f),\]
de sorte que \( g'(k,-k)=0 \) si \( k<0 \), et r{\'e}solvons le probl{\`e}me 
\[
\phi \lrcorner F=0,\quad \overline{\partial }_{H}\phi =\psi ,\quad \pi _{P^{1}}\phi =g';\]
 j'affirme que \( \phi  \) a ses coefficients de Fourier n{\'e}gatifs triviaux.
En effet, puisque la solution est calcul{\'e}e de proche en proche, il est clair
que \( \phi (k,-k+4)=0 \) si \( k<0 \) ; compte tenu des poids de l'action
de \( S^{1} \) calcul{\'e}s en (\ref{poids-S1}) et de la positivit{\'e} des coefficients,
on sait que 
\[
\psi _{\overline{1}\overline{2}}^{1}(0,6)=\psi _{\overline{1}\overline{2}}^{2}(0,6)=\psi _{\overline{1}\overline{2}}^{1}(2,4)=0;\]
 testant {\`a} nouveau le diagramme permettant de calculer \( \phi  \), nous d{\'e}duisons
tout d'abord 
\[
\phi _{\overline{1}}^{1}(0,4)=\phi _{\overline{1}}^{2}(0,4)=\phi _{\overline{2}}^{1}(0,4)=0;\]
 finalement, dans l'{\'e}quation \( \overline{\partial }\phi =\psi  \), le terme
\( \psi _{\overline{1}\overline{2}}^{1}(2,4) \) est d{\'e}termin{\'e} par \( \phi _{\overline{2}}^{1}(0,4) \)
et \( \phi _{\overline{1}}^{1}(2,2) \), donc son annulation implique 
\[
\phi _{\overline{1}}^{1}(2,2)=0\]
 qui {\'e}tait le dernier terme strictement n{\'e}gatif {\`a} annuler. 

\noindent \textbf{Second cas :} \( K\geq L-4 \). Dans ce cas, le noyau de l'op{\'e}rateur
(\ref{phi-seul}) est constitu{\'e} des \( \phi  \) tels que 
\begin{equation}
\label{phi-B}
\phi \lrcorner F=0,\quad \overline{\partial }_{H}\phi =0,\quad \pi _{P^{1}}\phi =0,
\end{equation}
qui correspondent exactement aux donn{\'e}es infinit{\'e}simales \( (\dot{g},\dot{\mathcal{Q}}) \)
dans \( \mathcal{B} \) ; ce noyau est donc nul pour \( K>L+4, \) et de dimension
2 pour \( L-4\leq K\leq L+4 \) ; un d{\'e}compte de dimension nous indique alors
que l'op{\'e}rateur (\ref{phi-seul}) est surjectif. 

Reste {\`a} montrer l'assertion sur les coefficients de Fourier positifs : on utilise
maintenant que le noyau de \( (\phi ,f)\rightarrow (\overline{\partial }_{H}\phi ,\pi _{P^{1}}(\phi -\overline{\partial }_{H}\#\overline{\partial }f)) \)
contient les 
\begin{equation}
\label{ker-phi-f}
(\phi =\overline{\partial }_{H}\#\overline{\partial }f,f).
\end{equation}
 Si \( \psi  \) est donn{\'e}e, {\`a} coefficients de Fourier positifs ou nuls, prenons
une solution \( \phi  \) de
\[
\phi \lrcorner F=0,\quad \overline{\partial }_{H}\phi =\psi ,\quad \pi _{P^{1}}\phi =g,\]
on va la modifier par un {\'e}l{\'e}ment de type (\ref{ker-phi-f}) de sorte de tuer
tous les coefficients strictement n{\'e}gatifs. Les premi{\`e}res {\'e}quations commutent
{\`a} l'action de \( S^{1} \), ce qui implique 
\[
\phi _{-}\lrcorner F=0,\quad \overline{\partial }_{H}\phi _{-}=0.\]

Dans le cas \( K\geq L, \) par le th{\'e}or{\`e}me \ref{th-im-dd}, il existe une fonction
r{\'e}elle \( f \) telle que \( \overline{\partial }_{H}\#\overline{\partial }f_{-}=-\phi _{-} \)
; alors le couple \( (\phi +\overline{\partial }_{H}\#\overline{\partial }f,f) \)
est solution des {\'e}quations voulues, et \( \phi +\overline{\partial }_{H}\#\overline{\partial }f \)
est {\`a} coefficients de Fourier positifs ou nuls.

Dans les cas \( K=L-2 \) et \( K=L-4 \), on peut appliquer le m{\^e}me raisonnement,
pourvu qu'on se ram{\`e}ne aux conditions prescrites par le th{\'e}or{\`e}me \ref{th-im-dd}.
Pour cela, notons qu'un \( \phi  \) satisfaisant (\ref{phi-B}) est enti{\`e}rement
d{\'e}termin{\'e} par les deux coefficients \( \phi ^{1}_{\overline{1}}(-K,K+4) \)
et \( \phi _{\overline{1}}^{2}(-K,K+4) \), tous les autres coefficients s'en
d{\'e}duisant de proche en proche ; cela signifie qu'{\'e}tant donn{\'e}e une solution \( \phi  \)
de 
\[
\phi \lrcorner F=0,\quad \overline{\partial }_{H}\phi =\psi ,\quad \pi _{P^{1}}\phi =g,\]
 on peut la modifier par un {\'e}l{\'e}ment satisfaisant (\ref{phi-B}) de sorte que
\( \phi _{\overline{1}}(-K,K+4)=0 \) (ce qui implique aussi \( \phi _{\overline{2}}^{1}(-K,K+4)=0 \))
; les conditions d'application du th{\'e}or{\`e}me \ref{th-im-dd} sont ainsi {\'e}tablies.\qed

\section{\label{sec10}Espace des infinis conformes de m{\'e}triques autoduales d'Einstein}

Notons 
\[
\mathcal{M}_{+}^{k}=\mathcal{B}_{+}^{k}\cap \mathcal{M}^{k}\]
 l'espace des m{\'e}triques conformes sur \( S^{3} \) qui sont les infinis conformes
de m{\'e}triques d'Einstein autoduales sur \( B^{4} \). Le but de cette section
est de montrer que \( \mathcal{M}_{+}^{k} \) est une sous-vari{\'e}t{\'e} de \( \mathcal{B}_{+}^{k} \):
il en r{\'e}sultera que l'espace tangent en la m{\'e}trique ronde est constitu{\'e} des
\( \dot{g} \) satisfaisant les conditions du th{\'e}or{\`e}me \ref{th-gQ-rem}, ce
qui fournit le th{\'e}or{\`e}me suivant.

\begin{thm}
\label{th-g-rem}L'espace \( \mathcal{M}^{k}_{+}/\Diff ^{k+1}_{0}S^{3} \) des
m{\'e}triques conformes \( [g] \), infinis conformes de m{\'e}triques autoduales d'Einstein
sur \( B^{4} \), modulo le groupe r{\'e}duit des diff{\'e}omorphismes, est pr{\`e}s de
la m{\'e}trique standard une vari{\'e}t{\'e} d'espace tangent constitu{\'e} des \( \dot{g}\in \mathcal{H}^{k} \)
tels que 
\[
\dot{g}\in \sum _{L\geq 0}[V^{L+4,L}]\otimes \hom _{SO_{3}}([V^{L+4,L}],Sym_{0}^{2}\mathbf{R}^{3}).\]

\end{thm}
\begin{rem}
Bien entendu, l'espace tangent {\`a} l'espace \( \mathcal{M}^{k}_{-}/\Diff ^{k+1}_{0}S^{3} \)
des infinis conformes de m{\'e}triques antiautoduales d'Einstein est 
\[
\sum _{K\geq 0}[V^{K,K+4}]\otimes \hom _{SO_{3}}([V^{K,K+4}],Sym_{0}^{2}\mathbf{R}^{3}),\]
 de sorte que, en vue du lemme \ref{lem-91}, on a bien en la m{\'e}trique ronde
la d{\'e}composition de l'espace tangent 
\[
T(\mathcal{M}^{k}/\Diff ^{k+1}_{0}S^{3})=T(\mathcal{M}^{k}_{+}/\Diff ^{k+1}_{0}S^{3})\oplus T(\mathcal{M}^{k}_{-}/\Diff ^{k+1}_{0}S^{3}),\]
 ce qui d{\'e}montre le th{\'e}or{\`e}me \ref{th-C}. 
\begin{rem}
L'espace \( \hom _{SO_{3}}([V^{K,K+4}],Sym_{0}^{2}\mathbf{R}^{3}) \) est de
dimension 1, comme on l'a vu dans la section \ref{sec91} ; la premi{\`e}re composante
de l'espace tangent, \( [V^{4,0}] \), de dimension 5, a une interpr{\'e}tation
simple : il s'agit des m{\'e}triques invariantes {\`a} droite sur \( S^{3} \), et les
m{\'e}triques d'Einstein autoduales correspondantes ont {\'e}t{\'e} explicit{\'e}es par Hitchin
\cite{Hit:95} ; la seconde composante, \( [V^{5,1}], \) de dimension 12, n'a
plus d'interpr{\'e}tation aussi claire.
\end{rem}
\end{rem}
Le th{\'e}or{\`e}me \ref{th-g-rem} est loin d'{\^e}tre une simple cons{\'e}quence du th{\'e}or{\`e}me
\ref{th-gQ-rem} : ce dernier fournit en effet un algorithme fabriquant {\`a} partir
d'une donn{\'e}e infinit{\'e}simale \( (\dot{g},\dot{\mathcal{Q}}) \) un couple \( [g,\mathcal{Q}] \)
qui est le bord d'une m{\'e}trique autoduale, mais il n'est pas clair que \( \dot{\mathcal{Q}}=0 \)
implique \( \mathcal{Q}=0 \). De plus, il ne semble pas y avoir de raison a
priori pour que l'intersection de \( \mathcal{B}_{+}^{k} \) et de \( \mathcal{M}^{k} \)
soit transverse. Nous allons employer ici une m{\'e}thode diff{\'e}rente, {\'e}tudiant les
d{\'e}formations d'une structure CR avec une structure de contact holomorphe ; on
va donc oublier le contactomorphisme et la fibration en \( P^{1} \) et montrer
que la structure est assez riche pour permettre de les r{\'e}cup{\'e}rer a posteriori.

\subsection{Th{\'e}orie de la d{\'e}formation}

Nous cherchons donc une d{\'e}formation \( \phi  \) de la structure CR, et une
forme de contact holomorphe \( \varpi ^{c} \) {\`a} valeurs dans \( L=\mathcal{O}(0,2) \)
; comme dans la section \ref{sec32}, puisque \( \varpi ^{c} \) est une petite
d{\'e}formation de \( \eta ^{c} \), on la cherche sous la forme 
\[
\varpi ^{c}=(\eta ^{c}+\xi )-\phi \lrcorner (\eta ^{c}+\xi ),\quad \xi \in \Omega '\otimes L,\quad \xi |_{T^{1,0}P_{2}^{1}}=0,\]
 et satisfaisant l'{\'e}quation (\ref{oc-hol}), c'est-{\`a}-dire 
\[
\overline{\partial }\varpi ^{c}+\phi \lrcorner \partial \varpi ^{c}+\alpha \wedge \varpi ^{c}=0\]
pour une (0,1)-forme \( \alpha  \). Cette {\'e}quation fait intervenir une connexion
sur \( L \), que nous choisissons {\'e}gale {\`a} la connexion provenant de la connexion
triviale sur \( \mathcal{O}_{P^{1}_{1}} \) et de la connexion standard sur
\( \mathcal{O}_{P^{1}_{2}}(2) \), de sorte que sa courbure est {\'e}gale {\`a} 
\[
F^{L}=-2i\omega _{2}.\]
 Si on a une solution \( (\phi ,\varpi ^{c}) \) des {\'e}quations, alors \( L \),
identifi{\'e} au quotient \( T'_{_{\phi }}/\ker \varpi ^{c} \), acquiert une structure
holomorphe, dont l'op{\'e}rateur \( \overline{\partial } \) s'{\'e}crit 
\[
\overline{\partial }+\phi \lrcorner \partial +\alpha ;\]
 l'int{\'e}grabilit{\'e} de cette structure holomorphe se traduit par l'{\'e}quation 
\begin{equation}
\label{eq-alpha}
\overline{\partial }\alpha +\phi \lrcorner (F^{L}+\partial \alpha )=0.
\end{equation}
 Rappelons d'autre part que, d'apr{\`e}s (\ref{etac+xi-hol0}), les {\'e}quations {\`a}
satisfaire par \( (\phi ,\varpi ^{c}) \) sont, en plus de la condition alg{\'e}brique
\( \phi \lrcorner F=0 \), 

\begin{eqnarray}
\overline{\partial }_{H}\phi +\frac{1}{2}[\phi ,\phi ] & = & 0,\label{eq-phi} \\
\overline{\partial }\xi -\partial '(\phi \lrcorner (\eta ^{c}+\xi ))+\phi \lrcorner \partial '(\eta ^{c}+\xi )+\alpha (\eta ^{c}+\xi ) & = & 0.\label{eq-xi} 
\end{eqnarray}
 Les trois {\'e}quations pr{\'e}c{\'e}dentes peuvent s'{\'e}crire \( P(\phi ,\alpha ,\xi )+Q(\phi ,\alpha ,\xi )=0, \)
o{\`u} \( P \) est la partie lin{\'e}aire des trois {\'e}quations pr{\'e}c{\'e}dentes, donn{\'e}e par
\( P(\phi ,\alpha ,\xi )=(\psi ,\beta ,\zeta ) \), avec 
\begin{eqnarray}
\psi  & = & \overline{\partial }_{H}\phi \nonumber \\
\beta  & = & \overline{\partial }\alpha +\phi \lrcorner F^{L}\label{for-P} \\
\zeta  & = & \overline{\partial }\xi -\partial '(\phi \lrcorner \eta ^{c})+\phi \lrcorner d\eta ^{c}+\alpha \eta ^{c}\nonumber 
\end{eqnarray}
 et \( Q \) est la partie quadratique restante.

D'apr{\`e}s la d{\'e}monstration du th{\'e}or{\`e}me \ref{th-gQ-rem}, on sait que, {\'e}tant donn{\'e}e
une m{\'e}trique infinit{\'e}simale \( \dot{g} \) satisfaisant la condition du th{\'e}or{\`e}me
\ref{th-g-rem}, il existe un contactomorphisme infinit{\'e}simal de l'espace des
twisteurs \( \mathcal{T} \) qui envoie la structure CR infinit{\'e}simale sur un
tenseur \( \phi _{1} \) {\`a} coefficients de Fourier positifs ; bien entendu,
d'apr{\`e}s le th{\'e}or{\`e}me \ref{th-fill-etac}, le contactomorphisme infinit{\'e}simal
envoie la structure de contact holomorphe sur un \( \xi _{1} \) {\`a} coefficients
positifs si bien qu'on se retrouve avec un triplet \( (\phi _{1},\xi _{1},\alpha _{1}) \),
{\`a} coefficients de Fourier positifs, solution des {\'e}quations infinit{\'e}simales 
\[
P(\phi _{1},\alpha _{1},\xi _{1})=0.\]
 On va maintenant montrer qu'{\`a} partir d'une telle donn{\'e}e infinit{\'e}simale, on
peut produire une solution, {\`a} coefficients de Fourier positifs, des {\'e}quations
(\ref{eq-alpha}), (\ref{eq-phi}) et (\ref{eq-xi}). L'{\'e}tape essentielle est
le lemme suivant.

\begin{lem}
\label{pas-de-H2}L'op{\'e}rateur \( P \), d{\'e}fini sur les triplets 
\[
(\phi ,\alpha ,\xi )\in \mathcal{H}^{k}(\Omega ^{0,1}\otimes T^{1,0})\times \mathcal{H}^{k+1}(\Omega ^{0,1})\times \mathcal{H}^{k-1}(\Omega ')\]
 tels que \( \phi \lrcorner F=0 \) et \( \xi |_{T^{1,0}P^{1}_{2}}=0 \), a
exactement pour image les triplets 
\[
(\psi ,\beta ,\zeta )\in \mathcal{H}^{k-1}(\Omega ^{0,2}\otimes T^{1,0})\times \mathcal{H}^{k}(\Omega ^{0,2})\times \mathcal{H}^{k-2}(\tilde{\Omega }^{1,1}),\]
 satisfaisant l'{\'e}quation 
\begin{equation}
\label{int-formelle}
\overline{\partial }\zeta =(\partial '\psi )\lrcorner \eta ^{c}+\beta \wedge \eta ^{c}.
\end{equation}

\end{lem}
Ce lemme est a priori surprenant, car on sait qu'avec une forme de Levi de signature
(1,1), les \( H^{1} \) sont de dimension infinie, donc on ne s'attend pas {\`a}
pouvoir atteindre tous les \( \zeta  \). La d{\'e}monstration sera faite dans la
section \ref{sec102}.

Revenons {\`a} notre th{\'e}orie de la d{\'e}formation : on est dans le cadre de la th{\'e}orie
standard de d{\'e}formation, et le lemme signifie que le \( H^{2} \) du complexe
de d{\'e}formation est nul, autrement dit il n'y a pas d'obstruction pour produire
une solution \( (\phi ,\alpha ,\xi ) \) du syst{\`e}me d'{\'e}quations (\ref{eq-alpha}),
(\ref{eq-phi}) et (\ref{eq-xi}) {\`a} partir d'une donn{\'e}e infinit{\'e}simale \( (\phi _{1},\alpha _{1},\xi _{1}) \).
Compte tenu de la forme des {\'e}quations (elles pr{\'e}servent la positivit{\'e} des coefficients
de Fourier), la solution \( (\phi ,\alpha ,\xi ) \) est {\`a} coefficients positifs
en m{\^e}me temps que \( (\phi _{1},\alpha _{1},\xi _{1}) \) : on a ainsi produit
une structure CR et une forme de contact holomorphe qui s'{\'e}tendent {\`a} l'int{\'e}rieur.
Si la donn{\'e}e initiale est compatible {\`a} la structure r{\'e}elle, il en est de m{\^e}me
pour \( (\phi ,\alpha ,\xi ) \).

Le seul point {\`a} traiter concerne la fibration en \( P^{1} \) que nous avons
a priori perdue par cette m{\'e}thode. La m{\'e}thode la plus naturelle pour la r{\'e}cup{\'e}rer
consiste {\`a} am{\'e}liorer la th{\'e}orie de la d{\'e}formation pr{\'e}c{\'e}dente pour assurer que
la forme de contact \( \varpi ^{c} \) satisfait de plus \( F\wedge \varpi ^{c}\wedge \overline{\varpi ^{c}}=0 \)
; en effet, si cette condition est satisfaite, il n'est pas difficile de voir
que 
\[
\Re (\ker (\varpi ^{c})\cap T^{1,0}_{\phi }+\overline{\ker (\varpi ^{c})\cap T^{1,0}_{\phi }})\]
 est une distribution int{\'e}grable et legendrienne ; par le th{\'e}or{\`e}me de stabilit{\'e}
de Reeb, une perturbation d'un feuilletage en sph{\`e}res \( S^{2} \) demeure une
fibration en sph{\`e}res \( S^{2} \), donc on obtient une fibration legendrienne
en sph{\`e}res \( S^{2} \), qui est n{\'e}cessairement isomorphe {\`a} la fibration initiale
par un contactomorphisme.

Cependant, au lieu de revenir sur la th{\'e}orie de la d{\'e}formation, nous pouvons
utiliser un argument plus rapide : l'extension {\`a} l'int{\'e}rieur de \( (\phi ,\alpha ,\xi ) \)
d{\'e}finit un domaine complexe \( D \) muni d'une forme de contact holomorphe
\( \varpi ^{c} \) ; on peut aussi {\'e}largir l{\'e}g{\`e}rement ce domaine, ainsi que
la forme de contact holomorphe, de l'autre c{\^o}t{\'e} de son bord \( \partial D \),
pour obtenir un domaine l{\'e}g{\`e}rement plus grand \( D' \) ; {\`a} pr{\'e}sent, la fibration
initiale en courbes rationnelles r{\'e}elles {\`a} fibr{\'e} normal \( \mathcal{O}(1)\oplus \mathcal{O}(1) \),
par le th{\'e}or{\`e}me de Kodaira, peut elle aussi {\^e}tre perturb{\'e}e, et cette famille
de \( P^{1} \) est donc l'espace des param{\`e}tres d'une vari{\'e}t{\'e} autoduale d'Einstein
\( X \), dont le bord \( \partial X \) est caract{\'e}ris{\'e} comme la famille des
\( P^{1} \) r{\'e}els qui sont tangents {\`a} la distribution de contact holomorphe
;  LeBrun \cite[proposition 1]{LeB:91} a d{\'e}montr{\'e} cette {\'e}quation est transverse,
ce qui signifie que le bord est toujours lisse, et que l'on peut le trouver
en d{\'e}formant le bord initial de la m{\'e}trique hyperbolique r{\'e}elle. Cette construction
permet donc de r{\'e}cup{\'e}rer {\`a} partir de \( (\phi ,\alpha ,\xi ) \) le bord \( \partial X \),
qui porte la m{\'e}trique conforme que nous voulions construire. Cela ach{\`e}ve la
d{\'e}monstration du th{\'e}or{\`e}me \ref{th-g-rem}.\qed

\subsection{\label{sec102}D{\'e}monstration du lemme \ref{pas-de-H2}}

On commence par se d{\'e}barrasser de \( \psi  \) et de \( \zeta  \) par le lemme
suivant.

\begin{lem}
L'op{\'e}rateur \( (\phi ,\alpha )\rightarrow (\overline{\partial }_{H}\phi ,\overline{\partial }\alpha +\phi \lrcorner F^{L}) \),
d{\'e}fini sur les \( (\phi ,\alpha ) \) tels que \( \phi \lrcorner F=0 \), est
surjectif.
\end{lem}
\begin{proof}
L'op{\'e}rateur \( \phi \rightarrow \overline{\partial }_{H}\phi  \) est surjectif
par le lemme \ref{lem-92}, donc il suffit de v{\'e}rifier que \( \alpha \rightarrow \overline{\partial }\alpha  \)
est surjectif, c'est-{\`a}-dire de montrer que \( H^{2}(\mathcal{O})=0 \). On peut
d{\'e}composer cette cohomologie suivant le poids \( k \) de l'action de \( S^{1} \)
comme dans la d{\'e}monstration du lemme \ref{lem-31}, et on d{\'e}duit 
\begin{eqnarray}
H^{2}(\mathcal{O}) & = & \oplus _{k}H^{2}(P_{1}^{1}\times P_{2}^{1},\mathcal{O}(k,-k))\nonumber \\
 & = & \oplus _{k}H^{1}(P_{1}^{1},\mathcal{O}(k))\otimes H^{1}(P_{2}^{1},\mathcal{O}(-k))=0.\label{cal-H1O} 
\end{eqnarray}

\end{proof}
De ce lemme, on d{\'e}duit que pour calculer l'image de l'op{\'e}rateur \( P \) d{\'e}fini
par (\ref{for-P}), on peut se contenter de regarder l'op{\'e}rateur \( (\phi ,\alpha ,\xi )\rightarrow \zeta  \).

Commen{\c c}ons par {\'e}tudier le simple probl{\`e}me de cohomologie \( \zeta =\overline{\partial }\xi  \).
Rappelons que \( \xi  \) est une section de 
\[
\Omega '\otimes L=(\mathbf{C}\eta \oplus \Omega ^{1,0})\otimes L,\]
 et l'op{\'e}rateur \( \overline{\partial } \) de \( \Omega ' \) a {\'e}t{\'e} calcul{\'e}
dans la formule (\ref{db-O'}), {\`a} savoir, dans la d{\'e}composition pr{\'e}c{\'e}dente de
\( \Omega ' \), 
\[
\overline{\partial }(\xi ^{0}\eta \oplus \xi ^{1,0})=\overline{\partial }\xi ^{0}\wedge \eta \oplus (\overline{\partial }\xi ^{1,0}+\xi ^{0}(\omega _{1}-\omega _{2}));\]
 ici, \( \xi ^{0} \) est une section de \( L=\mathcal{O}(0,2) \) ; {\`a} nouveau,
la cohomologie se d{\'e}compose suivant les poids \( k \) de l'action de \( S^{1} \)
et on notera \( H_{\geq 0} \) la cohomologie selon les poids positifs.

On aura {\'e}galement recours, comme dans la section \ref{sec7}, aux d{\'e}compositions
suivant les repr{\'e}sentations de \( Sp_{1}Sp_{1} \) ; de ce point de vue, le
fibr{\'e} \( L \) appara{\^\i}t comme le fibr{\'e} homog{\`e}ne associ{\'e} {\`a} la repr{\'e}sentation
\( \ell  \) de \( U_{1} \) avec poids 2 ; rappelons {\'e}galement que l'on a d{\'e}compos{\'e}
\( \mathfrak {sp}_{1}\oplus \mathfrak {sp}_{1}=\mathfrak {u}_{1}\oplus \mathfrak {m} \),
et on notera \( \mathfrak {m}^{1,0}=\mathfrak {m}^{1,0}_{1}\oplus \mathfrak {m}^{1,0}_{2} \)
la d{\'e}composition des formes suivant les deux \( P^{1} \).

\begin{lem}
\label{lem-H1}On a les {\'e}galit{\'e}s 
\begin{eqnarray*}
H^{1}_{\geq 0}(L) & = & \oplus _{k\geq 4}H^{0}(P^{1}_{1},\mathcal{O}(k))\otimes H^{1}(P^{1}_{2},\mathcal{O}(-k+2)),\\
H^{1}_{\geq 0}(\Omega ^{1}_{P^{1}_{1}}\otimes L) & = & H^{1}(P^{1}_{1},\Omega ^{1}_{P^{1}_{1}})\otimes H^{0}(P^{1}_{2},L)\oplus H^{1}_{\geq 0}(\Omega ^{1}_{P^{1}_{1}}\otimes L)_{0},\\
H^{1}_{\geq 0}(\Omega ^{1}_{P^{1}_{1}}\otimes L)_{0} & = & \oplus _{k\geq 4}H^{0}(P^{1}_{1},\mathcal{O}(k-2))\otimes H^{1}(P^{1}_{2},\mathcal{O}(-k+2)),\\
H^{1}_{\geq 0}(\Omega ^{1}_{P^{1}_{2}}\otimes L) & = & \oplus _{k\geq 2}H^{0}(P^{1}_{1},\mathcal{O}(k))\otimes H^{1}(P^{1}_{2},\mathcal{O}(-k)).
\end{eqnarray*}
 Dans la d{\'e}composition suivant les repr{\'e}sentations de \( Sp_{1}Sp_{1} \), la
cohomologie \( H^{1}_{\geq 0}(L) \) est concentr{\'e}e sur les repr{\'e}sentations
\( V^{K,L} \) avec \( K=L+4 \), et est repr{\'e}sent{\'e}e par les \( \left\langle w,\rho (g^{-1})v\right\rangle  \)
avec 
\[
v\otimes w\in V^{K,L}\otimes \hom (V(K,-L),\mathfrak {m}^{0,1}_{2}\otimes \ell );\]
 la cohomologie \( H^{1}_{\geq 0}(\Omega ^{1}\otimes L)_{0} \) est concentr{\'e}e
sur les repr{\'e}sentations \( V^{K,L} \) avec \( K=L+2 \), et est repr{\'e}sent{\'e}e
par les \( \left\langle w,\rho (g^{-1})v\right\rangle  \) avec 
\[
v\otimes w\in V^{K,L}\otimes \hom (V(K,-L),\mathfrak {m}^{1,0}\otimes \mathfrak {m}^{0,1}_{2}\otimes \ell ).\]

\end{lem}
\begin{proof}
Le calcul des \( H^{1} \), d{\'e}compos{\'e}s suivant les poids de l'action de \( S^{1} \),
est similaire {\`a} celui effectu{\'e} en (\ref{cal-H1O}), sachant que \( L=\mathcal{O}(0,2) \).
Le point ici est l'assertion sur les repr{\'e}sentations (on remarquera que les
dimensions des espaces de repr{\'e}sentations indiqu{\'e}s co{\"\i}ncident bien avec les
dimensions des \( H^{1} \)) : le r{\'e}sultat se lit sur les figures, comme dans
la section \ref{sec7} (voir aussi la figure \ref{fig-4} ci-dessous).
\end{proof}
Finissons {\`a} pr{\'e}sent la d{\'e}monstration du lemme \ref{pas-de-H2} : {\'e}tant donn{\'e}
\( \zeta \in \Omega ^{0,1}\otimes \Omega '\otimes L \), tel que \( \overline{\partial }\zeta =0 \),
on cherche \( \phi  \), \( \xi  \) (tel que \( \xi |_{T^{1,0}_{P^{1}_{2}}}=0 \))
et \( \alpha  \) tels que 
\[
\zeta =\overline{\partial }\xi -\partial '(\phi \lrcorner \eta ^{c})+\phi \lrcorner d\eta ^{c}+\alpha \wedge \eta ^{c};\]
 on peut trouver \( \xi  \) et \( \alpha =\overline{\partial }f \) qui tuent
toute la partie de \( \zeta  \) orthogonale {\`a} la cohomologie calcul{\'e}e dans
le lemme \ref{lem-H1}, il nous reste {\`a} utiliser \( \phi  \) et la partie cohomologie
de \( \alpha  \) pour tuer le reste.

Commen{\c c}ons par noter que la composante \( H^{1}(P^{1}_{1},\Omega ^{1}_{P^{1}_{1}})\otimes H^{0}(P^{1}_{2},L) \)
obtenue dans le lemme est repr{\'e}sent{\'e}e par les \( \sigma \omega _{1} \), o{\`u}
\( \sigma \in H^{0}(P^{1}_{2},L) \) ; comme 
\begin{equation}
\label{H1-3}
\overline{\partial }(\sigma \eta )=\sigma (\omega _{1}-\omega _{2}),
\end{equation}
 on voit que quitte {\`a} modifier \( \zeta  \) par \( \sigma \omega _{2} \) (qui
appartient {\`a} l'espace \( H^{1}_{\geq 0}(\Omega ^{1}_{P^{1}_{2}}\otimes L) \)
pour le poids \( k=2 \)) on peut supposer que cette composante s'annule.

Une des cons{\'e}quences importantes du lemme \ref{lem-H1} est que les probl{\`e}mes
de cohomologie pour les composantes de \( \zeta  \) sur \( \Omega ^{0,1}\wedge \eta \otimes L \)
et \( \Omega ^{1,1}\otimes L \) sont compl{\`e}tement d{\'e}connect{\'e}s, puisque les
cohomologies correspondantes apparaissent dans des repr{\'e}sentations diff{\'e}rentes.
On peut donc traiter s{\'e}par{\'e}ment ces deux probl{\`e}mes : ils sont trait{\'e}s dans les
deux lemmes suivants, ce qui ach{\`e}ve la d{\'e}monstration du lemme \ref{pas-de-H2}.\qed

Le premier lemme concerne les repr{\'e}sentations {\`a} \( K=L+4 \).

\begin{lem}
\label{lem-H10}L'image de l'application \( \phi \rightarrow -\partial '(\phi \lrcorner \eta ^{c})+\phi \lrcorner d\eta ^{c} \),
d{\'e}finie sur les \( \phi \in \Omega ^{0,1}_{P^{1}_{2}}\otimes T^{1,0}_{P^{1}_{2}} \)
satisfaisant \( \overline{\partial }_{H}\phi =0 \), contient la cohomologie
\( H^{1}_{\geq 0}(L) \).
\end{lem}
\begin{proof}
On remarquera que les conditions alg{\'e}briques \( \phi \lrcorner F=0 \) et \( \phi \lrcorner F^{L}=0 \)
sont une cons{\'e}quence de \( \phi \in \Omega ^{0,1}_{P^{1}_{2}}\otimes T^{1,0}_{P^{1}_{2}} \).
On a vu dans le lemme \ref{lem-H1} qu'un {\'e}l{\'e}ment \( \zeta ^{0} \) de \( H^{1}(L) \)
peut se repr{\'e}senter par des \( v\otimes w \) avec \( w\in \hom (V(K,-K+4),\mathfrak {m}^{0,1}_{2}\otimes \ell ) \)
; cela nous am{\`e}ne {\`a} chercher un \( \phi  \) v{\'e}rifiant la m{\^e}me propri{\'e}t{\'e} : plus
pr{\'e}cis{\'e}ment, soit \( \phi \in \Omega ^{0,1}_{P^{1}_{2}}\otimes T^{1,0}_{P^{1}_{2}} \),
dans la composante {\`a} poids \( K \) de l'action de \( S^{1} \), alors on calcule
\begin{equation}
\label{for-Rdphi}
R\lrcorner (-\partial '(\phi \lrcorner \eta ^{c})+\phi \lrcorner d\eta ^{c})=-iK\phi \lrcorner \eta ^{c};
\end{equation}
 {\'e}tant donn{\'e} \( \zeta ^{0}\in \Omega ^{0,1}_{P^{1}_{2}}\otimes L \) tel que
\( \overline{\partial }\zeta ^{0}=0 \), on d{\'e}finit \( \phi \in \Omega ^{0,1}_{P^{1}_{2}}\otimes T^{1,0}_{P^{1}_{2}} \)
par 
\[
-iK\phi \lrcorner \eta ^{c}=\zeta ^{0},\]
 alors on a bien \( \overline{\partial }_{H}\phi =0 \), et 
\[
-\partial '(\phi \lrcorner \eta ^{c})+\phi \lrcorner d\eta ^{c}=\zeta ^{0}\eta +\zeta ^{1,1};\]
 la composante \( \zeta ^{1,1} \) sur \( \Omega ^{1,1}\otimes L \) est nulle
en cohomologie, puisque, d'apr{\`e}s le lemme \ref{lem-H1}, celle-ci est concentr{\'e}e
sur les repr{\'e}sentations {\`a} \( K=L+2 \).
\end{proof}
Le second lemme concerne les repr{\'e}sentations {\`a} \( K=L+2 \). 

\begin{lem}
\label{lem-H11}L'image de l'application \( (\phi ,\alpha )\rightarrow -\partial '(\phi \lrcorner \eta ^{c})+\phi \lrcorner d\eta ^{c}+\alpha \eta ^{c} \),
d{\'e}finie sur les couples \( (\phi ,\alpha ) \) tels que \( \phi \lrcorner F=0 \),
\( \phi \lrcorner F^{L}=0, \) \( \overline{\partial }_{H}\phi =0 \), et \( \overline{\partial }\alpha +\phi \lrcorner F^{L}=0 \),
contient \( H^{1}_{\geq 0}(\Omega ^{1}_{P^{1}_{1}}\otimes L)_{0} \) et \( H^{1}_{\geq 0}(\Omega ^{1}_{P^{1}_{2}}\otimes L) \).
\end{lem}
\begin{proof}
Ce lemme est l{\'e}g{\`e}rement plus difficile {\`a} d{\'e}montrer que le pr{\'e}c{\'e}dent. Nous raisonnons
de mani{\`e}re alg{\'e}brique, comme dans la section \ref{sec7}, dans la repr{\'e}sentation
\( V^{K,L} \) avec \( K=L+2 \) ; par abus de notation, nous identifions \( \phi  \)
avec \( \left\langle \phi ,\rho (g^{-1})v\right\rangle  \) o{\`u} \( v\in V^{K,L} \)
et on prend \( \phi \in \hom (V(K,-K+4),\mathfrak {m}^{0,1}\otimes \mathfrak {m}_{1,0}) \)
de la forme (comme dans le lemme pr{\'e}c{\'e}dent) 
\[
\phi =\phi _{\overline{2}}^{2}e^{\overline{2}}e_{2};\]
\begin{figure}[hbt]
{\par\centering \includegraphics{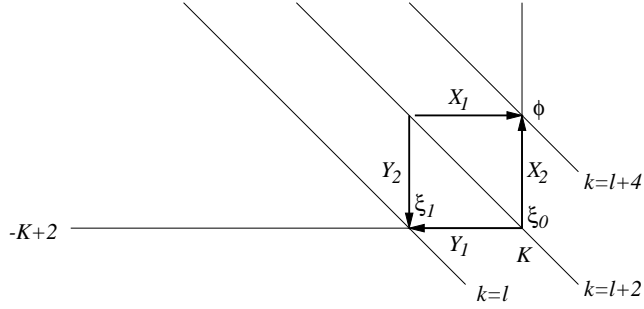} \par}

\caption{\label{fig-4}la cohomologie est concentr{\'e}e en \protect\( (K,-K+2)\protect \)}
\end{figure}
un tel \( \phi  \) v{\'e}rifie automatiquement \( \overline{\partial }_{H}\phi =0 \)
car \( \phi _{\overline{2}}^{2}\circ \rho _{*}(Y_{1})=0 \) si \( \phi _{\overline{2}}^{2} \)
est concentr{\'e} sur \( V(K,-K+4) \), voir la figure \ref{fig-4} ; la forme homog{\`e}ne
\( \eta ^{c} \) s'{\'e}crit alg{\'e}briquement \( \eta ^{c}=(e^{1}-e^{2})\sigma \in \mathfrak {m}^{1,0}\otimes \ell  \),
o{\`u} \( \sigma \in \ell  \) est un vecteur de base fix{\'e} ; par cons{\'e}quent, on
a 
\[
-\phi \lrcorner \eta ^{c}=\phi _{\overline{2}}^{2}e^{\overline{2}}\sigma ;\]
 d'un autre c{\^o}t{\'e}, remarquons que si \( \xi ^{0}\in \hom (V(K,-K+2),\mathbf{C}) \),
alors \( \overline{\partial }(\xi ^{0}\sigma )=\xi ^{0}\circ \rho _{*}(Y_{2})e^{\overline{2}}\sigma  \)
(il n'y a pas de d{\'e}riv{\'e}e sur \( e^{\overline{1}} \) car \( \xi ^{0}\circ \rho _{*}(Y_{1})=0 \)
pour \( \xi ^{0} \) concentr{\'e} sur \( V(K,-K+2) \)) ; compte tenu de la relation
\( -Y_{2}X_{2}=[X_{2},Y_{2}]=H_{2}=-K+4 \) sur \( V(K,-K+4) \), on d{\'e}duit
qu'en posant 
\[
\xi ^{0}=-i\frac{K}{K-4}\phi _{\overline{2}}^{2}\circ \rho _{*}(X_{2}),\]
 on obtient 
\[
R\lrcorner (\overline{\partial }\xi _{0}\sigma +K\phi \lrcorner \eta ^{c})=0,\]
 c'est-{\`a}-dire que, d'apr{\`e}s (\ref{for-Rdphi}), l'image de \( (\phi ,\xi ^{0}\sigma ) \)
est concentr{\'e}e sur \( \Omega ^{1,1}\otimes L \) ; cette image est d'ailleurs
{\'e}gale {\`a} 
\[
-\partial (\phi \lrcorner \eta ^{c})+\xi ^{0}(\omega _{1}-\omega _{2}),\]
 soit alg{\'e}briquement 
\begin{eqnarray*}
 &  & (1-\frac{K}{K-4})\phi _{\overline{2}}^{2}\circ \rho _{*}(X_{2})e^{2}\wedge e^{\overline{2}}\sigma \\
 &  & \qquad +\left( \phi _{\overline{2}}^{2}\circ \rho _{*}(X_{1})e^{1}\wedge e^{\overline{2}}\sigma +\frac{K}{K-4}\phi _{\overline{2}}^{2}\circ \rho _{*}(X_{2})e^{1}\wedge e^{\overline{1}}\sigma \right) ;
\end{eqnarray*}
 le second morceau de cette somme est une (0,1)-forme \( \overline{\partial } \)-ferm{\'e}e
{\`a} valeurs dans \( \Omega ^{1,0}_{P^{1}_{1}}\otimes L \), mais \( H^{1}_{\geq 0}(\Omega ^{1,0}_{P^{1}_{1}}\otimes L) \)
{\'e}tant concentr{\'e} sur \( \hom (V(K,-K+2),\mathfrak {m}^{0,1}\otimes e^{1}) \),
elle est \( \overline{\partial } \)-exacte (ce fait se lit d'ailleurs imm{\'e}diatement
sur la figure \ref{fig-4}), si bien qu'on obtient finalement 
\[
(1-\frac{K}{K-4})\phi _{\overline{2}}^{2}\circ \rho _{*}(X_{2})e^{2}\wedge e^{\overline{2}}\sigma +\overline{\partial }\xi ^{1},\]
 o{\`u} \( \xi ^{1} \) est concentr{\'e} sur \( V(K-2,-K+2) \) ; on en d{\'e}duit que
toute la cohomologie \( H^{1}_{\geq 0}(\Omega ^{1}_{P^{1}_{2}}\otimes L) \)
peut {\^e}tre tu{\'e}e par le choix de \( (\phi ,\xi ) \) comme ci-dessus ; comme on
dispose {\'e}galement de \( \alpha \eta ^{c} \) tel que \( \overline{\partial }\alpha =0 \),
de sorte que 
\[
\alpha \in H_{\geq 0}^{1}(\mathcal{O})=\oplus _{k\geq 2}H^{0}(P^{1}_{1},\mathcal{O}(k))\otimes H^{1}(P^{1}_{2},\mathcal{O}(-k)),\]
 on voit imm{\'e}diatement que tout \( H^{1}_{\geq 0}(\Omega ^{1}_{P^{1}_{1}}\otimes L)_{0}\oplus H^{1}_{\geq 0}(\Omega ^{1}_{P^{1}_{2}}\otimes L) \)
est atteint.
\end{proof}

\section{\label{sec-11}M{\'e}triques quaternion-k{\"a}hl{\'e}riennes}

Nous passons maintenant au cas de la dimension sup{\'e}rieure.

\subsection{Structures de contact quaternioniennes}

Rappelons bri{\`e}vement qu'une structure de contact quaternionienne sur une vari{\'e}t{\'e}
\( X^{4m-1} \) est l'analogue quaternionien d'une structure CR : il s'agit
d'une distribution \( H\subset TX \) de codimension 3, munie d'une structure
quaternionienne conforme (une \( CSp_{m-1}Sp_{1} \)-structure), c'est-{\`a}-dire
d'une m{\'e}trique \( \gamma  \) et de trois structures presque complexes \( I_{1} \),
\( I_{2} \), \( I_{3} \) satisfaisant les relations de commutation des quaternions
(\( I_{1}I_{2}I_{3}=-1 \)) et pr{\'e}servant la m{\'e}trique ; les structures complexes
ne sont d{\'e}finies qu'{\`a} l'action pr{\`e}s de \( Sp_{1} \), et la m{\'e}trique {\`a} un facteur
conforme pr{\`e}s. De plus, on demande la compatibilit{\'e} suivante avec la distribution
\( H \) : localement, il existe une 1-forme \( \eta =(\eta _{1},\eta _{2},\eta _{3}) \)
{\`a} valeurs dans \( \mathbf{R}^{3} \), telle que \( H=\ker \eta  \) et 
\[
d\eta _{i}(X,Y)=\gamma (I_{i}X,Y).\]
 Il est d{\'e}montr{\'e} dans \cite{Biq:00} (voir aussi le survey \cite{Biq:99}) qu'en
dimension \( 4m-1\geq 11 \), une telle donn{\'e}e admet un espace de twisteurs
\( \mathcal{T}^{4m+1} \) qui est une vari{\'e}t{\'e} CR int{\'e}grable munie d'une structure
de contact holomorphe ; en dimension 7, l'existence de l'espace des twisteurs
n'est pas automatique---c'est une condition d'int{\'e}grabilit{\'e} {\`a} imposer {\`a} la structure
(similaire {\`a} la condition d'int{\'e}grabilit{\'e} des structures CR). L'exemple standard
est celui de la sph{\`e}re \( S^{4m-1} \), dont la structure de contact quaternionienne
est fourni par les espaces horizontaux de la fibration de Hopf 
\[
\begin{array}{ccc}
S^{3} & \longrightarrow  & S^{4m-1}\\
 &  & \downarrow \\
 &  & \mathbf{H}P^{m-1}
\end{array}.\]

Le cas de l'espace hyperbolique quaternionien \( \mathbf{H}H^{m} \) est explicite,
et similaire au cas \( \mathbf{H}H^{1}=\mathbf{R}H^{4} \) expliqu{\'e}e dans la
section \ref{sec12} : l'espace des twisteurs de l'espace projectif quaternionien
\( \mathbf{H}P^{m} \) est \( P^{2m+1} \), avec projection twistorielle 
\[
p([z^{1}:z^{2}:\cdots :z^{2m+1}:z^{2m+2}])=[z^{1}+jz^{2}:\cdots :z^{2m+1}+jz^{2m+2}],\]
et l'espace hyperbolique quaternionien se r{\'e}alise comme le domaine 
\[
\mathbf{H}H^{m}=\{|q^{1}|^{2}+\cdots +|q^{m}|^{2}<|q^{m+1}|^{2}\}\subset \mathbf{H}P^{m},\]
 avec pour espace des twisteurs 
\[
\mathcal{N}=\{|z^{1}|^{2}+\cdots +|z^{2m}|^{2}<|z^{2m+1}|^{2}+|z^{2m+2}|^{2}\}\subset P^{2m+1},\]
 muni de la structure r{\'e}elle standard et de la structure de contact holomorphe
\begin{eqnarray*}
\eta ^{c} & = & (z^{1}dz^{2}-z^{2}dz^{1})+\cdots +(z^{2m-1}dz^{2m}-z^{2m}dz^{2m-1})\\
 &  & \phantom {(z^{1}dz^{2}-z^{2}dz^{1})+\cdots }-(z^{2m+1}dz^{2m+2}-z^{2m+2}dz^{2m+1});
\end{eqnarray*}
 le bord \( \mathcal{T}=\partial \mathcal{N} \) appara{\^\i}t {\`a} nouveau comme le
fibr{\'e} \( \mathcal{O}(-1,1) \) sur \( P^{2m-1}\times P^{1} \), avec courbure
\[
iF=\omega _{P^{2m-1}}-\omega _{P^{1}}.\]

\subsection{D{\'e}monstration du th{\'e}or{\`e}me \ref{th-E}}

Le probl{\`e}me d'extension des structures de contact quaternioniennes (en dimension
7, on doit supposer que l'espace des twisteurs existe) est ramen{\'e} par la construction
twistorielle au probl{\`e}me de prolonger dans \( \mathcal{N} \) la d{\'e}formation
de la structure CR et de la structure de contact holomorphe de \( \mathcal{T} \).

Pour r{\'e}aliser ce remplissage, on veut montrer que toute structure CR sur \( \mathcal{T} \),
proche de la structure standard, peut {\^e}tre mise {\`a} coefficients de Fourier positifs
apr{\`e}s action d'un contactomorphisme. Comme dans la section \ref{sec5}, on a
un groupe bien d{\'e}fini de contactomophismes qui agit sur l'espace des structures
CR. On analyse l'action infinit{\'e}simale.

\begin{lem}
\label{lem-111}Pour \( m>1 \), l'image de l'op{\'e}rateur \( \overline{\partial }_{H}\#\overline{\partial } \)
agissant sur les fonctions complexes de classe \( \mathcal{H}^{k+2} \) {\`a} coefficients
de Fourier strictement n{\'e}gatifs est l'espace des \( \phi \in \mathcal{H}^{k}(\Omega ^{0,1}\otimes T^{1,0}) \),
{\`a} coefficients de Fourier strictement n{\'e}gatifs, satisfaisant \( \overline{\partial }_{H}\phi =0 \)
et \( \phi \lrcorner F=0 \).
\end{lem}
\begin{proof}
Soit un entier \( k>0 \), le probl{\`e}me \( \overline{\partial }_{H}\#\overline{\partial }f=\phi  \)
en poids \( -k \) est {\'e}quivalent au probl{\`e}me \( \overline{\partial }\#\overline{\partial }f=\phi  \),
o{\`u} maintenant sur \( P^{2m-1}\times P^{1} \) on voit \( f \) comme une section
de \( \mathcal{O}(-k,k) \) et \( \phi  \) de \( \Omega ^{0,1}\otimes T^{1,0}\otimes \mathcal{O}(-k,k) \).
Or on a l'annulation suivante de la cohomologie, pour \( m>1 \) et \( k>0 \)
: 
\[
H^{1}(P^{2m-1}\times P^{1},T^{1,0}\otimes \mathcal{O}(-k,k))=H^{1}(P^{2m-1}\times P^{1},\mathcal{O}(-k,k))=0,\]
qui se d{\'e}duit imm{\'e}diatement de l'annulation bien connue, pour \( m>1 \),
\[
H^{1}(P^{2m-1},T^{1,0}\otimes \mathcal{O}(\ell ))=0.\]
 D{\'e}duisons-en le lemme : l'annulation du premier \( H^{1} \) indique que l'on
peut r{\'e}soudre le probl{\`e}me \( \overline{\partial }X=\phi  \) ; alors 
\[
\overline{\partial }(\flat X)=\overline{\partial }(X\lrcorner F)=\phi \lrcorner F=0,\]
 donc, par l'annulation du second \( H^{1} \), on peut r{\'e}soudre \( \overline{\partial }f=\flat X \),
donc \( \overline{\partial }\#\overline{\partial }f=\phi  \).
\end{proof}
Quand on se restreint aux fonctions r{\'e}elles, le lemme signifie que le suppl{\'e}mentaire
\( W \) de l'image de \( \overline{\partial }_{H}\#\overline{\partial } \)
d{\'e}fini dans le lemme \ref{lem-64} est constitu{\'e} de tenseurs \( \phi  \) qui
n'ont que des coefficients de Fourier positifs ; comme dans le corollaire \ref{cor-for-W},
une petite d{\'e}formation CR de \( \phi  \) peut {\^e}tre mise en jauge de Coulomb
\( \varphi ^{*}\phi \in W \), donc {\`a} coefficients de Fourier positifs, et par
cons{\'e}quent, comme dans le lemme \ref{lem-fill-J}, est le bord d'une d{\'e}formation
complexe de \( \mathcal{N} \).

Notons qu'on obtient probablement une autre d{\'e}monstration de ce fait en appliquant
le th{\'e}or{\`e}me d'extension de Kiremidjian \cite{Kir:79}, {\`a} condition de montrer
l'annulation de \( H^{2}_{c}(\mathcal{N},T^{1,0}) \). La d{\'e}monstration utilisant
le lemme \ref{lem-111} est plut{\^o}t plus directe, en particulier elle {\'e}vite le
recours au th{\'e}or{\`e}me de Nash-Moser.

La d{\'e}monstration du th{\'e}or{\`e}me \ref{th-fill-etac} reste {\'e}galement valable dans
notre cas, et donne l'extension de la structure de contact holomorphe. Le th{\'e}or{\`e}me
\ref{th-E} en r{\'e}sulte.\qed

\begin{rem}
Le lecteur attentif aura not{\'e} que la m{\'e}thode utilis{\'e}e dans le lemme \ref{lem-111}
donne une voie diff{\'e}rente pour calculer, en dimension 5, les obstructions obtenues
dans la section \ref{sec7}. J'ai pr{\'e}f{\'e}r{\'e} le calcul par les d{\'e}compositions harmoniques,
plus complet : en particulier, les calculs de la section \ref{sec102} sont
imm{\'e}diats gr{\^a}ce {\`a} ces d{\'e}compositions ; de plus, les espaces tangents aux m{\'e}triques
remplissables, d{\'e}termin{\'e}s par les th{\'e}or{\`e}mes \ref{th-gQ-rem} et \ref{th-g-rem},
ne s'expriment simplement qu'{\`a} l'aide de ces m{\^e}mes d{\'e}compositions.
\end{rem}
\bibliographystyle{alpha}
\bibliography{ade}

\begin{thebibliography}{Rum94}

\bibitem[All98]{All:98}
D.~Allcock.
\newblock An isoperimetric inequality for the {H}eisenberg groups.
\newblock {\em Geom. Funct. Anal.}, 8(2):219--233, 1998.

\bibitem[BD91]{Bla-Duc:91}
J.~Bland and T.~Duchamp.
\newblock Moduli for pointed convex domains.
\newblock {\em Invent. Math.}, 104(1):61--112, 1991.

\bibitem[BE90]{Bur-Eps:90}
D.~M. Burns and C.~L. Epstein.
\newblock Embeddability for three-dimensional {C}{R}-manifolds.
\newblock {\em J. Amer. Math. Soc.}, 3(4):809--841, 1990.

\bibitem[Bes87]{Bes:87}
A.~L. Besse.
\newblock {\em Einstein manifolds}.
\newblock Springer-Verlag, Berlin, 1987.

\bibitem[Biq99]{Biq:99}
O.~Biquard.
\newblock Quaternionic contact structures.
\newblock In {\em Proceedings of the Second Meeting on Quaternionic Structures
  in Mathematics and Physics, Rome}, 1999.

\bibitem[Biq00]{Biq:00}
O.~Biquard.
\newblock M{\'e}triques d'{E}instein asymptotiquement sym{\'e}triques.
\newblock {\em Ast{\'e}risque}, 265, 2000.

\bibitem[Bla94]{Bla:94}
J.~Bland.
\newblock Contact geometry and {C}{R} structures on ${S}\sp 3$.
\newblock {\em Acta Math.}, 172(1):1--49, 1994.

\bibitem[CL95]{Che-Lee:95}
J.~H. Ch\^eng and J.~M. Lee.
\newblock A local slice theorem for $3$-dimensional {C}{R} structures.
\newblock {\em Amer. J. Math.}, 117(5):1249--1298, 1995.

\bibitem[Eps92]{Eps:92}
C.~L. Epstein.
\newblock C{R}-structures on three-dimensional circle bundles.
\newblock {\em Invent. Math.}, 109(2):351--403, 1992.

\bibitem[FS74]{Fol-Ste:74}
G.~B. Folland and E.~M. Stein.
\newblock Estimates for the $\bar \partial \sb{b}$ complex and analysis on the
  {H}eisenberg group.
\newblock {\em Comm. Pure Appl. Math.}, 27:429--522, 1974.

\bibitem[Gau86]{Gau:86}
P.~Gauduchon.
\newblock Twisteurs et applications harmoniques en dimension $4$.
\newblock In {\em S\'eminaire de Th\'eorie Spectrale et G\'eom\'etrie, No.\ 4,
  Ann\'ee 1985--1986}, pages 35--93. Univ. Grenoble I, Saint, 1986.

\bibitem[GL91]{Gra-Lee:91}
C.~R. Graham and J.~M. Lee.
\newblock Einstein metrics with prescribed conformal infinity on the ball.
\newblock {\em Adv. Math.}, 87(2):186--225, 1991.

\bibitem[Ham77]{Ham:77}
R.~S. Hamilton.
\newblock Deformation of complex structures on manifolds with boundary. {I}.
  {T}he stable case.
\newblock {\em J. Differential Geometry}, 12(1):1--45, 1977.

\bibitem[Hit95]{Hit:95}
N.~J. Hitchin.
\newblock Twistor spaces, {E}instein metrics and isomonodromic deformations.
\newblock {\em J. Differential Geom.}, 42(1):30--112, 1995.

\bibitem[Hit97]{Hit:97}
N.~J. Hitchin.
\newblock Einstein metrics and the eta-invariant.
\newblock {\em Boll. Un. Mat. Ital. B (7)}, 11(2, suppl.):95--105, 1997.

\bibitem[Kir79]{Kir:79}
G.~K. Kiremidjian.
\newblock A direct extension method for {C}{R} structures.
\newblock {\em Math. Ann.}, 242:1--19, 1979.

\bibitem[LeB82]{LeB:82}
C.~LeBrun.
\newblock {$\mathcal{H}$}-space with a cosmological constant.
\newblock {\em Proc. Roy. Soc. London Ser. A}, 380(1778):171--185, 1982.

\bibitem[LeB84]{LeB:84}
C.~LeBrun.
\newblock Twistor {C}{R} manifolds and three-dimensional conformal geometry.
\newblock {\em Trans. Amer. Math. Soc.}, 284(2):601--616, 1984.

\bibitem[LeB85]{LeB:85}
C.~LeBrun.
\newblock Foliated {C}{R} manifolds.
\newblock {\em J. Differential Geom.}, 22(1):81--96, 1985.

\bibitem[LeB89]{LeB:89}
C.~LeBrun.
\newblock {Quaternionic-K{\"a}hler manifolds and conformal geometry}.
\newblock {\em Math. Ann.}, 284:353--376, 1989.

\bibitem[LeB91]{LeB:91}
C.~LeBrun.
\newblock On complete quaternionic-{K}\"ahler manifolds.
\newblock {\em Duke Math. J.}, 63(3):723--743, 1991.

\bibitem[Lem81]{Lem:81}
L.~Lempert.
\newblock La m{\'e}trique de {K}obayashi et la repr{\'e}sentation des domaines
  sur la boule.
\newblock {\em Bull. Soc. Math. France}, 109:427--474, 1981.

\bibitem[Lem92]{Lem:92}
L.~Lempert.
\newblock On three-dimensional {C}auchy-{R}iemann manifolds.
\newblock {\em J. Amer. Math. Soc.}, 5(4):923--969, 1992.

\bibitem[Ped86]{Ped:86}
H.~Pedersen.
\newblock Einstein metrics, spinning top motions and monopoles.
\newblock {\em Math. Ann.}, 274(1):35--59, 1986.

\bibitem[Rum90]{Rum:90}
M.~Rumin.
\newblock Un complexe de formes diff\'erentielles sur les vari\'et\'es de
  contact.
\newblock {\em C. R. Acad. Sci. Paris S\'er. I Math.}, 310(6):401--404, 1990.

\bibitem[Rum94]{Rum:94}
M.~Rumin.
\newblock Formes diff\'erentielles sur les vari\'et\'es de contact.
\newblock {\em J. Differential Geom.}, 39(2):281--330, 1994.

\bibitem[Sal82]{Sal:82}
S.~Salamon.
\newblock {Quaternionic K{\"a}hler manifolds}.
\newblock {\em Invent. Math.}, 67:143--171, 1982.

\end{thebibliography}

\medskip{}
{\par\raggedleft Institut de Recherche Math{\'e}matique Avanc{\'e}e,\par}

{\par\raggedleft Universit{\'e} Louis Pasteur et CNRS,\par}

{\par\raggedleft 7 rue Ren{\'e} Descartes, F-67084 Strasbourg Cedex\par}

\medskip{}
{\par\raggedleft \texttt{olivier.biquard@math.u-strasbg.fr}\par}

\end{document}